\documentclass[12pt]{amsart}
\DeclareFontFamily{OML}{script}{}
\DeclareFontShape{OML}{script}{m}{it}
{ <5-20> rsfs10 }{}
\DeclareMathAlphabet{\mathscript}{OML}{script}{m}{it}

\renewcommand{\mathcal}[1]{{\mathscript #1}\hspace{0.2ex}}
\usepackage{color}
\ifx\red\undefined
\newcommand{\red}{\color{red}}

\fi
\usepackage[ansinew]{inputenc}
\usepackage[encapsulated]{CJK}

\usepackage{graphicx,graphics}
\usepackage{cite}
\usepackage{pifont}
\usepackage{amsthm}
\usepackage{subfigure}
\usepackage{amscd}
\usepackage{amsmath}
\usepackage{latexsym}
\usepackage{amsfonts}
\usepackage{amssymb}
\usepackage{color}
\usepackage{multicol}
\usepackage{hyperref}
\hypersetup{hypertex=true,
            colorlinks=true,
            linkcolor=blue,
            anchorcolor=blue,
            citecolor=blue}
\usepackage{bm}
\allowdisplaybreaks[4]

\newcommand{\Rmnum}[1]{\uppercase\expandafter{\romannumeral #1}}

\renewcommand{\epsilon}{\varepsilon}

\ifx\text\undefined
\newcommand{\text}{\mbox}
\fi
\ifx\operatorname\undefined
\newcommand{\operatorname}{\mathop}
\fi
\newcommand\be{\begin{equation}}
\newcommand\ee{\end{equation}}
\newcommand\bea{\begin{eqnarray}}
\newcommand\eea{\end{eqnarray}}
\newcommand\beaa{\begin{eqnarray*}}
\newcommand\eeaa{\end{eqnarray*}}

\newenvironment{eqa}{\begin{equation}%
  \begin{array}{rcl}}{\end{array}\end{equation}}

\newcommand\beqa{\begin{eqa}}
\newcommand\eeqa{\end{eqa}}
\numberwithin{equation}{section}

\renewcommand{\tilde}{\widetilde}

\newtheorem{thm}{Theorem}[section]

\newtheorem{lem}[thm]{Lemma}

\newtheorem{rem}{Remark}[section]

\newcommand{\void}[1]{}

\numberwithin{equation}{section}

\begin{document}\begin{CJK}{UTF8}{gkai}
\title[Diffusion]{Equilibrium-diffusion limit of the radiation model}
\author{Lei Li}
\date{\today}
\address[Lei Li]{School of Mathematics and Statistics, Fuyang Normal University, Fuyang, 236037, P. R. China}
\email{202410001@fynu.edu.cn}
\thanks{Corresponding author: Lei Li}

\maketitle
\begin{abstract}
We justify rigorously the equilibrium-diffusion limit of the model consists of a radiative transfer satisfied by the specific intensity of radiation coupled to a diffusion equation satisfied by the material temperature. For general initial data, we construct the existence of the solution to the coupled model in $\mathbb{T}^{3}$ by the Hilbert expansion and prove the convergence of the solutions to the limiting system in the equilibrium-diffusion regime. Moreover, the initial layer for the radiative density and the temperature are constructed to get the strong convergence in $L^\infty$ norm. We also get the convergence rates about the intensity of radiation and temperature in this paper.
\end{abstract}
\section{Introduction}\label{yi}
In this paper, we consider the asymptotic analysis of a coupled model arising in radiative transfer in $\mathbb{T}^{3}$ as
\begin{equation}\label{research equations}\left\{
\begin{split}
&\epsilon^{2}\partial_{t}f^{\epsilon}+\epsilon \vec w\cdot \nabla_{x}f^{\epsilon}+{\epsilon}^{2}(f^{\epsilon}-\overline {f^{\epsilon}})+f^{\epsilon}=B(\theta^{\epsilon})\ \ \mathrm{in} \ \ \mathbb{T}^{3}\times\mathbb{S}^{2},\\
&\epsilon^{2}\partial_{t}\theta^{\epsilon}+\epsilon^{2}\mathrm{div}(\vec{u}\theta^{\epsilon})-\epsilon^{2}\Delta_{x} \theta^{\epsilon}
=\overline {f^{\epsilon}}-B(\theta^{\epsilon})\ \ \mathrm{in} \ \ \mathbb{T}^{3}, \\
&f^{\epsilon}(0, \vec{x}, \vec{w})=h(\vec{x}, \vec{w})\ \ \mathrm{for}\ \ \ ~~  \vec{x}\in\mathbb{T}^{3},\\
&\theta^{\epsilon}(0, \vec{x})=\theta^{0}(\vec{x})\ \ \mathrm{for}\ \  \vec{x}\in\mathbb{T}^{3},
\end{split}\right.
\end{equation}
where $\mathbb{S}^{2}$ is the unit sphere of $\mathbb{R}^{3}$, 
\begin{equation}\label{def average}
\overline {f^{\epsilon}}(t, \vec{x})=\frac{1}{|\mathbb{S}^{2}|}\int_{\mathbb{S}^{d-1}}f^{\epsilon}(t, \vec{x}, \vec{w})\mathrm{d}\vec{w} ,
\end{equation}
the small parameter $\epsilon<1$ represents some small positive parameter derived from physical quantities.
The unknowns are the (nonnegative) specific intensity of radiation $f^{\epsilon}(t, \vec{x}, \vec{w})$, which depends on variables of space $\vec{x}\in\mathbb{T}^{3}$ and direction $\vec{w}\in \mathbb{S}^{2}$, and the material temperature $\theta^{\epsilon}(t, \vec x)$. $f^\epsilon$ satisfies a kinetic equation and $\theta^{\epsilon}$ satisfies a drift-diffusion equation with a given velocity field $\vec{u}$. Furthermore, we assume that $B(\theta^{\epsilon})=(\theta^{\epsilon})^{4}$ under the grey assumption, and $g$ and $\varphi$ are given positive data. Our purpose is to investigate the behavior of solutions $(f^{\epsilon}, \theta^{\epsilon})$ to the Cauchy problem \eqref{research equations} under the appropriate conditions as the positive small parameter $\epsilon$ tends to zero.

The theory of radiation hydrodynamics finds a wide range of application, including such diverse astrophysical phenomena in stellar atmosphere, nonlinear stellar pulsation, supernova explosions, stellar winds, and many others. The general equations of radiation hydrodynamics are a system of a transport equation (Boltzmann equation) coupled with the compressible Euler equations (see \cite{bu-8, bu-9, bu-10}) or Navier-Stokes-Fourier equations (see \cite{bu-10}). Therefore, the study of the mathematical theory of radiation hydrodynamics is very challenging in general. To get a simpler model, Godillon-Lafitte and Goudon \cite{jianmolaiyuanwenzhangb} assumed the density and the velocity of the fluids are given (with a constant density) and considered only a diffusion equation for the material temperature to obtain the system as
\begin{equation}\label{source term equation}\left\{
\begin{split}
&\frac{1}{c}\partial_{t}f+\vec w\cdot \nabla_{x}f=Q(f, \theta)\ \ \mathrm{in} \ \ \mathbb{R}^{+}\times \mathbb{R}^{3}\times \mathbb{S}^{2},\\
&\frac{3}{2}k\Big(\partial_{t}\theta+\mathrm{div}_{x}(\vec{u}\theta)-d\Delta_{x} \theta\Big)=-\int_{S^{2}}\frac{\Lambda}{\gamma}Q(f, \theta)\mathrm{d}\vec{w}\ \ \mathrm{in} \ \ \mathbb{R}^{+}\times \mathbb{R}^{3}, \\
\end{split}\right.
\end{equation}
where $c$ denotes the velocity of light, $\vec{u}=(u_{1}, u_{2}, u_{3}): (t, \vec{x})\in \mathbb{R}^{+}\times \mathbb{R}^{3}\rightarrow \mathbb{R}^{3}$ is a given velocity field, and $d$ and $k$ are the positive diffusion constant and Boltzmann constant, respectively. For simplicity, we use the equivalent notation $\langle f \rangle=\overline{f}$. The source term is
\begin{equation}\nonumber
\begin{aligned}
Q(f, \theta)&=\frac{1}{l_{s}}\Big(\frac{\langle\sigma_{s}\Lambda^{2}f\rangle}{\Lambda^{3}}-\Big\langle\frac{\sigma_{s}}{\Lambda^{2}}\Big\rangle\Lambda f\Big)+\frac{1}{l_{a}}\sigma_{a}\Big(\frac{B(\theta)}{\Lambda^{3}}-\Lambda f\Big)
\\&:=\frac{1}{l_{s}}Q_{s}(f)+\frac{1}{l_{a}}Q_{a}(f, \theta)
\end{aligned}
\end{equation}
where $l_{s}$ is the scattering mean free path, $l_{a}$ is the absorption mean free path, $\sigma_{a}$ and $\sigma_{s}$ represent the absorption and scattering coefficients, respectively. We consider the grey assumption and the black body emission, and then $B(\theta)=\frac{\sigma}{\pi}\theta^{4}$, where $\sigma$ is the Stefan-Boltzmann constant. The weights $\Lambda$ and $\gamma$ related to Doppler corrections are defined by
\begin{equation}\nonumber
\Lambda(t, \vec{x}, \vec{w})=\frac{1-\vec{u}(t, \vec{x})\cdot \vec{w}/c}{\sqrt{{1-|\vec{u}(t, \vec{x})|^{2}/c^{2}}}}
\end{equation}
and
\begin{equation}\nonumber
\gamma(t, \vec{x})=\frac{1}{\sqrt{1-|\vec{u}(t, \vec{x})|^{2}/c^{2}}}.
\end{equation}

Next, we perform a dimensionless analysis of system \eqref{source term equation} as in \cite{jianmolaiyuanwenzhangb} by setting
\begin{equation}\label{equation}\nonumber\left\{
\begin{split}
&t=\overline Tt_{\star}=\frac{L}{u_{\infty}}t_{\star},~~\vec{x}=L\vec{x}_{\star},\\
&\vec{u}(t, \vec{x})=u_{\infty}\vec{u}_{\star}(t_{\star}, \vec{x}_{\star}),~~\theta=\frac{n_{\infty}u_{\infty}^{2}}{3k}\theta_{\star}(t_{\star}, \vec{x}_{\star}), \\
&f(t, \vec{x}, \vec{w})=\sigma\tau_{\infty}^{4}f_{\star}(t_{\star}, \vec{x}_{\star}, \vec{w}),~~B(\theta)=\sigma\tau_{\infty}^{4}B_{\star}(\theta_{\star}),\
\end{split}\right.
\end{equation}
where $L$ denotes characteristic of the flow behavior, and $u_{\infty}, n_{\infty}$ and  $\tau_{\infty}$ denote the reference hydrodynamical velocity, density and temperature respectively. The diffusion coefficient is assumed  to scale to $d=u_{\infty}L\mathcal{D}$ with $\mathcal{D}>0$. The following dimensionless parameters are introduced in \cite{jianmolaiyuanwenzhangb}
\begin{equation}\nonumber
\mathcal{C}=\frac{c}{u_{\infty}},\;\; \mathcal{P}=\frac{2\sigma\tau_{\infty}^{4}}{n_{\infty}u_{\infty}^{3}},\;\; \mathcal{L}_{s}=\frac{L}{l_{s}},\;\; \mathcal{L}_{a}=\frac{L}{l_{a}},\;\; \mathcal{D}=\frac{d}{u_{\infty}L}.
\end{equation}
Then the following dimensionless system is obtained by dropping the stars
\begin{equation}\label{equationdi}\left\{
\begin{split}
&\frac{1}{\mathcal{C}}\partial_{t}f+ \vec{w}\cdot\nabla_{x}f=\mathcal{L}_{s}\Big(\frac{\langle\sigma_{s}\Lambda^{2}f\rangle}{\Lambda^{3}}
-\Big\langle\frac{\sigma_{s}}{\Lambda^{2}}\Big\rangle\Lambda f\Big)+\mathcal{L}_{a}\sigma_{a}\Big(\frac{B(\theta)}{\Lambda^{3}}-\Lambda f\Big),\\
&\partial_{t}\theta+\mathrm{div}_{x}(\vec{u}\theta)-\mathcal{D}\Delta_{x}\theta=-\mathcal{P}\mathcal{L}_{s}
\Big\langle\frac{\Lambda}{\gamma}\Big(\frac{\langle\sigma_{s}\Lambda^{2}f\rangle}{\Lambda^{3}}-\Big\langle\frac{\sigma_{s}}{\Lambda^{2}}\Big\rangle\Lambda f\Big)\Big\rangle\\&\hspace{4.7cm}-\mathcal{P}\mathcal{L}_{a}\Big\langle\frac{\Lambda}{\gamma}\sigma_{a}\Big(\frac{B(\theta)}{\Lambda^{3}}-\Lambda f\Big)\Big\rangle,\
\end{split}\right.
\end{equation}
where
\begin{equation}\nonumber
\Lambda(t, \vec{x}, \vec{w})=\gamma(t, \vec{x})\Big(1-\frac{\vec{w}\cdot \vec{u}(t, \vec{x})}{\mathcal{C}}\Big), \;\;\;\;\gamma(t, \vec{x})=\frac{1}{\sqrt{1-\frac{|\vec{u}(t, \vec{x})|^{2}}{\mathcal{C}^{2}}}}.
\end{equation}

For the nonrelativistic flows with a moderate amount of radiation in the flow with the following scalings,
\begin{equation}\nonumber
\mathcal{C}=\frac{1}{\epsilon},\;\;\;\;\frac{\mathcal{P}}{\mathcal{C}}=1,\;\;\;\;\mathcal{D}=1.
\end{equation}
Godillon-Lafitte and Goudon \cite{jianmolaiyuanwenzhangb} identified both the equilibrium and nonequilibrium regimes for system \eqref{equationdi} with Doppler corrections in the whole space. The equilibrium regime corresponds to
\begin{equation}\label{equ1}
\mathcal{L}_{s}=\epsilon,\;\;\;\; \mathcal{L}_{a}=\frac{1}{\epsilon}
\end{equation}
 and the nonequilibrium regime corresponds to
\begin{equation}\label{nonequ1}
\mathcal{L}_{s}=\frac{1}{\epsilon},\;\;\;\; \mathcal{L}_{a}=\epsilon.
\end{equation}

In the first case, they proved that the limiting system is derived only by the material temperature satisfying a nonlinear drift-diffusion equation, and in the latter case, the radiation temperature and the material temperature are coupled by the nonlinear drift-diffusion system.

We introduce some results about the single neutron transport equation in the following. In \cite{shijiu}, Bardos, Golse, Perthame and Sentis considered the existence and the diffusive limit of the following equation
\begin{equation}\label{single-neutron-c1}\left\{
\begin{split}
&\frac{\partial f^{\epsilon}}{\partial t}+\frac{\vec w\cdot\nabla_{x} f^{\epsilon}}{\epsilon}+\frac{\sigma(\overline {f^{\epsilon}})}{\epsilon^{2}}(f^\epsilon -\overline {f^{\epsilon}})=0, (t, \vec x, \vec w)\in\mathbb{R}^{+}\times\Omega\times\mathbb{S}^{N},\\
&f^{\epsilon}(t, \vec x_{0}, \vec w)=k, (t, \vec x_{0}, \vec w)\in\mathbb{R}^{+}\times\partial\Omega\times\mathbb{S}^{N}, \vec w\cdot\vec n<0,\\
&f^{\epsilon}(0, \vec x, \vec w)=h(\vec x), (\vec x, \vec w)\in\Omega\times\mathbb{S}^{N}.\\
\end{split}\right.
\end{equation}
Here, $f^{\epsilon}=f^{\epsilon}(t, \vec x, \vec w)$, $\sigma$ is a function of $\overline {f^{\epsilon}}$, $\Omega\subset\mathbb{R}^{N+1}$, $\mathbb{S}^{N}$ is the unit sphere of $\mathbb{R}^{N+1}$, $k$ is a constant and $h=h(\vec x)$ only depend on $\vec x$ and is independent of $\vec w$. The existence of \eqref{single-neutron-c1} follows from the study of the spectrum of the transport operator, which has been carried out by many authors \cite{ershi, ershiyi, ershier, ershisan, ershisi}. Furthermore, the diffusion limit of \eqref{single-neutron-c1} has been studied extensively \cite{lionssidawenzhang, bu-8, yi-bu1, yi-bu2, ershiyi, ershisan}. 

There have been extensive research efforts about the boundary layer theory with geometric correction; see \cite{convex-domain-123, 3-d-diffusive-limit-convex-domain, unsteady-neutron-transport-2-d-unit-disk, in-flow-boundry-neutron-transport, 3-d-zuinanqingxing, yuanhuanneutrontransport, geometric-CMP}.

For the equilibrium-regime coupled model, we refer to \cite{jianmolaiyuanwenzhangb} again. However, the mathematical results are not complete. Because Godillon-Lafitte and Goudon \cite{jianmolaiyuanwenzhangb} didn't give the proof process of the well-posedness of the equilibrium-regime coupled model in $\mathbb{R}^{3}$. They only focus on the diffusive limit analysis and the numerical calculation. In our paper, we will focus on both the well-posedness and the diffusion limit of the equilibrium-regime coupled model. However, our model is simpler without Doppler corrections in contrast to \cite{jianmolaiyuanwenzhangb}. Furthermore, we consider the coupled model in $\mathbb{T}^{3}$ with general initial data, which will be more challenging. Recently, there are many results about the following kind of equilibrium-regime coupled model:
\begin{equation}\label{research equations-dayu0-c1-equilibrium-jia}\left\{
\begin{split}
&{\epsilon}^{2}\partial_{t}f^{\epsilon}+\epsilon \vec w\cdot \nabla_{x}f^{\epsilon}+f^{\epsilon}=B(\theta^{\epsilon}),\ \ (t, \vec x, \vec w)\in \mathbb{R}^+\times \Omega\times \mathbb{S}^{d-1},\\
&\epsilon^{2}\partial_{t}\theta^{\epsilon}-\epsilon^{2}\mathcal{D}\Delta_{x} \theta^{\epsilon}
=\overline {f^{\epsilon}}-B(\theta^{\epsilon}),\ \ \ \ (t, \vec x)\in \mathbb{R}^+\times \Omega, \\
\end{split}\right.
\end{equation}
where $\mathcal{D}>0$ and $\Omega\subset\mathbb{R}^{d}$. We note that \eqref{research equations-dayu0-c1-equilibrium-jia} is also a simpler model without Doppler corrections in contrast to the corresponding model in \cite{jianmolaiyuanwenzhangb}. Klar and Schmeiser \cite{bu-1} showed the existence and diffusive limit of the similar model of (\ref{research equations-dayu0-c1-equilibrium-jia}). Moreover, they investigated the nonlinear Milne problem describing the boundary layer and proved the existence result. Finally, they presented numerical results for different physical situations. It is noted that Ghattassi, Huo and Masmoudi \cite{traceth} proved the global existence of weak solutions for the similar model of (\ref{research equations-dayu0-c1-equilibrium-jia}) and the convergence of the weak solutions to a nonlinear diffusion model under the diffusive limit, which extended the work by Klar and Schmeiser \cite{bu-1}. Other works on similar models can be found in \cite{bu-2, bu-3, bu-4, bu-5}. More recently, Ghattassi, Huo and Masmoudi \cite{traceth-1} study the diffusive limit of the steady state of (\ref{research equations-dayu0-c1-equilibrium-jia}) for non-homogeneous Dirichlet boundary conditions in a bounded domain with flat boundaries. Considering the arising of the boundary layers, a composite approximate solution is constructed using asymptotic analysis by resorting to the results of \cite{traceth-2} on the boundary layer problem to the steady state of (\ref{research equations-dayu0-c1-equilibrium-jia}). This result \cite{traceth-1} extends their previous work \cite{traceth} for the well-prepared boundary data case to the general boundary data. Furthermore, Ghattassi, Huo and Masmoudi \cite{traceth-3} study the diffusive limit of the steady state of (\ref{research equations-dayu0-c1-equilibrium-jia}) in a unit disk with non-flat boundaries by considering the geometric corrections \cite{geometric-CMP}. Note that the models (\ref{research equations-dayu0-c1-equilibrium-jia}) considered in \cite{traceth-1, traceth-2, traceth-3} are the steady case of equilibrium regime. The main purpose of the present paper is to identify the unsteady case of equilibrium regime for the general initial data in $\mathbb{T}^{3}$ which is the difference between our work and \cite{traceth-1, traceth-2, traceth-3}. Similar models has been considered in \cite{traceth} for the well-prepared data, but we focus on the general initial data and overcome the difficulties arising from the initial layers. We also note that the system \eqref{research equations} has extra terms ${\epsilon}^{2}(f^{\epsilon}-\overline {f^{\epsilon}})$ and $\epsilon^{2}\mathrm{div}(\vec{u}\theta^{\epsilon})$ in contrast to \eqref{research equations-dayu0-c1-equilibrium-jia} which has many physical meanings and takes many mathematical calculation difficulties for us. 

For the non-equilibrium regime of \eqref{equationdi}, Ju, Li and Zhang have some works in \cite{yuanmoxing-bu, yuanmoxing-bu-jia, yuanmoxing-bu-jia1, yuanmoxing-bu-jia2}. In \cite{yuanmoxing-bu}, they get the non-equilibrium regime of \eqref{equationdi} for the well-prepared data in a unit disk. In \cite{yuanmoxing-bu-jia}, Li and Zhang get the non-equilibrium regime of \eqref{equationdi} for the general initial data in $\mathbb{T}^{3}$. Then, Li \cite{yuanmoxing-bu-jia1} extends the result in \cite{yuanmoxing-bu} to the general initial and boundary data by the Hilbert expansion. In \cite{yuanmoxing-bu-jia2}, Li and Zhang get the diffusion limit for the steady case of the model in \cite{yuanmoxing-bu}. 

We mention that one can take the P1 hypothesis to get the compressible Euler/Navier-Stokes-P1 approximation model when the distribution of photons is almost isotropic (see \cite{bu-10}). Recently, there are some investigations on the singular limits of such approximation model. For example, Jiang, Li and Xie \cite{bu-7} studied the nonrelativistic limit problem for the Navier-Stokes-Fourier-P1 approximation model as the reciprocal of the light speed tends to zero. Jiang, Ju and Liao \cite{bu-6} showed the diffusive limit of the compressible Euler-P1 approximation model arising in radiation hydrodynamics as the Mach number tends to zero. Furthermore, they also considered \cite{bu-6-jia0} the diffusive limit of the compressible Euler-P1 approximation model arising in radiation hydrodynamics for the fixed Mach number case. Note that Ju and Liao \cite{bu-6-jia1} considered the equilibrium and non-equilibrium diffusion limits for a baratropic model of radiative flow, recently. In \cite{jia-jia-4}, Li and Zhang extend the results \cite{bu-6} to the MHD-P1 model. In \cite{jia-jia-4-jia}, Li and Zhang showed the diffusive limit of the compressible Navier-Stokes-Fourier-P1 approximation model arising in radiation hydrodynamics as the Mach number tends to zero with partial general data. We refer to \cite{DN1, bujia-1, bujia-3, bujia-4, bujia-5, bujia-6, bujia-7} for more results of radiation models at low Mach number. 

The rigorous proof of both equilibrium and non-equilibrium diffusion limits for the NSF-R model in the framework of weak solutions has been given by Ducomet and Ne\v{c}asov\'{a} \cite{bujia-4, bujia-10} by using the relative entropy method for well-prepared initial data. Danchin and Ducomet \cite{bujia-2} established the existence of global strong solutions to the isentropic Navier-Stokes-P1 model with small enough initial data in critical regularity spaces. However, to our best knowledge, the rigorous verification of the diffusion limit of NSF-R model is still unsolved. Recently , Ju, Li and Zhang \cite{yuanmoxing-bu-jia3jia} justify rigorously the non-equilibrium-diffusion limit of the compressible Euler model coupled with a radiative transfer equation arising
in radiation hydrodynamics with non-relativistic source terms. Furthermore, Li and Zhang get the convergence rates in \cite{yuanmoxing-bu-jia3}. Then, Li and Zhang extend the results in \cite{yuanmoxing-bu-jia3jia, yuanmoxing-bu-jia3} to the MHD-R model \cite{yuanmoxing-bu-jia4}.

The main purpose of the present paper is to identify the equilibrim asymptotic regime of system \eqref{equationdi} in $\mathbb{T}^{3}$. In the present paper, we will ignore the Doppler effects by assuming that $\gamma=\Lambda\equiv 1$ for the nonrelativistic flow. For simplicity, we will set $\sigma_{s}=\sigma_{a}=1.$ Moreover, we ignore the influence of other constants and we set them to be one. As a result, it is equivalent to proving the limit of solution $(f^{\epsilon}, \theta^{\epsilon})$ to the Cauchy problem \eqref{research equations} as $\epsilon$ tends to zero. 

\begin{thm}\label{result}
Assume that $\vec u\in C^{N+2}([0, \infty)\times\mathbb{T}^{3})$, $h(\vec{x}, \vec{w})\in L^{\infty}(\mathbb{S}^{2}; H^{N+2}(\mathbb{T}^{3}))$ and $\theta^{0}(\vec{x})\in H^{N+2}(\mathbb{T}^{3})$ with $h(\vec{x}, \vec{w})>0, \theta^{0}(\vec x)\geq a$ and $\|h-(\theta^{0})^{4}\|_{L^{\infty}(\mathbb{S}^{2}; H^{N+2}(\mathbb{T}^{3}))}\leq \eta$, where $N\geq 9$, $a$ and $\eta$ are positive constants and $\eta$ is sufficiently small.  There exists a $T>0$ independent of $\epsilon$, such that the problem (\ref{research equations}) has a unique solution $(f^{\epsilon}, \theta^{\epsilon})\in C([0, T]; L^{2}(\mathbb{S}^{2}; H^{2}(\mathbb{T}^{3})))\times C([0, T]; H^{2}(\mathbb{T}^{3}))\cap L^{2}(0, T; H^{3}(\mathbb{T}^{3}))$. Furthermore, the unique solution $(f^{\epsilon}, \theta^{\epsilon})$ satisfies
\begin{equation}\label{jielunbudengshi 1}
\|f^{\epsilon}-f_{0}-f_{I, 0}\|_{L^{\infty}((0, T)\times\mathbb{T}^{3}\times\mathbb{S}^{2})}=O(\epsilon),
\end{equation} and
\begin{equation}\label{jielunbudengshi 2}
\|\theta^{\epsilon}-\theta_{0}-\theta_{I, 0}\|_{C^{0}([0, T]; H^{2}(\mathbb{T}^{3}))}=O(\epsilon),
\end{equation}
where $(f_{0}, \theta_{0})\in (C([0, T]; H^{N+2}(\mathbb{T}^{3}))\cap L^{2}(0, T; H^{N+3}(\mathbb{T}^{3})))\times (C([0, T]; H^{N+2}(\mathbb{T}^{3}))\cap L^{2}(0, T; H^{N+3}(\mathbb{T}^{3})))$ satisfies
\begin{equation}\label{equations about f0 epsilon and theta}\left\{
\begin{split}
&f_{0}=\overline {f_{0}}=\theta_{0}^{4},\\
&\partial_{t}(\theta_{0}+\theta_{0}^{4})+\mathrm{div}(\vec{u}\theta_{0})-\Delta_{x}\Big(\theta_{0}+\frac{1}{3}\theta_{0}^{4}\Big)=0,~~\mathrm{in} \ \ (0,T)\times\mathbb{T}^{3},\\
&\theta_{0}(0, \vec{x})=\theta_{0}^{0}(\vec{x}),\;\;\overline {f_{0}}(0, \vec{x})=(\theta_{0}^{0})^{4}~~\mathrm{in}  \ \ \mathbb{T}^{3},
\end{split}\right.
\end{equation}
and the zeroth-order initial layer $(f_{I, 0}, \theta_{I, 0})\in C^{1}([0, \infty); L^{\infty}(\mathbb{S}^{2}; H^{N+2}(\mathbb{T}^{3})))\times C^{1}([0, \infty); \\H^{N+2}(\mathbb{T}^{3}))$ is defined as
\begin{equation}\label{inicoup}\left\{
\begin{split}
&\frac{\partial f_{I, 0}}{\partial\tau}+f_{I, 0}=B(\theta_{0}^{0}+\theta_{I, 0})-B(\theta_{0}^{0}),\\
&\frac{\partial \theta_{I, 0}}{\partial\tau}=(\overline{f_{I, 0}}-B(\theta_{0}^{0}+\theta_{I, 0})+B(\theta_{0}^{0})),\\
&f_{I, 0}(0, \vec{x}, \vec{w})=h(\vec{x}, \vec{w})-f_{0}(0, \vec x),\\
&\theta_{I, 0}(0, \vec{x}, \vec{w})=\theta^{0}(\vec x)-\theta_{0}^{0}(\vec x),\\
&\lim_{\tau\rightarrow\infty}f_{I, 0}(\tau, \vec{x}, \vec{w})=0,\ \ \lim_{\tau\rightarrow\infty}\theta_{I, 0}(\tau, \vec{x}, \vec{w})=0,
\end{split}\right.
\end{equation}
where
$\tau=\frac{t}{\epsilon^{2}}$.
\end{thm}

\begin{rem}
We can verify that $f_{0}(0, \vec x)=(\theta_{0}^{0})^{4}(\vec x), \theta_{0}^{0}(\vec x)\in H^{N+2}(\mathbb{T}^{3})$ and $\theta_{0}^{0}\geq b>0$ for some positive constant $b$ depending on $a$ through studying the initial layer coupled model \eqref{inicoup}.
\end{rem}

\begin{rem}
The close to equilibrium state of the initial data $\|h-(\theta^{0})^{4}\|_{L^{\infty}(\mathbb{S}^{2}; H^{N+2}(\mathbb{T}^{3}))}\\\leq \eta$, for small $\eta$, which will help us to establish the global existence for initial layers $(f_{I, k}, \theta_{I, k})$, $0\leq k\leq N$. Furthermore, the small estimates \eqref{ineq-theta00-jia1} and \eqref{ineq-theta00-jia2} can help us to get the estimates \eqref{bd11}, \eqref{bd13}, \eqref{bd14} and \eqref{bd15} contained the $O(\eta)$ terms which will be absorbed by the left term. So, the estimates \eqref{bd00} can be established by applying the Gronwall's inequality. 
\end{rem}

\begin{rem}
One of the main steps of getting the convergence rates \eqref{jielunbudengshi 1} and \eqref{jielunbudengshi 2} is to get the uniform estimate \eqref{bd00}. To get the error estimate \eqref{bd00}, we need to show the uniform estimates \eqref{bd0-jia5}, \eqref{bd0-jia7} and \eqref{bd0-jia8} for $\alpha=0,$ $\alpha=1$ and $\alpha=2$ iteratively to overcome the difficulties arising from the bad term $C\int_{0}^{t}\int_{\mathbb{T}^{3}}\int_{\mathbb{S}^{2}}(|D_{x}^{\alpha-1}\theta_{r}|^{2}
+|D_{x}^{\alpha-2}\theta_{r}|^{2})
\mathrm{d}\vec w\mathrm{d}\vec x\mathrm{d}s$ in \eqref{bd0-jia1}.
\end{rem}

We give the outline of the proof of Theorem \ref{result}. We first expand $(f^{\epsilon}, \theta^{\epsilon})$ as $f^{\epsilon}=\sum_{k=0}^{N}\epsilon^{k}f_{k}+\sum_{k=0}^{N}\epsilon^{k}f_{I, k}+f_{r}$ and $\theta^{\epsilon}=\sum_{k=0}^{N}\epsilon^{k}\theta_{k}+\sum_{k=0}^{N}\epsilon^{k}\theta_{I, k}+\theta_{r}$ where $f_{k}, f_{I, k}$ and $f_{r}$ are the interior expansion, initial layers expansion and remainder part of $f^\epsilon$, and $\theta_{k}, \theta_{I, k}$ and $\theta_{r}$ are the interior expansion, initial layers expansion and remainder part of $\theta^\epsilon$. In the following, we prove the well-posedness of the above expansion of interior, initial layer and boundary layer parts, and get the estimates about $\epsilon$ by careful calculations. Then, we can get the well-posedness of \eqref{research equations} by adding the interior, initial layer parts and the remainder part together. The main step of getting the convergence rates of $(f^{\epsilon}, \theta^{\epsilon})$ is estimating the remainders $f_{r}$ and $\theta_{r}$ about $\epsilon$.

Throughout this paper, $C$ denotes a certain positive constant and $C(\cdot)$ and $C_{1}(\cdot)$ are the positive constants depending on the quantity $"\cdot"$. $D_{x}$ and $D_{t}$ represent the operator $\partial_{x_{i}}, i=1, 2, 3$ and $\partial_{t}$, respectively. $H^{k}(\mathbb{T}^{3})=W^{k, 2}(\mathbb{T}^{3})$ denotes the usual Lebesgue space on $\mathbb{T}^{3}$ with norm $\|\cdot\|_{H^{k}(\mathbb{T}^{3})}$. We denote by $L^{p}(0, T; H^{k}(\mathbb{T}^{3}))$ $(1\leq p\leq\infty)$ the space of $L^{p}$ functions on $(0, T)$ with values in $H^{k}(\mathbb{T}^{3})$, and $C^{0}(I, H^{k}(\mathbb{T}^{3}))$ standards for the space of continuous functions on the interval $I$ with values in $H^{k}(\mathbb{T}^{3})$. For simplicity, we use $L_{T}^{2}, L_{\vec w}^{2}, H_{\vec x}^{k}$ to denote $L^{2}(0, T), L^{2}(\mathbb{S}^{2}), H^{k}(\mathbb{T}^{3})$, respectively. For $k\geq 0$ and $k\in\mathbb{Z}$, $C^{k}([0, \infty)\times\mathbb{T}^{3})$ denotes the set satisfying $D_{x}^{l}D_{t}^{r}f\in C^{0}([0, \infty)\times\mathbb{T}^{3})$ where $0\leq l+r\leq k$.

The paper is organized as follows. In Section \ref{er}, we discuss the well-posedness of neutron transport equation and prove the first part of Theorem \ref{result}. In Section \ref{san}, we give the Hilbert expansion of \eqref{research equations}. In Section \ref{si}, we will give the proof of Theorem \ref{result}.

\section{Well-posedness of equilibrium regime neutron transport equation}\label{er}
In this section, we consider the well-posedness of equilibrium regime neutron transport equation, which can be written as
\begin{equation}\label{research equation}\left\{
\begin{split}
&{\epsilon}^{2}\frac{\partial I}{\partial t}+\epsilon \vec{w}\cdot \nabla_{x}I+{\epsilon}^{2}(I-\overline I)+I=F(t, \vec{x}, \vec{w})~~\mathrm{in} \ \ (0, \infty)\times\mathbb{T}^{3}\times \mathbb{S}^{2},\\
&I(0, \vec{x}, \vec{w})=h(\vec{x}, \vec{w})~~\mathrm{in} \ \ \mathbb{T}^{3}\times \mathbb{S}^{2}.\\
\end{split}\right.
\end{equation}

Define the $L^{2}$ and $L^{\infty}$ norms in $\mathbb{T}^{3}\times \mathbb{S}^{2}$ by
\begin{equation}\nonumber
\|I\|_{L^{2}(\mathbb{T}^{3}\times \mathbb{S}^{2})}=\Big(\int_{\mathbb{T}^{3}}\int_{\mathbb{S}^{2}}|I(\vec{x}, \vec{w})|^{2}\mathrm{d}\vec{w}\mathrm{d}\vec{x}\Big)^{\frac{1}{2}},
\end{equation}

\begin{equation}\nonumber
\|I\|_{L^{\infty}(\mathbb{T}^{3}\times \mathbb{S}^{2})}=\sup_{(\vec{x}, \vec{w})\in\mathbb{T}^{3}\times \mathbb{S}^{2}}|I(\vec{x}, \vec{w})|.\;\;
\end{equation}
Similar notations also applies to the space $[0, T]\times\mathbb{T}^{3}\times\mathbb{S}^{2}.$

\subsection{$L^{\infty}$ well-posedness}
We will show the existence and $L^{\infty}$ estimate of the transport equation (\ref{research equation}) in the following.
\vskip 2mm

\begin{lem}\label{main result}
Assume $F(t, \vec{x}, \vec{w})\in L^{\infty}((0, \infty)\times\mathbb{T}^{3}\times \mathbb{S}^{2}), h(\vec{x}, \vec{w})\in L^{\infty}(\mathbb{T}^{3}\times \mathbb{S}^{2})$ and $F, h\geq 0$. Then there exists a unique nonnegative solution $I(t, \vec{x}, \vec{w})\in L^{\infty}((0, \infty)\times\mathbb{T}^{3}\times \mathbb{S}^{2})$ of the transport equation (\ref{research equation}) which satisfies
\begin{equation}\label{estimate result}
\|I\|_{L^{\infty}((0, \infty)\times\mathbb{T}^{3}\times \mathbb{S}^{2})}\leq C(\|h\|_{L^{\infty}(\mathbb{T}^{3}\times \mathbb{S}^{2})}+\|F\|_{L^{\infty}((0, \infty)\times \mathbb{T}^{3}\times \mathbb{S}^{2})}).
\end{equation}

\end{lem}
\noindent {\bf Proof.} We construct an approximating sequence $\{I_{k}\}_{k=0}^{\infty}$ which are also periodic, where $I^{0}=0$ and

\begin{equation}\label{approximating sequence equation}\left\{
\begin{split}
&{\epsilon}^{2}\frac{\partial I^{k}}{\partial t}+\epsilon \vec w\cdot \nabla_{x}I^{k}+(1+\epsilon^{2})I^{k}
-\epsilon^{2}\overline {I^{k-1}}=F(t, \vec{x}, \vec{w})~~\mathrm{in} \ \ (0, \infty)\times\mathbb{T}^{3}\times \mathbb{S}^{2},\\
&I^{k}(0, \vec{x}, \vec{w})=h(\vec{x}, \vec{w})~~\mathrm{in} \ \ \mathbb{T}^{3}\times \mathbb{S}^{2}.\\
\end{split}\right.
\end{equation}

The characteristics $(T(s), X(s), W(s))$ of equation (\ref{approximating sequence equation}) which goes through $(t, \vec{x}, \vec{w})$ is defined by
\begin{equation}\label{equation}\left\{
\begin{split}
&(T(0), X(0), W(0))=(t, \vec{x}, \vec{w}),\\
&\frac{dT(s)}{ds}=\epsilon^{2},\\
&\frac{dX(s)}{ds}=\epsilon W(s),\\
&\frac{dW(s)}{ds}=0,\\
\end{split}\right.
\end{equation}
which implies
\begin{equation}\nonumber\label{equation}\left\{
\begin{split}
&T(s)=t+\epsilon^{2}s,\\
&X(s)=\vec{x}+(\epsilon \vec{w})s,\\
&W(s)=\vec{w}.\\
\end{split}\right.
\end{equation}
Then,
\begin{equation}\label{similar-1}
\begin{split}
I^{k}(t, \vec{x}, \vec{w})=&h(\vec{x}-\frac{(\epsilon t\vec{w})}{\epsilon^{2}}, \vec{w})e^{-(1+\epsilon^{2})\frac{t}{\epsilon^{2}}} \\&+\int_{0}^{\frac{t}{\epsilon^{2}}}(\epsilon^{2}\overline {I^{k-1}}+F)(\epsilon^{2}s, \vec{x}-\epsilon(\frac{t}{\epsilon^{2}}-s)\vec{w}, \vec{w})e^{-(1+\epsilon^{2})(\frac{t}{\epsilon^{2}}-s)}ds.
\end{split}
\end{equation}
Since $I^{0}=0$ and $h, g, F\geq 0$, we obtain $I^{1}\geq 0$. Then, we get $I^{k}\geq 0$ for any $k\geq 0$ iteratively. Set $v_{k}=I^{k}-I^{k-1}$ for $k\geq 1.$ Then, we have
\begin{equation}\nonumber
\begin{split}
v_{k+1}(t, \vec{x}, \vec{w})=\int_{0}^{\frac{t}{\epsilon^{2}}}\epsilon^{2}\overline {v_{k}}(\epsilon^{2}s, \vec{x}-\epsilon(\frac{t}{\epsilon^{2}}-s)\vec{w}, \vec{w})e^{-(1+\epsilon^{2})(\frac{t}{\epsilon^{2}}-s)}ds.
\end{split}
\end{equation}
Since $\|\overline {v_{k}}\|_{L^{\infty}((0, \infty)\times\mathbb{T}^{3}\times \mathbb{S}^{2})}\leq \|v_{k}\|_{L^{\infty}((0, \infty)\times\mathbb{T}^{3}\times \mathbb{S}^{2})}$, we obtain
\begin{equation}\nonumber
\begin{split}
\|v_{k+1}\|_{L^{\infty}((0, \infty)\times\mathbb{T}^{3}\times \mathbb{S}^{2})}\leq\frac{1-e^{-(1+\epsilon^{2})\frac{t}{\epsilon^{2}}}}{1+\epsilon^{2}}\|\epsilon^{2}v_{k}\|_{L^{\infty}((0, \infty)\times\mathbb{T}^{3}\times \mathbb{S}^{2})}.
\end{split}
\end{equation}
It follows that
\begin{equation}\nonumber
\|v_{k+1}\|_{L^{\infty}((0, \infty)\times\mathbb{T}^{3}\times \mathbb{S}^{2})}\leq\frac{1}{1+\epsilon^2}\|\epsilon^{2}v_{k}\|_{L^{\infty}((0, \infty)\times\mathbb{T}^{3}\times \mathbb{S}^{2})},
\end{equation}
which implies $\{v_{k}\}$ is a contraction sequence. According to the fact that $v_{1}=I^{1}$, we have
\begin{equation}\nonumber
\|v_{k+1}\|_{L^{\infty}((0, \infty)\times\mathbb{T}^{3}\times \mathbb{S}^{2})}\leq(\frac{1}{1+\epsilon^2})^{k}\|\epsilon^{2}I^{1}\|_{L^{\infty}((0, \infty)\times\mathbb{T}^{3}\times\mathbb{S}^{2})}
\end{equation}
for $k\geq 1.$ Therefore, $I^{k}$ converges strongly in $L^{\infty}((0, \infty)\times\mathbb{T}^{3}\times \mathbb{S}^{2})$ to a limit solution $I$ which satisfies
\begin{equation}\label{infinity estimate about limit solution}
\|I\|_{L^{\infty}((0, \infty)\times\mathbb{T}^{3}\times \mathbb{S}^{2})}\leq\sum_{k=1}^{\infty}\|\epsilon^{2}v_{k}\|_{L^{\infty}((0, \infty)\times\mathbb{T}^{3}\times \mathbb{S}^{2})}\leq (1+\epsilon^{2})\|I^{1}\|_{L^{\infty}((0, \infty)\times\mathbb{T}^{3}\times \mathbb{S}^{2})}.
\end{equation}
Since
\begin{equation}\nonumber
\begin{split}
I^{1}(t, \vec{x}, \vec{w})=h(\vec{x}-\frac{(\epsilon t\vec{w})}{\epsilon^{2}}, \vec{w})e^{-(1+\epsilon^{2})\frac{t}{\epsilon^{2}}} +\int_{0}^{\frac{t}{\epsilon^{2}}}F(\epsilon^{2}s, \vec{x}-\epsilon(\frac{t}{\epsilon^{2}}-s)\vec{w}, \vec{w})e^{-(1+\epsilon^{2})(\frac{t}{\epsilon^{2}}-s)}ds,
\end{split}
\end{equation}
we obtain
\begin{equation}\label{the first term about the sequence infinity estimate}
\|I^{1}\|_{L^{\infty}((0, \infty)\times\mathbb{T}^{3}\times \mathbb{S}^{2})}\leq \|h\|_{L^{\infty}(\mathbb{R}^{3}\times \mathbb{S}^{2})}+\|F\|_{L^{\infty}((0, \infty)\times\mathbb{T}^{3}\times \mathbb{S}^{2})}.
\end{equation}

According to (\ref{infinity estimate about limit solution}) and (\ref{the first term about the sequence infinity estimate}), we can easily derive the existence and estimates. Furthermore, we can get the uniqueness of equation (\ref{research equation}) according to the $L^{\infty}$ energy estimate.  \hfill$\Box$

\begin{lem}\label{main result-7}
For any $0<T_{0}<\infty$ and $m\geq 0$, if $F(t, \vec{x}, \vec{w})\in C^{m}([0, T_{0}]\times\mathbb{T}^{3}\times\mathbb{S}^{2}), h(\vec{x}, \vec{w})\in C^{m}(\mathbb{T}^{3}\times \mathbb{S}^{2})$, we can obtain that the solution $I$ in Lemma \ref{main result} also satisfies that $I\in C^{m}([0, T_{0}]\times \mathbb{T}^{3}\times\mathbb{S}^{2})$.
\end{lem}
\noindent {\bf Proof.} We can apply the operator $D_{x}^{\alpha}$ to (\ref{similar-1}) by the fact that $D_{x}^{\alpha}F(t, \vec{x}, \vec{w})\in C^{0}([0, T_{0}]\times\mathbb{T}^{3}\times \mathbb{S}^{2}), D_{x}^{\alpha}h(\vec{x}, \vec{w})\in C^{0}(\mathbb{T}^{3}\times \mathbb{S}^{2})$ where $0\leq \alpha\leq m$.

We can finally obtain that
\begin{equation}\nonumber
\|D_{x}^{\alpha}v_{k+1}\|_{L^{\infty}((0, T_{0})\times\mathbb{T}^{3}\times \mathbb{S}^{2})}\leq(\frac{1}{1+\epsilon^2})^{k}\|\epsilon^{2}D_{x}^{\alpha}I^{1}\|_{L^{\infty}((0, T_{0})\times\mathbb{T}^{3}\times \mathbb{S}^{2})},
\end{equation}
which implies $\{D_{x}^{\alpha}v_{k}\}$ is a contraction sequence. We can conclude that $D_{x}^{\alpha}I^{k}$ converges strongly in $C^{0}([0, T_{0}]\times\mathbb{T}^{3}\times \mathbb{S}^{2})$ to a limit solution $I^{\alpha}$ which satisfies
\begin{equation}\label{infinity estimate about limit solution-jia}
\|D_{x}^{\alpha}I\|_{L^{\infty}((0, T_{0})\times\mathbb{T}^{3}\times \mathbb{S}^{2})}\leq\sum_{k=1}^{\infty}\|\epsilon^{2}D_{x}^{\alpha}v_{k}\|_{L^{\infty}((0, T_{0})\times\mathbb{T}^{3}\times \mathbb{S}^{2})}\leq(1+\epsilon^{2})\|D_{x}^{\alpha}I\|_{L^{\infty}((0, T_{0})\times\mathbb{T}^{3}\times \mathbb{S}^{2})}.
\end{equation}
Since $I^{k}$ converges strongly in $C^{0}([0, T_{0}]\times\mathbb{T}^{3}\times \mathbb{S}^{2})$ to a limit solution $I$, we can verify that $I_{\alpha}=D_{x}^{\alpha}I$, which implies $D_{x}^{\alpha}I\in C^{0}([0, T_{0}]\times\mathbb{T}^{3}\times \mathbb{S}^{2})$ where $0\leq\alpha\leq m$.
Since
\begin{equation}\nonumber
\begin{split}
D_{x}^{\alpha}I^{1}(t, \vec{x}, \vec{w})=&D_{x}^{\alpha}h(\vec{x}-\frac{(\epsilon t\vec{w})}{\epsilon^{2}}, \vec{w})e^{-(1+\epsilon^{2})\frac{t}{\epsilon^{2}}} \\&+\int_{0}^{\frac{t}{\epsilon^{2}}}D_{x}^{\alpha}F(\epsilon^{2}s, \vec{x}-\epsilon(\frac{t}{\epsilon^{2}}-s)\vec{w}, \vec{w})e^{-(1+\epsilon^{2})(\frac{t}{\epsilon^{2}}-s)}ds,
\end{split}
\end{equation}
we obtain
\begin{equation}\label{the first term about the sequence infinity estimate-jia}
\|D_{x}^{\alpha}I^{1}\|_{L^{\infty}((0, T_{0})\times\mathbb{T}^{3}\times \mathbb{S}^{2})}\leq \|D_{x}^{\alpha}h\|_{L^{\infty}(\mathbb{T}^{3}\times \mathbb{S}^{2})}+\|D_{x}^{\alpha}F\|_{L^{\infty}((0, T_{0})\times\mathbb{T}^{3}\times \mathbb{S}^{2})}.
\end{equation}
Therefore, we have
\begin{equation}\label{estimate result-1-1}
\|D_{x}^{\alpha}I\|_{L^{\infty}((0, T_{0})\times\mathbb{T}^{3}\times \mathbb{S}^{2})}\leq C(\|D_{x}^{\alpha}h\|_{L^{\infty}(\mathbb{T}^{3}\times \mathbb{S}^{2})}+\|D_{x}^{\alpha}F\|_{L^{\infty}((0, T_{0})\times\mathbb{T}^{3}\times \mathbb{S}^{2})}).
\end{equation}
which implies $D_{x}^{\alpha}I\in C^{0}([0, T_{0}]\times\mathbb{T}^{3}\times \mathbb{S}^{2})$ where $0\leq \alpha\leq m$. \hfill$\Box$

\begin{lem}\label{main result-7-bu}
For any $0<T_{0}<\infty$, if $F(t, \vec{x}, \vec{w})\in L^{2}(0, T_{0}; L^{2}(\mathbb{S}^{2}; H^{m}(\mathbb{T}^{3}))), h(\vec{x}, \vec{w})\in L^{2}(\mathbb{S}^{2}; H^{m}(\mathbb{T}^{3}))$, there exists a unique nonnegative solution $I(t, \vec{x}, \vec{w})\in L^{\infty}(0, T_{0}; L^{2}(\mathbb{S}^{2}; \\H^{m}(\mathbb{T}^{3})))$ of the transport equation (\ref{research equation}) which satisfies
\begin{equation}\label{L2 energy estimate}
\begin{aligned}
\|I\|_{L^{\infty}(0, T_{0}; L^{2}(\mathbb{S}^{2}; H^{m}(\mathbb{T}^{3})))}\leq C(\|h\|_{L^{2}(\mathbb{S}^{2}; H^{m}(\mathbb{T}^{3}))}+\frac{1}{\epsilon}\|F\|_{L^{2}(0, T_{0}; L^{2}(\mathbb{S}^{2}; H^{m}(\mathbb{T}^{3})))})
\end{aligned}
\end{equation}
\end{lem}
\noindent {\bf Proof.} We can construct sequences $\{F_{k}\}_{k=1}^{\infty}\subset C^{m}([0, T_{0}]\times\mathbb{T}^{3}\times\mathbb{S}^{2})$ and $\{h_{k}\}_{k=1}^{\infty}\subset C^{m}(\overline{\mathbb{T}^{3}}\times\mathbb{S}^{2})$ satisfying $\lim_{k\rightarrow\infty}\|F_{k}-F\|_{L^{2}(0, T_{0}; L^{2}(\mathbb{S}^{2}; H^{m}(\mathbb{T}^{3})))}=0$ and $\lim_{k\rightarrow\infty}\|h_{k}-h\|_{L^{2}(\mathbb{S}^{2}; H^{m}(\mathbb{T}^{3}))}\\=0$. Then, we consider the following equations
\begin{equation}\label{approximating sequence equation-bu}\left\{
\begin{split}
&{\epsilon}^{2}\frac{\partial I_{k}}{\partial t}+\epsilon \vec w\cdot \nabla_{x}I_{k}+\epsilon^{2}(I_{k}
-\overline {I_{k}})+I_{k}=F_{k}(t, \vec{x}, \vec{w}),~~\mathrm{in} \ \ (0, T_{0})\times\mathbb{T}^{3}\times \mathbb{S}^{2},\\
&I_{k}(0, \vec{x}, \vec{w})=h_{k}(\vec{x}, \vec{w}),~~\mathrm{in} \ \ \mathbb{T}^{3}\times \mathbb{S}^{2}.\\
\end{split}\right.
\end{equation}
According to Lemma \ref{main result-7}, we have $D_{x}^{\alpha}I_{k}\in C^{0}([0, T_{0}]\times\mathbb{T}^{3}\times\mathbb{S}^{2})$ where $0\leq l\leq m$. Applying the operator on both sides of $(\ref{approximating sequence equation-bu})_{1}$, multiplying $D_{x}^{l}I_{k}$ on both sides of $(\ref{approximating sequence equation-bu})_{1}$ and integrating on $(0, t)\times\mathbb{T}^{3}\times\mathbb{S}^{2}$ for $t\in [0, T_{0}]$, then, we get
\begin{equation}\nonumber
\begin{aligned}
&\int_{\mathbb{T}^{3}}\int_{\mathbb{S}^{2}}|D_{x}^{\alpha}I_{k}|^{2}(t)\mathrm{d}\vec x\mathrm{d}\vec w+\int_{0}^{t}\int_{\mathbb{T}^{3}}\int_{\mathbb{S}^{2}}|D_{x}^{\alpha}(I_{k}-\overline {I_{k}})|^{2}\mathrm{d} t\mathrm{d}\vec x\mathrm{d}\vec w+\frac{1}{\epsilon^{2}}\int_{0}^{t}\int_{\mathbb{T}^{3}}\int_{\mathbb{S}^{2}}|D_{x}^{\alpha}I_{k}|^{2}\\=&\int_{\mathbb{T}^{3}}\int_{\mathbb{S}^{2}}|D_{x}^{l}h_{k}|^{2}\mathrm{d}\vec x\mathrm{d} \vec w
+\frac{1}{\epsilon^{2}}\int_{0}^{t}\int_{\mathbb{T}^{3}}\int_{\mathbb{S}^{2}}D_{x}^{\alpha}F_{k}D_{x}^{\alpha}I_{k}\mathrm{d}t\mathrm{d}\vec x\mathrm{d}\vec w
\\\leq& \int_{\mathbb{T}^{3}}\int_{\mathbb{S}^{2}}|D_{x}^{\alpha}h_{k}|^{2}\mathrm{d}\vec x\mathrm{d} \vec w
+\frac{1}{2\epsilon^{2}}\int_{0}^{t}\int_{\mathbb{T}^{3}}\int_{\mathbb{S}^{2}}|D_{x}^{\alpha}F_{k}|^{2}\mathrm{d}t\mathrm{d}\vec x\mathrm{d}\vec w+\frac{1}{2\epsilon^{2}}\int_{0}^{t}\int_{\mathbb{T}^{3}}\int_{\mathbb{S}^{2}}|D_{x}^{\alpha}I_{k}|^{2}\mathrm{d}t\mathrm{d}\vec x\mathrm{d}\vec w.
\end{aligned}
\end{equation}
So,
\begin{equation}\nonumber
\begin{aligned}
\sup_{t\in[0, T_{0}]}\int_{\mathbb{T}^{3}}\int_{\mathbb{S}^{2}}|D_{x}^{\alpha}I_{k}|^{2}(t)\mathrm{d}\vec x\mathrm{d}\vec w\leq \int_{\mathbb{T}^{3}}\int_{\mathbb{S}^{2}}|D_{x}^{\alpha}h_{k}|^{2}\mathrm{d}\vec x\mathrm{d} \vec w
+\frac{1}{2\epsilon^{2}}\int_{0}^{T_{0}}\int_{\mathbb{T}^{3}}\int_{\mathbb{S}^{2}}|D_{x}^{\alpha}F_{k}|^{2}\mathrm{d}t\mathrm{d}\vec x\mathrm{d}\vec w.
\end{aligned}
\end{equation}
Letting $k\rightarrow\infty$, we have
\begin{equation}\nonumber
\begin{aligned}
\sup_{t\in[0, T_{0}]}\int_{\mathbb{T}^{3}}\int_{\mathbb{S}^{2}}|I^{\alpha}|^{2}(t)\mathrm{d}\vec x\mathrm{d}\vec w\leq \int_{\mathbb{T}^{3}}\int_{\mathbb{S}^{2}}|D_{x}^{\alpha}h|^{2}\mathrm{d}\vec x\mathrm{d} \vec w
+\frac{1}{2\epsilon^{2}}\int_{0}^{T_{0}}\int_{\mathbb{T}^{3}}\int_{\mathbb{S}^{2}}|D_{x}^{\alpha}F|^{2}\mathrm{d}t\mathrm{d}\vec x\mathrm{d}\vec w,
\end{aligned}
\end{equation}
where $I^\alpha$ is the weak limit of $D_{x}^{\alpha}I_{k}$ in $L^{\infty}(0, T_{0}; L^{2}(\mathbb{S}^{2}\times\mathbb{T}^{3}))$ in the weak-$\ast$ sense up to a subsequence. We can deduce that $I^{\alpha}=D_{x}^{\alpha}I$ by the fact that $I_{k}\rightarrow I$ in $L^{\infty}(0, T_{0}; L^{2}(\mathbb{S}^{2}\times\mathbb{T}^{3}))$ in the weak-$\ast$ sense up to the same subsequence. We can also verify that $I(0, \vec x, \vec w)=h(\vec x, \vec w)$ by referring to suitable test functions and integrating by parts. The method is standard, we omit the details here. So, we have proved the existence and estimate \eqref{L2 energy estimate}. The uniqueness of equation (\ref{research equation}) can be ensured by (\ref{L2 energy estimate}). \hfill$\Box$

\begin{lem}\label{lianxujieguo}
For the solution $I(t, \vec{x}, \vec{w})$ we obtained in Lemma \ref{main result-7-bu}, it further satisfies that $I(t, \vec{x}, \vec{w})\in C^{0}([0, T_{0}]; L^{2}(\mathbb{S}^{2}; H^{m}(\mathbb{T}^{3})))$.
\end{lem}
\noindent\textbf{Proof.} Define the Friedrichs mollifier $\phi_{\delta}$. Applying the operator on both sides of \eqref{research equation}, we derive 
\begin{equation}\label{research equation-m1}\left\{
\begin{split}
&{\epsilon}^{2}\frac{\partial I_{\delta}}{\partial t}+\epsilon \vec{w}\cdot \nabla_{x}I_{\delta}+{\epsilon}^{2}(I_{\delta}-\overline {I_{\delta}})+I_{\delta}=F_{\delta}(t, \vec{x}, \vec{w})~~\mathrm{in} \ \ (0, \infty)\times\mathbb{T}^{3}\times \mathbb{S}^{2},\\
&I_{\delta}(0, \vec{x}, \vec{w})=h_{\delta}(\vec{x}, \vec{w})~~\mathrm{in} \ \ \mathbb{T}^{3}\times \mathbb{S}^{2}.\\
\end{split}\right.
\end{equation}
Define $I_{\delta, \delta'}=I_{\delta}-I_{\delta'}$. Then, we can write the equation satisfied by $I_{\delta, \delta'}$ as follows
\begin{equation}\label{research equation-m2}\left\{
\begin{split}
&{\epsilon}^{2}\frac{\partial I_{\delta, \delta'}}{\partial t}+\epsilon \vec{w}\cdot \nabla_{x}I_{\delta, \delta'}+{\epsilon}^{2}(I_{\delta, \delta'}-\overline {I_{\delta, \delta'}})+I_{\delta, \delta'}=F_{\delta, \delta'}(t, \vec{x}, \vec{w})~~\mathrm{in} \ \ (0, T_{0})\times\mathbb{T}^{3}\times \mathbb{S}^{2},\\
&I_{\delta, \delta'}(0, \vec{x}, \vec{w})=h_{\delta, \delta'}(\vec{x}, \vec{w})~~\mathrm{in} \ \ \mathbb{T}^{3}\times \mathbb{S}^{2},\\
\end{split}\right.
\end{equation}
where $h_{\delta, \delta'}=h_{\delta}-h_{\delta'}$ and $F_{\delta, \delta'}=F_{\delta}-F_{\delta'}$.

According to the estimate \eqref{L2 energy estimate}, we have
\begin{equation}\label{L2 energy estimate-m2}
\begin{aligned}
\|I_{\delta, \delta'}\|_{L^{\infty}(0, T_{0}; L^{2}(\mathbb{S}^{2}; H^{m}(\mathbb{T}^{3})))}\leq C(\|h_{\delta, \delta'}\|_{L^{2}(\mathbb{S}^{2}; H^{m}(\mathbb{T}^{3}))}+\frac{1}{\epsilon}\|F_{\delta, \delta'}\|_{L^{2}(0, T_{0}; L^{2}(\mathbb{S}^{2}; H^{m}(\mathbb{T}^{3})))}). 
\end{aligned}
\end{equation}
As $\delta, \delta'\rightarrow 0$, $h_{\delta, \delta'}\rightarrow 0$ and $F_{\delta, \delta'}\rightarrow 0$ in $L^{2}(\mathbb{S}^{2}; H^{m}(\mathbb{T}^{3}))$ and $L^{2}(0, T_{0}; L^{2}(\mathbb{S}^{2}; H^{m}(\mathbb{T}^{3})))$, respectively, we derive that $I_{\delta, \delta'}\rightarrow 0$ in $L^{\infty}(0, T_{0}; L^{2}(\mathbb{S}^{2}; H^{m}(\mathbb{T}^{3})))$. Then, we confirm that $I_{\delta}$ is a Cauchy sequence in $C^{0}([0, T_{0}]; L^{2}(\mathbb{S}^{2}; H^{m}(\mathbb{T}^{3})))$ which further implies $I\in C^{0}([0, T_{0}]; L^{2}(\mathbb{S}^{2}; H^{m}(\mathbb{T}^{3})))$. \hfill$\Box$

\section{Hilbert Expansion}\label{san}
\subsection{Interior expansion}
First, we give the interior expansion of $\theta^{\epsilon}$ as follows:
\begin{equation}\label{expansion of theta epsilon}
\theta^{\epsilon}\sim\sum_{k=0}^{N}\epsilon^{k}\theta_{k}.
\end{equation}
Then, we get
\begin{equation}\label{expansion of B(theta epsilon)}
\begin{split}
B(\theta^{\epsilon})&\sim(\theta_{0}+\epsilon\theta_{1}+\cdot\cdot\cdot)^{4}
\\&=\theta_{0}^{4}+4\theta_{0}^{3}\theta_{1}\epsilon+(4\theta_{2}\theta_{0}^{3}+6\theta_{1}^{2}\theta_{0}^{2})\epsilon^{2}+\cdot\cdot\cdot
\\&=B_{0}(\theta_{0})+\epsilon B_{1}(\theta_{0}; \theta_{1})+\epsilon^{2}B_{2}(\theta_{0}; \theta_{2})+\cdot\cdot\cdot
\end{split}
\end{equation}
where 
\begin{equation}\label{defBk}
B_{k}(\theta_{0}; \theta_{k}):=\sum_{\begin{matrix}
  i+j+l+m=k, \\
  i,j,l,m\geq 0
\end{matrix}}\theta_{i}\theta_{j}\theta_{l}\theta_{m}.
\end{equation}
Define the interior expansion of $f^{\epsilon}(t, \vec{x}, w)$ as
\begin{equation}\label{expansion about f epsilon including k term}
f^{\epsilon}(t, \vec x, \vec w)\sim \sum_{k=0}^{N}{\epsilon}^{k}f_{k}(t, \vec x, \vec w).
\end{equation}
Define
\begin{equation}\label{suanzi1}
\mathcal{L}_{1}(f^{\epsilon}, \theta^{\epsilon}):=\epsilon^{2}\partial_{t}f^{\epsilon}+\epsilon \vec w\cdot \nabla_{x}f^{\epsilon}+{\epsilon}^{2}(f^{\epsilon}-\overline {f^{\epsilon}})+f^{\epsilon}-B(\theta^{\epsilon})
\end{equation}
and
\begin{equation}\label{suanzi2}
\mathcal{L}_{2}(f^{\epsilon}, \theta^{\epsilon}):=\epsilon^{2}\partial_{t}\theta^{\epsilon}+\epsilon^{2}\mathrm{div}(\vec{u}\theta^{\epsilon})
-\epsilon^{2}\Delta_{x} \theta^{\epsilon}
-(\overline {f^{\epsilon}}-B(\theta^{\epsilon})).
\end{equation}
Plugging \eqref{expansion of theta epsilon} and \eqref{expansion about f epsilon including k term} into the above formulas gives
\begin{equation}\label{resi1}
\begin{aligned}
\mathcal{L}_{1}(\sum_{k=0}^{N}{\epsilon}^{k}f_{k}(t, \vec x, \vec w), \sum_{k=0}^{N}\epsilon^{k}\theta_{k})=&\sum_{k=0}^{N}{\epsilon}^{k}(\partial_{t}f_{k-2}+\vec{w}\cdot\nabla_{x}f_{k-1}
+f_{k-2}-\overline{f_{k-2}}+f_{k}-B_{k}(\theta_{0};\theta_{k}))
\\&+\epsilon^{N+1}\partial_{t}f_{N-1}+\epsilon^{N+2}\partial_{t}f_{N}+\epsilon^{N+1}\vec{w}\cdot\nabla_{x}f_{N}
\\&-\sum_{k=N+1}^{4N}\epsilon^{k}\sum_{\begin{matrix}
  i+j+l+m=k, \\
  0\leq i,j,l,m\leq N
\end{matrix}}\theta_{i}\theta_{j}\theta_{l}\theta_{m}
\end{aligned}
\end{equation}
and
\begin{equation}\label{resi2}
\begin{aligned}
\mathcal{L}_{2}(\sum_{k=0}^{N}{\epsilon}^{k}f_{k}(t, \vec x, \vec w), \sum_{k=0}^{N}\epsilon^{k}\theta_{k})=&\sum_{k=0}^{N}{\epsilon}^{k}(\partial_{t}\theta_{k-2}
+\mathrm{div}(\vec{u}\theta_{k-2})-\Delta_{x}\theta_{k-2}
-(\overline{f_{k}}-B_{k}(\theta_{0};\theta_{k})))\\&+\epsilon^{N+1}\partial_{t}\theta_{N-1}
+\epsilon^{N+2}\partial_{t}\theta_{N}+\epsilon^{N+1}\Delta_{x}\theta_{N-1}+\epsilon^{N+2}\Delta_{x}\theta_{N}
\\&+\epsilon^{N+1}\mathrm{div}(\vec{u}\theta_{N-1})+\epsilon^{N+2}\mathrm{div}(\vec{u}\theta_{N})\\&+\sum_{k=N+1}^{4N}\epsilon^{k}\sum_{\begin{matrix}
  i+j+l+m=k, \\
  0\leq i,j,l,m\leq N
\end{matrix}}\theta_{i}\theta_{j}\theta_{l}\theta_{m},
\end{aligned}
\end{equation}
where $(f_{k}, \theta_{k})$, $k<0$ are taken to be zero.

Collecting the terms with the same order, we take
\begin{equation}\label{fk-2}
\partial_{t}f_{k-2}+\vec{w}\cdot\nabla_{x}f_{k-1}
+f_{k-2}-\overline{f_{k-2}}+f_{k}-B_{k}(\theta_{0};\theta_{k})=0,
\end{equation}

\begin{equation}\label{thetak-2}
\partial_{t}\theta_{k-2}
+\mathrm{div}(\vec{u}\theta_{k-2})-\Delta_{x}\theta_{k-2}
-(\overline{f_{k}}-B_{k}(\theta_{0};\theta_{k}))=0,
\end{equation}
for any $k=0, 1,\cdot\cdot\cdot, N$. From the above two equations we obtain
\begin{equation}
\partial_{t}\theta_{k-2}+\mathrm{div}(\vec u\theta_{k-2})-\Delta_{x}\theta_{k-2}=-\partial_{t}\overline{f_{k-2}}-\langle\vec{w}\cdot\nabla_{x}f_{k-1}\rangle.
\end{equation}
By \eqref{fk-2}, an iteration process on the above equation leads to 
\begin{equation}
\begin{aligned}
&\partial_{t}\theta_{k-2}+\mathrm{div}(\vec u\theta_{k-2})-\Delta_{x}\theta_{k-2}\\=&-\partial_{t}\overline{f_{k-2}}
-\langle\vec{w}\cdot\nabla_{x}(B_{k-1}(\theta_{0};\theta_{k-1})
+\overline{f_{k-3}}-f_{k-3}-\vec{w}\cdot\nabla_{x}f_{k-2}-\partial_{t}f_{k-3})\rangle
\\=&-\partial_{t}\overline{f_{k-2}}
+\langle\vec{w}\cdot\nabla_{x}f_{k-3}+(\vec{w}\cdot\nabla_{x})^{2}f_{k-2}+\vec{w}\cdot\nabla_{x}\partial_{t}f_{k-3}\rangle
\\=&-\partial_{t}\langle(B_{k-2}(\theta_{0};\theta_{k-2})+\overline{f_{k-4}}-f_{k-4}-\vec{w}\cdot\nabla_{x}f_{k-3}
-\partial_{t}f_{k-4})\rangle+\langle\vec{w}\cdot\nabla_{x}f_{k-3}\\&+(\vec{w}\cdot\nabla_{x})^{2}(B_{k-2}(\theta_{0};\theta_{k-2})
+\overline{f_{k-4}}-f_{k-4}-\vec{w}\cdot\nabla_{x}f_{k-3}
-\partial_{t}f_{k-4})
+\vec{w}\cdot\nabla_{x}\partial_{t}f_{k-3}\rangle.
\end{aligned}
\end{equation}
Consequently, equations \eqref{fk-2}-\eqref{thetak-2} can be rewritten as
\begin{equation}\label{thetak}
\partial_{t}(\theta_{k}+B_{k}(\theta_{0};\theta_{k}))+\mathrm{div}(\vec u\theta_{k})-\Delta_{x}\Big(\theta_{k}+\frac{1}{3}B_{k}(\theta_{0};\theta_{k})\Big)=F_{k-1},
\end{equation}
\begin{equation}\label{fk}
f_{k}=B_{k}(\theta_{0};\theta_{k})+\overline{f_{k-2}}-f_{k-2}-\partial_{t}f_{k-2}-\vec{w}\cdot\nabla_{x}f_{k-1},
\end{equation}
where
\begin{equation}
\begin{aligned}
F_{k-1}=&-\partial_{t}\langle\overline{f_{k-2}}-f_{k-2}+\vec{w}\cdot\nabla_{x}f_{k-1}-\partial_{t}f_{k-2}\rangle
+\langle\vec{w}\cdot\nabla_{x}f_{k-1}+(\vec{w}\cdot\nabla_{x})^{2}(\overline{f_{k-2}}-f_{k-2})
\\&-(\vec{w}\cdot\nabla_{x})^{3}f_{k-1}
-(\vec{w}\cdot\nabla_{x})^{2}\partial_{t}f_{k-2}+\vec{w}\cdot\nabla_{x}\partial_{t}f_{k-1}\rangle,
\end{aligned}
\end{equation}
for any $k=0, 1, \cdot\cdot\cdot, N$.

Note that $\theta_{k}$ is solved by  \eqref{thetak} and the solution is plugged into \eqref{fk} to get $f_{k}$. In addition, for $k=0$, equation \eqref{thetak} becomes
\begin{equation}
\partial_{t}(\theta_{0}+B(\theta_{0}))+\mathrm{div}(\vec{u}\theta_{0})-\Delta_{x}\Big(\theta_{0}+\frac{1}{3}B(\theta_{0})\Big)=0.
\end{equation}
For $k\geq 1$, equation \eqref{thetak} can be rewritten as the following linear parabolic equation
\begin{equation}
\begin{aligned}
\partial_{t}(\theta_{k}+4\theta_{0}^{3}\theta_{k})+\mathrm{div}(\vec u\theta_{k})-\Delta_{x}\Big(\theta_{k}+\frac{4}{3}\theta_{0}^{3}\theta_{k}\Big)=\tilde{F}_{k-1},
\end{aligned}
\end{equation}
where
\begin{equation}
\tilde{F}_{k-1}=F_{k-1}-\frac{4}{3}\Delta_{x}\Big(\sum_{\begin{matrix}
  i+j+l+m=k, \\
  i,j,l,m\geq 1
\end{matrix}}\theta_{i}\theta_{j}\theta_{l}\theta_{m}\Big).
\end{equation}
From the above equations, residuals \eqref{resi1} and \eqref{resi2} are given by
\begin{equation}\label{resi1-jia}
\begin{aligned}
\mathcal{L}_{1}(\sum_{k=0}^{N}{\epsilon}^{k}f_{k}, \sum_{k=0}^{N}\epsilon^{k}\theta_{k})=&\epsilon^{N+1}\partial_{t}f_{N-1}+\epsilon^{N+2}\partial_{t}f_{N}+\epsilon^{N+1}\vec{w}\cdot\nabla_{x}f_{N}\\&-\sum_{k=N+1}^{4N}\epsilon^{k}\sum_{\begin{matrix}
  i+j+l+m=k, \\
  0\leq i,j,l,m\leq N
\end{matrix}}\theta_{i}\theta_{j}\theta_{l}\theta_{m}
\end{aligned}
\end{equation}
and
\begin{equation}\label{resi2-jia}
\begin{aligned}
\mathcal{L}_{2}(\sum_{k=0}^{N}{\epsilon}^{k}f_{k}(t, \vec x, \vec w), \sum_{k=0}^{N}\epsilon^{k}\theta_{k})=&\epsilon^{N+1}\partial_{t}\theta_{N-1}
+\epsilon^{N+2}\partial_{t}\theta_{N}+\epsilon^{N+1}\Delta_{x}\theta_{N-1}+\epsilon^{N+2}\Delta_{x}\theta_{N}
\\&+\epsilon^{N+1}\mathrm{div}(\vec{u}\theta_{N-1})+\epsilon^{N+2}\mathrm{div}(\vec{u}\theta_{N})\\&+\sum_{k=N+1}^{4N}\epsilon^{k}\sum_{\begin{matrix}
  i+j+l+m=k, \\
  0\leq i,j,l,m\leq N
\end{matrix}}\theta_{i}\theta_{j}\theta_{l}\theta_{m},
\end{aligned}
\end{equation}
where the right hand side are both formally of order $\epsilon^{N+1}$.

\subsection{Initial layer expansion.}
In order to determine the initial condition for $\overline {f_{k}}$ and $\theta_{k}$, we resort to the initial layer expansion. We define the new variable $\tau$ by making the scaling transform for $(f^{\epsilon}(t, \vec{x}, \vec w), \theta^{\epsilon}(t, \vec x))\rightarrow (f^{\epsilon}(\tau, \vec{x}, \vec w), \theta^{\epsilon}(\tau, \vec x))$ with $\tau\in [0, \infty)$ as
\begin{equation}\nonumber
\tau=\frac{t}{\epsilon^{2}},
\end{equation}
which implies
\begin{equation}\nonumber
\frac{\partial f^{\epsilon}}{\partial t}=\frac{1}{\epsilon^{2}}\frac{\partial f^{\epsilon}}{\partial\tau}, \ \ \frac{\partial \theta^{\epsilon}}{\partial t}=\frac{1}{\epsilon^{2}}\frac{\partial \theta^{\epsilon}}{\partial\tau}.
\end{equation}
In this new variable, the system $(\ref{research equations})$ can be rewritten as
\begin{equation}\label{research equation in the initial layer}\left\{
\begin{split}
&\frac{\partial f^{\epsilon}}{\partial\tau}+\epsilon \vec w\cdot \nabla_{x}f^{\epsilon}+\epsilon^{2}(f^{\epsilon}-\overline {f^{\epsilon}})+f^{\epsilon}=B(\theta^{\epsilon}),\\
&\frac{1}{\epsilon^{2}}\frac{\partial \theta^{\epsilon}}{\partial\tau}+\mathrm{div}(\vec{u}\theta^{\epsilon})-\Delta_{x}\theta^{\epsilon}
=\frac{1}{\epsilon^{2}}(\overline {f^{\epsilon}}-B(\theta^{\epsilon})).\\
\end{split}\right.
\end{equation}
We take the ansatz as follows:
\begin{equation}\label{expansion about the initial layer}
f^{\epsilon}\thicksim \sum_{k=0}^{N}\epsilon^{k}(f_{k}(t, \vec x, \vec w)+f_{I, k}(\tau, \vec{x}, \vec{w})),\;\;\theta^{\epsilon}\thicksim\sum_{k=0}^{N}\epsilon^{k}(\theta_{k}(t, \vec x)+\theta_{I, k}(\tau, \vec{x})).
\end{equation}
Taking the above ansatz into \eqref{suanzi1} and \eqref{suanzi2} gives
\begin{equation}\nonumber
\begin{aligned}
&\mathcal{L}_{1}\Big(\sum_{k=0}^{N}{\epsilon}^{k}(f_{k}(t, \vec x, \vec w)+f_{I, k}(\tau, \vec x, \vec w)), \sum_{k=0}^{N}\epsilon^{k}(\theta_{k}(t, \vec x)+\theta_{I, k}(\tau, \vec x))\Big)\\=&\mathcal{L}_{1}\Big(\sum_{k=0}^{N}{\epsilon}^{k}f_{k}, \sum_{k=0}^{N}\epsilon^{k}\theta_{k}\Big)-\sum_{k=0}^{N}{\epsilon}^{k}(B_{k}(\theta_{0}+\theta_{I, 0};\theta_{k}+\theta_{I, k})-B_{k}(\theta_{0};\theta_{k}))
\\&+\sum_{k=0}^{N}{\epsilon}^{k}(B_{k}(A_{0}+\theta_{I, 0};A_{k}+\theta_{I, k})-B_{k}(A_{0};A_{k}))\\&-\sum_{k=N+1}^{4N}\epsilon^{k}\sum_{\begin{matrix}
  i+j+l+m=k, \\
  0\leq i,j,l,m\leq N
\end{matrix}}(\theta_{i}+\theta_{I, i})(\theta_{j}+\theta_{I, j})(\theta_{l}+\theta_{I, l})(\theta_{m}+\theta_{I, m})
\end{aligned}
\end{equation}
\begin{equation}\label{resi1-zaijia}
\begin{aligned}
\\&+\sum_{k=0}^{N}(\partial_{\tau}f_{I, k}+\vec{w}\cdot\nabla_{x}f_{I, k-1}+f_{I, k}+f_{I, k-2}-\overline{f_{I, k-2}}-(B_{k}(A_{0}+\theta_{I, 0};A_{k}+\theta_{I, k})-B_{k}(A_{0};A_{k})))\\&+\epsilon^{N+1}(f_{I, N-1}-\overline{f_{I, N-1}})+\epsilon^{N+2}(f_{I, N}-\overline{f_{I, N}})+\epsilon^{N+1}\vec{w}\cdot\nabla_{x}f_{I, N}
\end{aligned}
\end{equation}
and
\begin{equation}\nonumber
\begin{aligned}
&\mathcal{L}_{2}\Big(\sum_{k=0}^{N}{\epsilon}^{k}(f_{k}(t, \vec x, \vec w)+f_{I, k}(\tau, \vec x, \vec w)), \sum_{k=0}^{N}\epsilon^{k}(\theta_{k}(t, \vec x)+\theta_{I, k}(\tau, \vec x))\Big)\\=&\mathcal{L}_{2}\Big(\sum_{k=0}^{N}{\epsilon}^{k}f_{k}, \sum_{k=0}^{N}\epsilon^{k}\theta_{k}\Big)+\sum_{k=0}^{N}{\epsilon}^{k}(B_{k}(\theta_{0}+\theta_{I, 0};\theta_{k}+\theta_{I, k})-B_{k}(\theta_{0};\theta_{k}))
\end{aligned}
\end{equation}
\begin{equation}\label{resi2-zaijia}
\begin{aligned}
\\&-\sum_{k=0}^{N}{\epsilon}^{k}(B_{k}(A_{0}+\theta_{I, 0};A_{k}+\theta_{I, k})-B_{k}(A_{0};A_{k}))\\&+\sum_{k=N+1}^{4N}\epsilon^{k}\sum_{\begin{matrix}
  i+j+l+m=k, \\
  0\leq i,j,l,m\leq N
\end{matrix}}(\theta_{i}+\theta_{I, i})(\theta_{j}+\theta_{I, j})(\theta_{l}+\theta_{I, l})(\theta_{m}+\theta_{I, m})
\\&+\sum_{k=0}^{N}{\epsilon}^{k}(\partial_{\tau}\theta_{I, k}+\mathrm{div}(\vec v_{k}\theta_{I, k-2})-\overline{f_{I, k}}-\Delta_{x}\theta_{I, k-2}+(B_{k}(A_{0}+\theta_{I, 0};A_{k}+\theta_{I, k})-B_{k}(A_{0};A_{k})))\\&+\epsilon^{N+1}(\mathrm{div}(\vec u\theta_{I, k-1})-\Delta_{x}\theta_{I, k-1})+\epsilon^{N+2}(\mathrm{div}(\vec u\theta_{I, k})-\Delta_{x}\theta_{I, k})+\sum_{k=0}^{N}\epsilon^{k}\mathrm{div}((\vec u-\vec v_{k})\theta_{I, k-2}),
\end{aligned}
\end{equation}
where $\theta_{k}(0)=\theta_{k}(t=0)$, $A_{k}=A_{k}(\tau, \vec x)$ and $v_{k}=v_{k}(\tau, \vec x)$, $k=0, 1, \cdot\cdot\cdot, N$ are given by
\begin{equation}\label{formula1}
A_{k}(\tau, \vec x)=\sum_{l=0}^{k}\epsilon^{l}\frac{\tau^{l}}{l!}\frac{\partial^{l}}{\partial t^{l}}\theta_{k-l}(0, \vec x)
\end{equation}
and 
\begin{equation}\label{formula2}
\vec{v}_{k}(\tau, \vec x)=\sum_{l=0}^{k}\epsilon^{l}\frac{\tau^{l}}{l!}\frac{\partial^{l}}{\partial t^{l}}\vec{u}(0, \vec x).
\end{equation}

Collecting terms of the same order in \eqref{resi1-zaijia}-\eqref{resi2-zaijia}, we take
\begin{equation}\label{ini-f}
\partial_{\tau}f_{I, k}+\vec{w}\cdot\nabla_{x}f_{I, k-1}+f_{I, k}+f_{I, k-2}-\overline{f_{I, k-2}}-(B_{k}(A_{0}+\theta_{I, 0};A_{k}+\theta_{I, k})-B_{k}(A_{0};A_{k}))=0
\end{equation}
and
\begin{equation}
\partial_{\tau}\theta_{I, k}+\mathrm{div}(\vec v_{k}f_{I, k-2})-\overline{f_{I, k}}-\Delta_{x}f_{I, k-2}+(B_{k}(A_{0}+\theta_{I, 0};A_{k}+\theta_{I, k})-B_{k}(A_{0};A_{k}))=0,
\end{equation}
for $k=0, 1, \cdot\cdot\cdot, N$.

\subsection{Construction of asymptotic expansion.}\label{construction}
The bridge between the interior solution with the initial layer is the initial condition of (\ref{research equations}). First, we have the following relations
\begin{equation}\nonumber
f_{0}(0, \vec{x}, \vec{w})+f_{I, 0}(0, \vec{x}, \vec{w})=h(\vec{x}, \vec{w}),\hspace{0.5 cm}
\end{equation}
\begin{equation}\nonumber
\theta_{0}(0, \vec{x})+\theta_{I, 0}(0, \vec{x})=\theta^{0}(\vec{x}),\hspace{-0.2 cm}
\end{equation}
\begin{equation}\nonumber
f_{k}(0, \vec{x}, \vec{w})+f_{I, k}(0, \vec{x}, \vec{w})=0,k\geq 1,\hspace{0.3cm}
\end{equation}
\begin{equation}\nonumber
\theta_{k}(0, \vec{x})+\theta_{I, k}(0, \vec{x})=0,k\geq 1.\hspace{-0.7 cm}
\end{equation}
The construction of $f_{k}, f_{I, k}$ and $\theta_{k}, \theta_{I, k}$ are as follows:

\noindent\textbf{Step 1.}Construction of zeroth-order terms.

The zeroth-order initial layer $(f_{I, 0}, \theta_{I, 0})$ is defined as
\begin{equation}\label{equations about f I0,theta I0}\left\{
\begin{split}
&\frac{\partial f_{I, 0}}{\partial\tau}+f_{I, 0}=B(\theta_{0}^{0}+\theta_{I, 0})-B(\theta_{0}^{0}),\\
&\frac{\partial \theta_{I, 0}}{\partial\tau}=(\overline{f_{I, 0}}-B(\theta_{0}^{0}+\theta_{I, 0})+B(\theta_{0}^{0})),\\
&f_{I, 0}(0, \vec{x}, \vec{w})=h(\vec{x}, \vec{w})-\overline {f_{0}}(0, \vec x),\\
&\theta_{I, 0}(0, \vec{x}, \vec{w})=\theta^{0}(\vec x)-\theta_{0}^{0}(\vec x),\\
&\lim_{\tau\rightarrow\infty}f_{I, 0}(\tau, \vec{x}, \vec{w})=0,\ \ \lim_{\tau\rightarrow\infty}\theta_{I, 0}(\tau, \vec{x}, \vec{w})=0,
\end{split}\right.
\end{equation}
where we have used the notation $\theta_{0}^{0}=\theta_{0}|_{t=0}$.
\begin{thm}\label{3.1}
For the given data $h(\vec{x}, \vec{w}), \theta^{0}(\vec{x}), \vec{u}^{0}(\vec x)$ in Theorem \ref{result}, the problem \eqref{equations about f I0,theta I0} has a unique solution $(f_{I, 0}, \theta_{I, 0})\in (C^{1}([0, \infty); L^{\infty}_{\vec w}H_{\vec x}^{N+2}))\cap L^{1}([0, \infty); L^{\infty}_{\vec w}H_{\vec x}^{4N+2})))\times (C^{1}([0, \infty); H_{\vec x}^{N+2})\cap L^{1}([0, \infty); H_{\vec x}^{N+2}))$. Furthermore, we have
\begin{equation}\label{ineq-theta00-jia0}
\|\theta_{0}^{0}(\vec x)\|_{H_{\vec x}^{N+2}}\leq C,
\end{equation}
\begin{equation}\label{ineq-theta00-jia1}
\|e^{\sigma\tau}f_{I, 0}\|_{L^{\infty}([0, \infty); L^{\infty}_{\vec w}H_{\vec x}^{N+2})}+\|e^{\sigma\tau}\theta_{I, 0}\|_{L^{\infty}([0, \infty); H_{\vec x}^{N+2})}\leq C\eta
\end{equation}
and
\begin{equation}\label{ineq-theta00-jia2}
\|e^{\sigma\tau}f_{I, 0}\|_{L^{1}([0, \infty); L^{\infty}_{\vec w}H_{\vec x}^{4N+2})}+\|e^{\sigma\tau}\theta_{I, 0}\|_{L^{1}([0, \infty); H_{\vec x}^{N+2})}\leq C\eta,
\end{equation}
where $\sigma>0$, suitably small.
\end{thm}
\noindent\textbf{Proof.} From the equations $\eqref{equations about f I0,theta I0}_{1}$ and $\eqref{equations about f I0,theta I0}_{2}$, we can derive 
\begin{equation}
\partial_{\tau}(\overline {f_{I, 0}}+\theta_{I, 0})=0.
\end{equation}
By the fact that $f_{I, 0}, \theta_{I, 0}\rightarrow 0$ as $\tau\rightarrow\infty$, we can obtain
\begin{equation}\label{daoshuweiling}
\overline {f_{I, 0}}+\theta_{I, 0}\equiv 0,
\end{equation}
which further implies
\begin{equation}
\overline {f_{0}}(0, \vec x)+\theta_{0}(0, \vec x)=\overline h+\theta^{0}:=l_{0}(\vec x).
\end{equation}
According to \eqref{fk} by taking $k=0$ and the fact that $\theta_{0}(0, \vec x)=\theta^{0}_{0}(\vec x)$, we can derive the following equality
\begin{equation}\label{biaoshi}
(\theta^{0}_{0})^{4}+\theta_{0}^{0}=l_{0}(\vec x).
\end{equation}

Define $G(\theta_{0}^{0})=(\theta^{0}_{0})^{4}+\theta_{0}^{0}:\mathbb{R}^{+}\rightarrow\mathbb{R}^{+}$. Since $\rho^{0}\geq a>0$, we have $G'=4(\theta_{0}^{0})^{3}+1>0$ for any $\theta_{0}^{0}\in\mathbb{R}^{+}$. So, $G$ is a one-to-one function. By the fact that $l_{0}\geq a^{2}>0$ and $G$ is a $C^\infty$ function about $\theta_{0}^{0}$, then the inverse function $\theta_{0}^{0}(\vec x)=G^{-1}(l_{0}(\vec x))\geq b>0$ where $b=b(a)$, $G^{-1}(l_{0}(\vec x))$ is a $C^\infty$ function about $l_{0}$. Since $l_{0}(\vec x)\in H_{\vec x}^{N+2}$, we have $\theta_{0}^{0}(\vec x)=G^{-1}(l_{0}(\vec x))\in H_{\vec x}^{N+2}$. Then, \eqref{ineq-theta00-jia0} follows directly with some constant $C$ depending on $\overline h$ and $\theta^{0}$.

According to \eqref{biaoshi} and the definition of $l_{0}(\vec x)$, we have
\begin{equation}\label{biaoshi1}
(\theta^{0}_{0})^{4}+\theta_{0}^{0}=\overline h+\theta^{0},
\end{equation}
which can be further written as 
\begin{equation}\label{biaoshi2}
(\theta^{0}_{0})^{4}-(\theta^{0})^{4}+\theta_{0}^{0}-\theta^{0}=\overline h-(\theta^{0})^{4}.
\end{equation}
Then,
\begin{equation}\label{biaoshi3}
\theta_{0}^{0}-\theta^{0}=\frac{\overline h-(\theta^{0})^{4}}{((\theta_{0}^{0})^{2}+(\theta^{0})^{2})(\theta_{0}^{0}+\theta^{0})+1}.
\end{equation}
So,
\begin{equation}\label{biaoshi4}
\|\theta_{0}^{0}-\theta^{0}\|_{H^{N+2}(\mathbb{T}^{3})}\leq C\eta.
\end{equation}
Furthermore, we have
\begin{equation}\label{biaoshi5}
\|h-(\theta_{0}^{0})^{4}\|_{H^{N+2}(\mathbb{T}^{3})}\leq \|h-(\theta^{0})^{4}\|_{H^{N+2}(\mathbb{T}^{3})}+\|(\theta_{0}^{0})^{4}-(\theta^{0})^{4}\|_{H^{N+2}(\mathbb{T}^{3})}\leq C\eta.
\end{equation}

According to the equality \eqref{daoshuweiling}, we can write $\eqref{equations about f I0,theta I0}_{2}$ in the following form
\begin{equation}
\frac{\partial \theta_{I, 0}}{\partial\tau}=(-\theta_{I, 0}-B(\theta_{0}^{0}+\theta_{I, 0})+B(\theta_{0}^{0})),
\end{equation}
which can be further rewritten as
\begin{equation}\label{wenduini}
\frac{\partial \theta_{I, 0}}{\partial\tau}+(1+4(\theta_{0}^{0})^{3})\theta_{I, 0}=-6(\theta_{0}^{0})^{2}\theta_{I, 0}^{2}-4\theta_{0}^{0}\theta_{I, 0}^{3}.
\end{equation}
Then, we have
\begin{equation}\label{jifenfangcheng0}
\begin{aligned}
\theta_{I, 0}(\tau)=e^{-(1+4(\theta_{0}^{0})^{3})\tau}(\theta^{0}-\theta_{0}^{0})-\int_{0}^{\tau}e^{-(1+4(\theta_{0}^{0})^{3})(\tau-s)}(6(\theta_{0}^{0})^{2}\theta_{I, 0}^{2}+4\theta_{0}^{0}\theta_{I, 0}^{3})(s)\mathrm{d}s.
\end{aligned}
\end{equation}

We define the energy radius $E_{13}=\{\theta_{I, 0}:\|\theta_{I, 0}\|_{L^{\infty}(0, \infty; H_{\vec x}^{N+2})}\leq \tilde\eta\}$, where $\tilde\eta=2\|e^{-(1+4(\theta_{0}^{0})^{3})\tau}(\theta^{0}-\theta_{0}^{0})\|_{L^{\infty}(0, \infty; H_{\vec x}^{N+2})}=O(\eta)$.

We give the linearized form of \eqref{jifenfangcheng0} as follows:
\begin{equation}\label{jifenfangcheng00}
\begin{aligned}
\theta_{I, 0}^{k+1}(\tau)=e^{-(1+4(\theta_{0}^{0})^{3})\tau}(\theta^{0}-\theta_{0}^{0})
-\int_{0}^{\tau}e^{-(1+4(\theta_{0}^{0})^{3})(\tau-s)}(6(\theta_{0}^{0})^{2}(\theta_{I, 0}^{k})^{2}+4\theta_{0}^{0}(\theta_{I, 0}^{k})^{3})(s)\mathrm{d}s,
\end{aligned}
\end{equation}
where we choose $\theta_{I, 0}^{0}=0$. Assume $\|\theta_{I, 0}^{k}\|_{L^{\infty}(0, \infty; H_{\vec x}^{N+2})}\leq\tilde\eta,$ then,
\begin{equation}\label{eneinitialo}
\begin{aligned}
\|\theta_{I, 0}^{k+1}\|_{L^{\infty}(0, \infty; H_{\vec x}^{N+2})}\leq\tilde\eta+C_{1}\tilde\eta^{2}+C_{2}\tilde\eta^{3}\leq 2\tilde\eta,
\end{aligned}
\end{equation}
for sufficiently small $\eta$.

Furthermore, we set $\tilde{\theta_{I, 0}^{k+1}}=\theta_{I, 0}^{k+1}-\theta_{I, 0}^{k}$. Then, we can derive
\begin{equation}\label{jifenfangcheng1}
\begin{aligned}
\tilde{\theta_{I, 0}^{k+1}}(\tau)=&-\int_{0}^{\tau}e^{-(1+4(\theta_{0}^{0})^{3})(\tau-s)}(6(\theta_{0}^{0})^{2}(\theta_{I, 0}^{k}+\theta_{I, 0}^{k-1})\tilde{\theta_{I, 0}^{k}}+4\theta_{0}^{0}(\theta_{I, 0}^{k}+\theta_{I, 0}^{k-1})\\&\times((\theta_{I, 0}^{k})^{2}+\theta_{I, 0}^{k}\theta_{I, 0}^{k-1}+(\theta_{I, 0}^{k-1})^{2})\tilde{\theta_{I, 0}^{k}}(s)\mathrm{d}s,
\end{aligned}
\end{equation}
which further implies
\begin{equation}\label{contraction-ini}
\begin{aligned}
\|\tilde{\theta_{I, 0}^{k+1}}\|_{L^{\infty}(0, \infty; H_{\vec x}^{N+2})}\leq (C_{1}'\tilde\eta+C_{2}'\tilde\eta^{2})\|\tilde{\theta_{I, 0}^{k}}\|_{L^{\infty}(0, \infty; H_{\vec x}^{N+2})}.
\end{aligned}
\end{equation}
For sufficiently small $\eta$, we have $0<(C_{1}'\tilde\eta+C_{2}'\tilde\eta^{2})<1$. Then, we can claim that $\{\theta_{I, 0}^{k}\}$ is a contraction sequence.

According to \eqref{eneinitialo} and \eqref{contraction-ini}, we have $\theta_{I, 0}^{k}\rightarrow \theta_{I, 0}$ in $E_{13}$ as $k\rightarrow\infty$, by the Banach fixed point theorem.

According to $\eqref{equations about f I0,theta I0}_{1}$ and $\eqref{equations about f I0,theta I0}_{3}$, we can write $f_{I, 0}$ in the following form
\begin{equation}\label{fI0biaodashi}
f_{I, 0}=e^{-\tau}(h-(\theta_{0}^{0})^{4})+\int_{0}^{\tau}e^{s-\tau}(B(\theta_{0}^{0}+\theta_{I, 0})-B(\theta_{0}^{0}))\mathrm{d}s,
\end{equation}
which further implies the existence solution $f_{I, 0}\in C^{1}([0, \infty); L^{\infty}_{\vec w}H_{\vec x}^{N+2})$. 

Furthermore, we have
\begin{equation}
\|f_{I, 0}\|_{L^{\infty}(0, \infty; H_{\vec x}^{N+2})}\leq C\eta.
\end{equation}

In the following, let us prove that $(f_{I, 0}, \theta_{I, 0})\in L^{1}([0, \infty); L^{\infty}_{\vec w}H_{\vec x}^{N+2})\times L^{1}([0, \infty); H_{\vec x}^{N+2})$.

According to \eqref{fI0biaodashi}, we have 
\begin{equation}\label{gron1}
\begin{aligned}
\|f_{I, 0}(\tau)\|_{L^{\infty}_{\vec w}H_{\vec x}^{N+2}}\leq& e^{-\tau}\|h-(\theta_{0}^{0})^{4}\|_{L^{\infty}_{\vec w}H_{\vec x}^{N+2}}+C\eta\int_{0}^{\tau}\|\theta_{I, 0}\|_{H_{\vec x}^{N+2}}\mathrm{d}s
\\\leq& C\eta e^{-\tau}+C\eta\int_{0}^{\tau}\|\theta_{I, 0}\|_{H_{\vec x}^{N+2}}\mathrm{d}s.
\end{aligned}
\end{equation}

According to \eqref{jifenfangcheng0}, we have
\begin{equation}\label{jifenfangcheng0}
\begin{aligned}
\theta_{I, 0}(\tau)=e^{-(1+4(\theta_{0}^{0})^{3})\tau}(\theta^{0}-\theta_{0}^{0})
-\int_{0}^{\tau}e^{-(1+4(\theta_{0}^{0})^{3})(\tau-s)}(6(\theta_{0}^{0})^{2}\theta_{I, 0}^{2}+4\theta_{0}^{0}\theta_{I, 0}^{3})(s)\mathrm{d}s,
\end{aligned}
\end{equation}

So,
\begin{equation}\label{gron2}
\begin{aligned}
\|\theta_{I, 0}(\tau)\|_{H_{\vec x}^{N+2}}\leq C\eta e^{-\kappa\tau}+C\eta\int_{0}^{\tau}e^{-\kappa(\tau-s)}\|\theta_{I, 0}\|_{H_{\vec x}^{N+2}}\mathrm{d}s,
\end{aligned}
\end{equation}
where $\kappa\in (0, 1)$ is a suitably small positive constant.

Add \eqref{gron1} and \eqref{gron2} together, we have
\begin{equation}\label{gron3}
\begin{aligned}
&\|f_{I, 0}(\tau)\|_{L^{\infty}_{\vec w}H_{\vec x}^{N+2}}+\|\theta_{I, 0}(\tau)\|_{H_{\vec x}^{N+2}}\\\leq& C\eta e^{-\kappa\tau}+C\eta\int_{0}^{\tau}e^{-\kappa(\tau-s)}(\|f_{I, 0}\|_{L^{\infty}_{\vec w}H_{\vec x}^{N+2}}+\|\theta_{I, 0}\|_{H_{\vec x}^{N+2}})\mathrm{d}s,
\end{aligned}
\end{equation}
which further implies 
\begin{equation}\label{gron4}
\begin{aligned}
&\|f_{I, 0}\|_{L^{1}([0, \infty); L^{\infty}_{\vec w}H_{\vec x}^{N+2})}+\|\theta_{I, 0}(\tau)\|_{L^{1}([0, \infty); H_{\vec x}^{N+2})}\leq C\eta,
\end{aligned}
\end{equation}
combining with the Gronwall's inequality.

Finally, we will prove the exponential decay of the above solutions.

For sufficiently small $\sigma>0$, multiplying $e^{\sigma\tau}$ on both sides of \eqref{equations about f I0,theta I0}, we have
\begin{equation}\label{equations about f I0,theta I0-zhishu}\left\{
\begin{split}
&\frac{\partial (e^{\sigma\tau}f_{I, 0})}{\partial\tau}+(1-\sigma)e^{\sigma\tau}f_{I, 0}=e^{\sigma\tau}(B(\theta_{0}^{0}+\theta_{I, 0})-B(\theta_{0}^{0})),\\
&\partial_{\tau}(e^{\sigma\tau}\theta_{I, 0})+\Big(4(\theta_{0}^{0})^{3}-\sigma\Big)e^{\sigma\tau}\theta_{I, 0}=e^{\sigma\tau}(\overline{f_{I, 0}}-6(\theta_{0}^{0})^{2}\theta_{I, 0}^{2}-4\theta_{0}^{0}\theta_{I, 0}^{3}-\theta_{I, 0}^{4}),\\
&f_{I, 0}(0, \vec{x}, \vec{w})=h(\vec{x}, \vec{w})-\overline {f_{0}}(0, \vec x),\\
&\theta_{I, 0}(0, \vec{x}, \vec{w})=\theta^{0}(\vec x)-\theta_{0}^{0}(\vec x),\\
&\lim_{\tau\rightarrow\infty}f_{I, 0}(\tau, \vec{x}, \vec{w})=0,\ \ \lim_{\tau\rightarrow\infty}\theta_{I, 0}(\tau, \vec{x}, \vec{w})=0.
\end{split}\right.
\end{equation}

Then, we can write 
\begin{equation}\label{fI0biaodashi-zhishu}
e^{\sigma\tau}f_{I, 0}=e^{-(1-\sigma)\tau}(h-(\theta_{0}^{0})^{4})+\int_{0}^{\tau}e^{-(1-\sigma)(\tau-s)}(B(\theta_{0}^{0}+\theta_{I, 0})-B(\theta_{0}^{0}))\mathrm{d}s
\end{equation}
and
\begin{equation}
\begin{aligned}
e^{\sigma\tau}\theta_{I, 0}(\tau)=&exp\Big\{-\Big(1+4(\theta_{0}^{0})^{3}-\sigma\Big)\tau\Big\}(\theta^{0}-\theta^{0}_{0})
\\&+\int_{0}^{\tau}exp\Big\{-\Big(1+4(\theta_{0}^{0})^{3}-\sigma\Big)(\tau-s)\Big\}(-6(\theta_{0}^{0})^{2}\theta_{I, 0}^{2}-4\theta_{0}^{0}\theta_{I, 0}^{3}-\theta_{I, 0}^{4})\mathrm{d}s.
\end{aligned}
\end{equation}
For suitably small $\sigma$, we can derive
\begin{equation}
\|e^{\sigma\tau}f_{I, 0}\|_{L^{\infty}([0, \infty); L^{\infty}_{\vec w}H_{\vec x}^{N+2})}+\|e^{\sigma\tau}\theta_{I, 0}\|_{L^{\infty}([0, \infty); H_{\vec x}^{N+2})}\leq C\eta
\end{equation}
and
\begin{equation}
\|e^{\sigma\tau}f_{I, 0}\|_{L^{1}([0, \infty); L^{\infty}_{\vec w}H_{\vec x}^{N+2})}+\|e^{\sigma\tau}\theta_{I, 0}\|_{L^{1}([0, \infty); H_{\vec x}^{N+2})}\leq C\eta
\end{equation}
by using similar calculations as above and Gronwall's inequality.

Then, collecting \eqref{thetak}, \eqref{fk} by taking $k=0$ and the fact that $B(\theta_{0})=\theta_{0}^{4}$, we have
\begin{equation}\label{equations about f0 epsilon and theta}\left\{
\begin{split}
&f_{0}=\overline {f_{0}}=\theta_{0}^{4},\\
&\partial_{t}(\theta_{0}+\theta_{0}^{4})+\mathrm{div}(\vec{u}\theta_{0})-\Delta_{x}\Big(\theta_{0}+\frac{1}{3}\theta_{0}^{4}\Big)=0,~~\mathrm{in} \ \ (0,T)\times\mathbb{T}^{3},\\
&\theta_{0}(0, \vec{x})=\theta_{0}^{0}(\vec{x}),\;\;\overline {f_{0}}(0, \vec{x})=(\theta_{0}^{0})^{4}~~\mathrm{in}  \ \ \mathbb{T}^{3}.
\end{split}\right.
\end{equation}
\begin{thm}\label{3.2}
For the given data $\theta_{0}^{0}(\vec{x})\in H_{\vec x}^{N+2}$ and $\theta_{0}^{0}\geq b$ where $b$ is a given positive constant, the problem \eqref{equations about f0 epsilon and theta} has a unique local solution on $[0, T]$, and $(\theta_{0}, f_{0})\in C^{0}([0, T]; H_{\vec x}^{N+2})\cap L^{2}(0, T; H_{\vec x}^{N+3})\times C^{0}([0, T]; H_{\vec x}^{N+2})\cap L^{2}(0, T; H_{\vec x}^{N+3})$. Furthermore, we have
\begin{equation}\label{ineq-theta00}
\|(\theta_{0}, f_{0})\|_{L_{T}^{\infty}H_{\vec x}^{N+2}}+\|(\theta_{0}, f_{0})\|_{L^{2}_{T}H_{\vec x}^{N+3}}\leq C
\end{equation}
and $\rho_{0}, \theta_{0}\geq b/2$.
\end{thm}
\noindent\textbf{Proof.}We can write \eqref{equations about f0 epsilon and theta} in the following form
\begin{equation}\label{equations about f0 epsilon and theta-1}\left\{
\begin{split}
&f_{0}=\overline {f_{0}}=\theta_{0}^{4},\\
&\partial_{t}\theta_{0}-\frac{4\theta_{0}^{3}}{(1+4\theta_{0}^{3})}\Delta_{x}\theta_{0}=-\frac{1}{(1
+4\theta_{0}^{3})}\vec{u}_{0}\cdot\nabla_{x}\theta_{0}-\frac{\theta_{0}}{(1+4\theta_{0}^{3})}\mathrm{div}_{x}\vec{u}_{0}
\\&+\frac{4\theta_{0}^{2}}{3(1+4\theta_{0}^{3})}|\nabla_{x}\theta_{0}|^{2}+\frac{4\theta_{0}^{3}}{3(1+4\theta_{0}^{3})}\nabla_{x}\theta_{0}\cdot\vec u_{0},~~\mathrm{in} \ \ (0,T)\times\mathbb{T}^{3}, \\
&\theta_{0}(0, \vec{x})=\theta_{0}^{0}(\vec{x}),\;\;\overline {f_{0}}(0, \vec{x})=(\theta_{0}^{0})^{4}~~\mathrm{in}  \ \ \mathbb{T}^{3}.
\end{split}\right.
\end{equation}
We shall study the Cauchy problem for $\eqref{equations about f0 epsilon and theta-1}$.

In the sequel, we construct a solution to $\eqref{equations about f0 epsilon and theta-1}$. For $k\geq 0$ and $k\in\mathbb{Z}$, we define $\theta^{k+1}$ inductively as the solution of the following linearized system:
\begin{equation}\label{research equations-4}\left\{
\begin{split}
&\partial_{t}\theta^{k+1}-\frac{4(\theta^{k})^{3}}{(1+4(\theta^{k})^{3})}\Delta_{x}\theta^{k+1}=-\frac{1}{(1
+4(\theta^{k})^{3})}\vec{u}\cdot\nabla_{x}\theta^{k}-\frac{\theta^{k}}{(1+4(\theta^{k})^{3})}\mathrm{div}_{x}\vec{u}
\\&+\frac{4(\theta^{k})^{2}}{3(1+4(\theta^{k})^{3})}|\nabla_{x}\theta^{k}|^{2}+\frac{4(\theta^{k})^{3}}{3(1+4(\theta^{k})^{3})}\nabla_{x}\theta^{k}\cdot\vec u~~\mathrm{in} \ \ (0, T^{k+1})\times\mathbb{T}^{3}\times \mathbb{S}^{2},\\
&\theta^{k+1}(0, \vec{x})=\theta_{0}^{0}(\vec{x})~~\mathrm{in} \ \ \mathbb{T}^{3}.\\
\end{split}\right.
\end{equation}
To avoid confusion, we write $\theta^{k}|_{k=0}=\theta^{00}$. 

Define the map
\begin{equation}\nonumber
S: \theta^{k}\in G\rightarrow \theta^{k+1}:=S(\theta^{k})
\end{equation}
where $\theta^{k}$ satisfies (\ref{research equations-4}).

Set
\begin{equation}\nonumber
\overline r_{k}=\|\theta^{k}\|_{L^{\infty}_{T^{k}}H_{\vec x}^{N+2}}.
\end{equation}

First, we note that the Proposition 2.7 in \cite{bu-7-jiajia} can be shown in a similar manner in $\mathbb{T}^{3}$.
Choosing $\theta^{00}=\theta^{0}$ for some $T^{1}\in(0, \infty),$ we get $V^{1}\in C^{0}([0, T^{1}]; H_{\vec x}^{N+2})$ and $\theta^{1}\in C^{0}([0, T^{1}]; H_{\vec x}^{N+2})\cap L^{2}(0, T^{1}; H_{\vec x}^{N+3})$ by the Proposition 2.7 in \cite{bu-7-jiajia}. Note that $T^{1}$ is the largest time of existence such that $\theta^{1}\geq M_{1}>0$ and $\overline r_{1}\leq r^{0}$, where $r^{0}$ is to be determined. Then, we can derive that $\theta^{2}\in C^{0}([0, T^{2}]; H_{\vec x}^{N+2})\cap L^{2}(0, T^{2}; H_{\vec x}^{N+3})$. We also note that $T^{2}$ is the largest time of existence such that $\rho^{2}, \theta^{2}\geq M_{2}>0$ and $\overline r_{2}\leq r^{0}$. Finally, we obtain that $\theta^{k+1}\in C^{0}([0, T^{k+1}]; H_{\vec x}^{N+2})\cap L^{2}(0, T^{k+1}; H_{\vec x}^{N+3})$ with $\theta^{k+1}\geq M_{k+1}>0$ and $\overline r_{k+1}\leq r^{0}$ by choosing sufficiently small satisfying $T^{k+1}\leq T^{k}$ for any $k\geq 0$ iteratively.

Set
\begin{equation}\nonumber
r_{k}=\|\theta^{k}\|_{L_{T}^{\infty}H_{\vec x}^{N+2}},
\end{equation}
where $T\in [0, T^{k}]$ is to be determined. We will prove that there exist $T$ independent of $k$ such that
\begin{equation}\nonumber
T^{k}\geq T>0,
\end{equation}
\begin{equation}\nonumber
\theta^{k}\geq M_{k}\geq M>0
\end{equation}
and
\begin{equation}\nonumber
r_{k}\leq r^{0}
\end{equation}
where $M$ and $r$ are independent of $k$ for any $k\geq 0$. Note that $r_{k}, M_{k}, \theta^{k}$ will disappear if $k<0$.

Applying the operator $D_{x}^{\gamma}$ to the equation $(\ref{research equations-4})_{2}$, multiplying the resulting equation with $2(D_{x}^{\gamma}\theta^{k+1})$ in $L^{2}((0, t)\times\mathbb{T}^{3})$, we have
\begin{equation}\label{Eq 1}
\begin{split}
&\int_{\mathbb{T}^{3}}|D_{x}^{\gamma}\theta^{k+1}|^{2}(t)\mathrm{d}\vec{x}-\int_{\mathbb{T}^{3}}|D_{x}^{\gamma}\theta^{0}_{0}|^{2}(t)\mathrm{d}\vec{x}
\\=&-\int_{0}^{t}\int_{\mathbb{T}^{3}}\frac{8(\theta^{k})^{3}}{(1+4(\theta^{k})^{3})}|D_{x}^{\gamma}\nabla_{x}\theta^{k+1}|^{2}\mathrm{d}\vec{x}\mathrm{d}s
+\int_{0}^{t}\int_{\mathbb{T}^{3}}\Big(\Big[D_{x}^{\gamma}, \frac{8(\theta^{k})^{3}}{(1+4(\theta^{k})^{3})}\Big]\Delta_{x}\theta^{k+1}\Big)D_{x}^{\gamma}\theta^{k+1}\mathrm{d}\vec{x}\mathrm{d}s
\\&+\int_{0}^{t}\int_{\mathbb{T}^{3}}D_{x}^{\gamma}\Big(-\frac{1}{(1
+4(\theta^{k})^{3})}\vec{u}\cdot\nabla_{x}\theta^{k}\Big)D_{x}^{\gamma}\theta^{k+1}\mathrm{d}\vec{x}\mathrm{d}s
+\int_{0}^{t}\int_{\mathbb{T}^{3}}D_{x}^{\gamma}\Big(-\frac{\theta^{k}}{(1+4(\theta^{k})^{3})}\mathrm{div}_{x}\vec{u}^{k}\Big)\\&\times D_{x}^{\gamma}\theta^{k+1}\mathrm{d}\vec{x}\mathrm{d}s
+\int_{0}^{t}\int_{\mathbb{T}^{3}}D_{x}^{\gamma}\Big(\frac{4(\theta^{k})^{2}}{3(1+4(\theta^{k})^{3})}|\nabla_{x}\theta^{k}|^{2}\Big)
D_{x}^{\gamma}\theta^{k+1}\mathrm{d}\vec{x}\mathrm{d}s
\\&+\int_{0}^{t}\int_{\mathbb{T}^{3}}D_{x}^{\gamma}\Big(\frac{4(\theta^{k})^{3}}{3(1+4(\theta^{k})^{3})}\nabla_{x}\theta^{k}\cdot\vec u\Big)D_{x}^{\gamma}\theta^{k+1}\mathrm{d}\vec{x}\mathrm{d}s:=D_{1}+D_{2}+D_{3}+D_{4}+D_{5}+D_{6},
\end{split}
\end{equation}
where $0\leq\gamma\leq N+2$ and $t\in [0, T]$. We also note that $[\cdot, \cdot]$ is the commutator operator, i.e.,
\begin{equation}\label{commutator}
[D_{x}^{\gamma}, f]g=\sum_{\alpha+\beta=\gamma, |\alpha|\geq 1}C_{\alpha, \beta}D_{x}^{\alpha}f\cdot D_{x}^{\beta}g.
\end{equation}

Then, we have
\begin{equation}\label{D1}
\begin{split}
D_{1}=&-\int_{0}^{t}\int_{\mathbb{T}^{3}}\frac{8(\theta^{k})^{3}}{(1+4(\theta^{k})^{3})}|D_{x}^{\gamma}\nabla_{x}\theta^{k+1}|^{2}
\mathrm{d}\vec{x}\mathrm{d}s
\\\leq&-\frac{M_{k}^{3}}{C(r_{k}^{3}+1)}\int_{0}^{t}\int_{\mathbb{T}^{3}}|D_{x}^{\gamma+1}\theta^{k+1}|^{2}
\mathrm{d}\vec{x}\mathrm{d}s
\end{split}
\end{equation}
and
\begin{equation}\label{D2}
\begin{split}
D_{2}=&\int_{0}^{t}\int_{\mathbb{T}^{3}}\Big(\Big[D_{x}^{\gamma}, \frac{8(\theta^{k})^{3}}{(1+4(\theta^{k})^{3})}\Big]\Delta_{x}\theta^{k+1}\Big)D_{x}^{\gamma}\theta^{k+1}\mathrm{d}\vec{x}\mathrm{d}s
\\\leq&\int_{0}^{t}\Big\|\Big(\frac{8(\theta^{k})^{3}}{(1+4(\theta^{k})^{3})}\Big)(s)\Big\|_{H_{\vec x}^{N+2}}
\|D_{x}\theta^{k+1}(s)\|_{H_{\vec x}^{N+2}}\|\theta^{k+1}(s)\|_{H_{\vec x}^{N+2}}\mathrm{d}s
\\\leq& C(1+r_{k}^{l})\Big(1+\frac{1}{M_{k}^{l}}\Big)\|D_{x}\theta^{k+1}(s)\|_{L^{2}_{t}H_{\vec x}^{N+2}}
\|\theta^{k+1}(s)\|_{L^{1}_{t}H_{\vec x}^{N+2}}
\\\leq& C(1+r_{k}^{2l})\Big(1+\frac{1}{M_{k}^{2l}}\Big)\delta^{-1}\|\theta^{k+1}(s)\|_{L^{2}_{t}H_{\vec x}^{N+2}}^{2}+
\delta\|D_{x}\theta^{k+1}(s)\|_{L^{1}_{t}H_{\vec x}^{N+2}}^{2}
\\\leq& C(1+r_{k}^{2l})\Big(1+\frac{1}{M_{k}^{2l}}\Big)\delta^{-1}\|\theta^{k+1}(s)\|_{L^{\infty}_{t}H_{\vec x}^{N+2}}^{2}t+
\delta\|D_{x}\theta^{k+1}(s)\|_{L^{2}_{t}H_{\vec x}^{N+2}}^{2},
\end{split}
\end{equation}
where $l$ is suitably large.

We can write $D_{3}$ in the following form
\begin{equation}\label{D3}
\begin{split}
D_{3}=&\int_{0}^{t}\int_{\mathbb{T}^{3}}D_{x}^{\gamma}\Big(-\frac{1}{(1+4(\theta^{k})^{3})}\vec u\cdot\nabla_{x}\theta^{k}\Big)D_{x}^{\gamma}\theta^{k+1}\mathrm{d}\vec{x}\mathrm{d}s
\\=&\int_{0}^{t}\int_{\mathbb{T}^{3}}\Big(-\frac{1}{(1+4(\theta^{k})^{3})}\vec u\cdot D_{x}^{\gamma}\nabla_{x}\theta^{k}\Big)D_{x}^{\gamma}\theta^{k+1}\mathrm{d}\vec{x}\mathrm{d}s\\&+\int_{0}^{t}\int_{\mathbb{T}^{3}}
\Big(\Big[D_{x}^{\gamma}, -\frac{1}{(1+4(\theta^{k})^{3})}\vec u\cdot\Big] \nabla_{x}\theta^{k}\Big)D_{x}^{\gamma}\theta^{k+1}\mathrm{d}\vec{x}\mathrm{d}s:=D_{31}+D_{32}.
\end{split}
\end{equation}
By integrating by parts, we write $D_{31}$ into the following two parts
\begin{equation}\label{D31}
\begin{split}
D_{31}=&\int_{0}^{t}\int_{\mathbb{T}^{3}}\Big(-\frac{1}{(1+4(\theta^{k})^{3})}\vec u\cdot D_{x}^{\gamma}\nabla_{x}\theta^{k}\Big)D_{x}^{\gamma}\theta^{k+1}\mathrm{d}\vec{x}\mathrm{d}s
\\=&\int_{0}^{t}\int_{\mathbb{T}^{3}}D_{x}\Big(-\frac{1}{(1+4(\theta^{k})^{3})}\vec u\cdot\Big) (D_{x}^{\gamma-1}\nabla_{x}\theta^{k})D_{x}^{\gamma}\theta^{k+1}\mathrm{d}\vec{x}\mathrm{d}s\\&+\int_{0}^{t}\int_{\mathbb{T}^{3}}
\Big(-\frac{1}{(1+4(\theta^{k})^{3})}\vec u\cdot D_{x}^{\gamma-1}\nabla_{x}\theta^{k}\Big)D_{x}^{\gamma+1}\theta^{k+1}\mathrm{d}\vec{x}\mathrm{d}s:=D_{311}+D_{312},
\end{split}
\end{equation}
where
\begin{equation}\label{D311}
\begin{split}
D_{311}=&\int_{0}^{t}\int_{\mathbb{T}^{3}}D_{x}\Big(-\frac{1}{(1+4(\theta^{k})^{3})}\vec u\cdot\Big) (D_{x}^{\gamma-1}\nabla_{x}\theta^{k})D_{x}^{\gamma}\theta^{k+1}\mathrm{d}\vec{x}\mathrm{d}s
\\\leq& \Big\|D_{x}\Big(-\frac{1}{(1+4(\theta^{k})^{3})}\vec u\Big)\Big\|_{L_{t}^{\infty}L_{\vec x}^{\infty}}\|\theta^{k}\|_{L^{\infty}_{t}H_{\vec x}^{N+2}}\|\theta^{k+1}\|_{L^{\infty}_{t}H_{\vec x}^{N+2}}t
\\\leq& C\Big\|D_{x}\Big(-\frac{1}{(1+4(\theta^{k})^{3})}\vec u\Big)\Big\|_{L_{t}^{\infty}L_{\vec x}^{\infty}}^{2}\|\theta^{k}\|_{L^{\infty}_{t}H_{\vec x}^{N+2}}^{2}t+\|\theta^{k+1}\|_{L^{\infty}_{t}H_{\vec x}^{N+2}}^{2}t
\\\leq& C\Big(1+\frac{1}{M_{k}^{l}}\Big)(1+r_{k}^{l})t+\|\theta^{k+1}\|_{L^{\infty}_{t}H_{\vec x}^{N+2}}^{2}t
\end{split}
\end{equation}
and
\begin{equation}\label{D312}
\begin{split}
D_{312}=&\int_{0}^{t}\int_{\mathbb{T}^{3}}\Big(-\frac{1}{(1+4(\theta^{k})^{3})}\vec u\cdot D_{x}^{\gamma-1}\nabla_{x}\theta^{k}\Big)D_{x}^{\gamma+1}\theta^{k+1}\mathrm{d}\vec{x}\mathrm{d}s
\\\leq& \Big\|-\frac{1}{(1+4(\theta^{k})^{3})}\vec u\Big\|_{L_{t}^{\infty}L_{\vec x}^{\infty}}\|\theta^{k}\|_{L^{\infty}_{t}H_{\vec x}^{N+2}}\sqrt{t}\|D_{x}\theta^{k+1}\|_{L^{2}_{t}H_{\vec x}^{N+2}}
\\\leq& C\delta^{-1}\Big\|-\frac{1}{(1+4(\theta^{k})^{3})}\vec u\Big\|_{L_{t}^{\infty}L_{\vec x}^{\infty}}^{2}\|\theta^{k}\|_{L^{\infty}_{t}H_{\vec x}^{N+2}}^{2}t+\delta\|D_{x}\theta^{k+1}\|_{L^{2}_{t}H_{\vec x}^{N+2}}^{2}
\\\leq& C\delta^{-1}\Big(1+\frac{1}{M_{k}^{l}}\Big)(1+r_{k}^{l})t+\delta\|D_{x}\theta^{k+1}\|_{L^{2}_{t}H_{\vec x}^{N+2}}^{2}.
\end{split}
\end{equation}
Finally, 
\begin{equation}\label{D32}
\begin{split}
D_{32}=&\int_{0}^{t}\int_{\mathbb{T}^{3}}
\Big(\Big[D_{x}^{\gamma}, -\frac{1}{(1+4(\theta^{k})^{3})}\vec u\cdot\Big] \nabla_{x}\theta^{k}\Big)D_{x}^{\gamma}\theta^{k+1}\mathrm{d}\vec{x}\mathrm{d}s
\\\leq& \Big\|-\frac{1}{(1+4(\theta^{k})^{3})}\vec u\Big\|_{L^{\infty}_{t}H_{\vec x}^{N+2}}\|\theta^{k}\|_{L^{\infty}_{t}H_{\vec x}^{N+2}}\|\theta^{k+1}\|_{L^{\infty}_{t}H_{\vec x}^{N+2}}t
\\\leq& C\Big\|-\frac{1}{(1+4(\theta^{k})^{3})}\vec u\Big\|_{L^{\infty}_{t}H_{\vec x}^{N+2}}^{2}\|\theta^{k}\|_{L^{\infty}_{t}H_{\vec x}^{N+2}}^{2}t+\|\theta^{k+1}\|_{L^{\infty}_{t}H_{\vec x}^{N+2}}^{2}t
\\\leq& C\Big(1+\frac{1}{M_{k}^{l}}\Big)(1+r_{k}^{l})t+\|\theta^{k+1}\|_{L^{\infty}_{t}H_{\vec x}^{N+2}}^{2}t.
\end{split}
\end{equation}
Then, we can derive the following estimate
\begin{equation}\label{D3-jia}
\begin{split}
D_{3}\leq C(1+\delta^{-1})\Big(1+\frac{1}{M_{k}^{l}}\Big)(1+r_{k}^{l})t+C\|\theta^{k+1}\|_{L^{\infty}_{t}H_{\vec x}^{N+2}}^{2}t+\delta\|D_{x}\theta^{k+1}\|_{L^{2}_{t}H_{\vec x}^{N+2}}^{2}.
\end{split}
\end{equation}
We can write $D_{4}$ in the following form
\begin{equation}\label{D4}
\begin{split}
D_{4}=&\int_{0}^{t}\int_{\mathbb{T}^{3}}D_{x}^{\gamma}\Big(-\frac{\theta^{k}}{(1+4(\theta^{k})^{3})}\mathrm{div}_{x}\vec{u}\Big) D_{x}^{\gamma}\theta^{k+1}\mathrm{d}\vec{x}\mathrm{d}s
\\=&\int_{0}^{t}\int_{\mathbb{T}^{3}}\Big(-\frac{\theta^{k}}{(1+4(\theta^{k})^{3})}D_{x}^{\gamma}\mathrm{div}_{x}\vec{u}\Big) D_{x}^{\gamma}\theta^{k+1}\mathrm{d}\vec{x}\mathrm{d}s\\&+\int_{0}^{t}\int_{\mathbb{T}^{3}}\Big(\Big[D_{x}^{\gamma}, -\frac{\theta^{k}}{(1+4(\theta^{k})^{3})}\Big]\mathrm{div}_{x}\vec{u}\Big) D_{x}^{\gamma}\theta^{k+1}\mathrm{d}\vec{x}\mathrm{d}s:=D_{41}+D_{42}.
\end{split}
\end{equation}
Furthermore, we can write $D_{41}$ in the following form
\begin{equation}\label{D41}
\begin{split}
D_{41}=&\int_{0}^{t}\int_{\mathbb{T}^{3}}\Big(-\frac{\theta^{k}}{(1+4(\theta^{k})^{3})}
D_{x}^{\gamma}\mathrm{div}_{x}\vec{u}\Big) D_{x}^{\gamma}\theta^{k+1}\mathrm{d}\vec{x}\mathrm{d}s
\\=&\int_{0}^{t}\int_{\mathbb{T}^{3}}D_{x}\Big(\frac{\theta^{k}}{(1+4(\theta^{k})^{3})}\Big)
(D_{x}^{\gamma-1}\mathrm{div}_{x}\vec{u}) D_{x}^{\gamma}\theta^{k+1}\mathrm{d}\vec{x}\mathrm{d}s
\\&+\int_{0}^{t}\int_{\mathbb{T}^{3}}\Big(\frac{\theta^{k}}{(1+4(\theta^{k})^{3})}
(D_{x}^{\gamma-1}\mathrm{div}_{x}\vec{u})\Big) D_{x}^{\gamma+1}\theta^{k+1}\mathrm{d}\vec{x}\mathrm{d}s:=D_{411}+D_{412},
\end{split}
\end{equation}
by integrating by parts.

For $D_{411}$, we can derive
\begin{equation}\label{D411}
\begin{aligned}
D_{411}=&\int_{0}^{t}\int_{\mathbb{T}^{3}}D_{x}\Big(\frac{\theta^{k}}{(1+4(\theta^{k})^{3})}\Big)
(D_{x}^{\gamma-1}\mathrm{div}_{x}\vec{u})D_{x}^{\gamma}\theta^{k+1}\mathrm{d}\vec{x}\mathrm{d}s
\\\leq& C\int_{0}^{t}\Big\|D_{x}\Big(\frac{\theta^{k}}{(1+4(\theta^{k})^{3})}\Big)\Big\|_{L^{\infty}_{\vec x}}
\|D_{x}^{\gamma-1}\mathrm{div}_{x}\vec{u}\|_{L^{2}_{\vec x}}\|D_{x}^{\gamma}\theta^{k+1}\|_{L^{2}_{\vec x}}\mathrm{d}s
\\\leq& C\Big(1+\frac{1}{M_{k}^{l}}\Big)(1+r_{k}^{l})\|\vec{u}\|_{L^{2}_{t}H_{\vec x}^{N+2}}
\|\theta^{k+1}\|_{L^{2}_{t}H_{\vec x}^{N+2}}
\end{aligned}
\end{equation}
by the H\"{o}lder inequality and Sobolev embedding theorem.

For $D_{412}$, we can similarly derive
\begin{equation}\label{D412}
\begin{aligned}
D_{412}=&\int_{0}^{t}\int_{\mathbb{T}^{3}}\Big(\frac{\theta^{k}}{(1+4(\theta^{k})^{3})}
(D_{x}^{\gamma-1}\mathrm{div}_{x}\vec{u})\Big) D_{x}^{\gamma+1}\theta^{k+1}\mathrm{d}\vec{x}\mathrm{d}s
\\\leq& C\int_{0}^{t}\Big\|\frac{\theta^{k}}{(1+4(\theta^{k})^{3})}\Big\|_{L^{\infty}_{\vec x}}
\|D_{x}^{\gamma-1}\mathrm{div}_{x}\vec{u}\|_{L^{2}_{\vec x}}\|D_{x}^{\gamma+1}\theta^{k+1}\|_{L^{2}_{\vec x}}\mathrm{d}s
\\\leq& C\Big(1+\frac{1}{M_{k}^{l}}\Big)(1+r_{k}^{l})\|\vec{u}\|_{L^{2}_{t}H_{\vec x}^{N+2}}
\|D_{x}\theta^{k+1}\|_{L^{2}_{t}H_{\vec x}^{N+2}}
\\\leq& C\Big(1+\frac{1}{M_{k}^{2l}}\Big)(1+r_{k}^{2l})\delta^{-1}\|\vec{u}\|_{L^{2}_{t}H_{\vec x}^{N+2}}^{2}+\delta\|D_{x}\theta^{k+1}\|_{L^{2}_{t}H_{\vec x}^{N+2}}^{2}.
\end{aligned}
\end{equation}
By the H\"{o}lder inequality, we obtain
\begin{equation}\label{D42}
\begin{aligned}
D_{42}=&\int_{0}^{t}\int_{\mathbb{T}^{3}}\Big(\Big[D_{x}^{\gamma}, -\frac{\theta^{k}}{(1+4(\theta^{k})^{3})}\Big]\mathrm{div}_{x}\vec{u}\Big) D_{x}^{\gamma}\theta^{k+1}\mathrm{d}\vec{x}\mathrm{d}s
\\\leq& C\int_{0}^{t}\Big\|\frac{\theta^{k}}{(1+4(\theta^{k})^{3})}\Big\|_{H_{\vec x}^{N+2}}
\|\mathrm{div}_{x}\vec{u}\|_{H_{\vec x}^{N+1}}\|\theta^{k+1}\|_{H_{\vec x}^{N+2}}\mathrm{d}s
\\\leq& C\Big(1+\frac{1}{M_{k}^{l}}\Big)(1+r_{k}^{l})\|\vec{u}\|_{L_{t}^{2}H_{\vec x}^{N+2}}\|\theta^{k+1}\|_{L_{t}^{2}H_{\vec x}^{N+2}}.
\end{aligned}
\end{equation}
Collecting \eqref{D4}, \eqref{D41}, \eqref{D411}, \eqref{D412} and \eqref{D42}, we can derive the following estimate
\begin{equation}\label{D4-jia}
\begin{aligned}
D_{4}\leq& C\Big(1+\frac{1}{M_{k}^{l}}\Big)(1+r_{k}^{l})\|\vec{u}\|_{L_{t}^{2}H_{\vec x}^{N+2}}(\|\theta^{k+1}\|_{L_{t}^{2}H_{\vec x}^{N+2}}+\|D_{x}\theta^{k+1}\|_{L_{t}^{2}H_{\vec x}^{N+2}})
\\\leq& C\Big(1+\frac{1}{M_{k}^{2l}}\Big)(1+r_{k}^{2l})(1+\delta^{-1})t+C\|\theta^{k+1}\|_{L_{t}^{\infty}H_{\vec x}^{N+2}}^{2}t+\delta\|D_{x}\theta^{k+1}\|_{L_{t}^{2}H_{\vec x}^{N+2}}^{2}.
\end{aligned}
\end{equation}
We can write $D_{5}$ in the following form
\begin{equation}\label{D5}
\begin{split}
D_{5}=&\int_{0}^{t}\int_{\mathbb{T}^{3}}D_{x}^{\gamma}\Big(\frac{4(\theta^{k})^{2}}{3(1+4(\theta^{k})^{3})}|\nabla_{x}\theta^{k}|^{2}\Big)
D_{x}^{\gamma}\theta^{k+1}\mathrm{d}\vec{x}\mathrm{d}s
\\=&\int_{0}^{t}\int_{\mathbb{T}^{3}}\Big(\frac{4(\theta^{k})^{2}}{3(1+4(\theta^{k})^{3})}D_{x}^{\gamma}|\nabla_{x}\theta^{k}|^{2}\Big)
D_{x}^{\gamma}\theta^{k+1}\mathrm{d}\vec{x}\mathrm{d}s\\&+\int_{0}^{t}\int_{\mathbb{T}^{3}}\Big(\Big[D_{x}^{\gamma}, \frac{4(\theta^{k})^{2}}{3(1+4(\theta^{k})^{3})}\Big]|\nabla_{x}\theta^{k}|^{2}\Big)
D_{x}^{\gamma}\theta^{k+1}\mathrm{d}\vec{x}\mathrm{d}s:=D_{51}+D_{52}.
\end{split}
\end{equation}
By integrating by parts, we write $D_{51}$ into the following two parts
\begin{equation}\label{D51}
\begin{split}
D_{51}=&\int_{0}^{t}\int_{\mathbb{T}^{3}}\Big(\frac{4(\theta^{k})^{2}}{3(1+4(\theta^{k})^{3})}D_{x}^{\gamma}|\nabla_{x}\theta^{k}|^{2}\Big)
D_{x}^{\gamma}\theta^{k+1}\mathrm{d}\vec{x}\mathrm{d}s
\\=&\int_{0}^{t}\int_{\mathbb{T}^{3}}D_{x}\Big(\frac{4(\theta^{k})^{2}}{3(1+4(\theta^{k})^{3})}\Big)D_{x}^{\gamma-1}|\nabla_{x}\theta^{k}|^{2}
D_{x}^{\gamma}\theta^{k+1}\mathrm{d}\vec{x}\mathrm{d}s\\&+\int_{0}^{t}\int_{\mathbb{T}^{3}}\Big(\frac{4(\theta^{k})^{2}}{3(1+4(\theta^{k})^{3})}D_{x}^{\gamma-1}|\nabla_{x}\theta^{k}|^{2}\Big)
D_{x}^{\gamma+1}\theta^{k+1}\mathrm{d}\vec{x}\mathrm{d}s:=D_{511}+D_{512},
\end{split}
\end{equation}
where
\begin{equation}\label{D511}
\begin{split}
D_{511}=&\int_{0}^{t}\int_{\mathbb{T}^{3}}D_{x}\Big(\frac{4(\theta^{k})^{2}}{3(1+4(\theta^{k})^{3})}\Big)D_{x}^{\gamma-1}|\nabla_{x}\theta^{k}|^{2}
D_{x}^{\gamma}\theta^{k+1}\mathrm{d}\vec{x}\mathrm{d}s
\\\leq& C\Big\|D_{x}\Big(\frac{4(\theta^{k})^{2}}{3(1+4(\theta^{k})^{3})}\Big)\Big\|_{L^{\infty}_{t}L^{\infty}_{\vec x}}\|D_{x}^{\gamma-1}|\nabla_{x}\theta^{k}|^{2}\|_{L^{2}_{t}L^{2}_{\vec x}}\|D_{x}^{\gamma}\theta^{k+1}\|_{L^{2}_{t}L^{2}_{\vec x}}
\\\leq& C\Big(1+\frac{1}{M_{k}^{l}}\Big)(1+r_{k}^{l})\|\theta^{k}\|_{L_{t}^{\infty}H_{\vec x}^{N+2}}^{2}\sqrt{t}\|\theta^{k+1}\|_{L_{t}^{\infty}H_{\vec x}^{N+2}}\sqrt{t}
\\\leq& C\Big(1+\frac{1}{M_{k}^{2l}}\Big)(1+r_{k}^{2l+2})t+\|\theta^{k+1}\|_{L_{t}^{\infty}H_{\vec x}^{N+2}}^{2}t
\end{split}
\end{equation}
and
\begin{equation}\label{D512}
\begin{split}
D_{512}=&\int_{0}^{t}\int_{\mathbb{T}^{3}}\Big(\frac{4(\theta^{k})^{2}}{3(1+4(\theta^{k})^{3})}D_{x}^{\gamma-1}|\nabla_{x}\theta^{k}|^{2}\Big)
D_{x}^{\gamma+1}\theta^{k+1}\mathrm{d}\vec{x}\mathrm{d}s
\\\leq& C\Big\|\frac{4(\theta^{k})^{2}}{3(1+4(\theta^{k})^{3})}\Big\|_{L_{t}^{\infty}L_{\vec x}^{\infty}}\|D_{x}^{\gamma-1}|\nabla_{x}\theta^{k}|^{2}\|_{L_{t}^{2}L_{\vec x}^{2}}\|D_{x}^{\gamma+1}\theta^{k+1}\|_{L_{t}^{2}L_{\vec x}^{2}}
\\\leq& C\Big(1+\frac{1}{M_{k}^{l}}\Big)(1+r_{k}^{l})\|\theta^{k}\|_{L_{t}^{\infty}H_{\vec x}^{N+2}}^{2}\sqrt{t}\|D_{x}\theta^{k+1}\|_{L_{t}^{2}H_{\vec x}^{N+2}}
\\\leq& C\delta^{-1}\Big(1+\frac{1}{M_{k}^{2l}}\Big)(1+r_{k}^{2l+2})t+\delta\|D_{x}\theta^{k+1}\|_{L_{t}^{2}H_{\vec x}^{N+2}}^{2}.
\end{split}
\end{equation}
For $D_{52}$, we have
\begin{equation}\label{D52}
\begin{split}
D_{52}=&\int_{0}^{t}\int_{\mathbb{T}^{3}}\Big(\Big[D_{x}^{\gamma}, \frac{4(\theta^{k})^{2}}{3(1+4(\theta^{k})^{3})}\Big]|\nabla_{x}\theta^{k}|^{2}\Big)
D_{x}^{\gamma}\theta^{k+1}\mathrm{d}\vec{x}\mathrm{d}s
\\\leq& C\Big\|\frac{4(\theta^{k})^{2}}{3(1+4(\theta^{k})^{3})}\Big\|_{L_{t}^{\infty}H_{\vec x}^{N+2}}\|\theta^{k}\|_{L_{t}^{\infty}H_{\vec x}^{N+2}}^{2}\|\theta^{k+1}\|_{L_{t}^{\infty}H_{\vec x}^{N+2}}t
\\\leq& C\Big\|\frac{4(\theta^{k})^{2}}{3(1+4(\theta^{k})^{3})}\Big\|_{L_{t}^{\infty}H_{\vec x}^{N+2}}^{2}\|\theta^{k}\|_{L_{t}^{\infty}H_{\vec x}^{N+2}}^{4}t+\|\theta^{k+1}\|_{L_{t}^{\infty}H_{\vec x}^{N+2}}^{2}t
\\\leq& C\Big(1+\frac{1}{M_{k}^{l}}\Big)(1+r_{k}^{l})t+\|\theta^{k+1}\|_{L_{t}^{\infty}H_{\vec x}^{N+2}}^{2}t
\end{split}
\end{equation}
Collecting \eqref{D5}, \eqref{D51}, \eqref{D511}, \eqref{D512} and \eqref{D52}, we have
\begin{equation}\label{D5-jia}
\begin{split}
D_{5}\leq C(1+\delta^{-1})\Big(1+\frac{1}{M_{k}^{l}}\Big)(1+r_{k}^{l})t+C\|\theta^{k+1}\|_{L_{t}^{\infty}H_{\vec x}^{N+2}}^{2}t+\delta\|D_{x}\theta^{k+1}\|_{L_{t}^{2}H_{\vec x}^{N+2}}^{2}.
\end{split}
\end{equation}

Similarly, we can write $D_{6}$ in the following form
\begin{equation}\label{D6}
\begin{split}
D_{6}=&\int_{0}^{t}\int_{\mathbb{T}^{3}}D_{x}^{\gamma}\Big(\frac{4(\theta^{k})^{3}}{3(1+4(\theta^{k})^{3})}\nabla_{x}\theta^{k}\cdot\vec u\Big)D_{x}^{\gamma}\theta^{k+1}\mathrm{d}\vec{x}\mathrm{d}s
\\=&\int_{0}^{t}\int_{\mathbb{T}^{3}}\Big(\frac{4(\theta^{k})^{3}}{3(1+4(\theta^{k})^{3})}\vec u\cdot D_{x}^{\gamma}\nabla_{x}\theta^{k}\Big)D_{x}^{\gamma}\theta^{k+1}\mathrm{d}\vec{x}\mathrm{d}s
\\&+\int_{0}^{t}\int_{\mathbb{T}^{3}}
\Big(\Big[D_{x}^{\gamma}, \frac{4(\theta^{k})^{3}}{3(1+4(\theta^{k})^{3})}\vec u\cdot\Big] \nabla_{x}\theta^{k}\Big)D_{x}^{\gamma}\theta^{k+1}\mathrm{d}\vec{x}\mathrm{d}s:=D_{61}+D_{62}.
\end{split}
\end{equation}
By integrating by parts, we write $D_{61}$ into the following two parts
\begin{equation}\label{D61}
\begin{split}
D_{61}=&\int_{0}^{t}\int_{\mathbb{T}^{3}}\Big(\frac{4(\theta^{k})^{3}}{3(1+4(\theta^{k})^{3})}\vec u\cdot D_{x}^{\gamma}\nabla_{x}\theta^{k}\Big)D_{x}^{\gamma}\theta^{k+1}\mathrm{d}\vec{x}\mathrm{d}s
\\=&\int_{0}^{t}\int_{\mathbb{T}^{3}}D_{x}\Big(\frac{4(\theta^{k})^{3}}{3(1+4(\theta^{k})^{3})}\vec u\cdot\Big) (D_{x}^{\gamma-1}\nabla_{x}\theta^{k})D_{x}^{\gamma}\theta^{k+1}\mathrm{d}\vec{x}\mathrm{d}s\\&+\int_{0}^{t}\int_{\mathbb{T}^{3}}\Big(\frac{4(\theta^{k})^{3}}{3(1+4(\theta^{k})^{3})}\vec u\cdot D_{x}^{\gamma-1}\nabla_{x}\theta^{k}\Big)D_{x}^{\gamma+1}\theta^{k+1}\mathrm{d}\vec{x}\mathrm{d}s:=D_{611}+D_{612},
\end{split}
\end{equation}
where
\begin{equation}\label{D611}
\begin{split}
D_{611}=&\int_{0}^{t}\int_{\mathbb{T}^{3}}D_{x}\Big(\frac{4(\theta^{k})^{3}}{3(1+4(\theta^{k})^{3})}\vec u\cdot\Big) (D_{x}^{\gamma-1}\nabla_{x}\theta^{k})D_{x}^{\gamma}\theta^{k+1}\mathrm{d}\vec{x}\mathrm{d}s
\\\leq& \Big\|D_{x}\Big(\frac{4(\theta^{k})^{3}}{3(1+4(\theta^{k})^{3})}\vec u\Big)\Big\|_{L^{\infty}_{t}L^{\infty}_{\vec x}}\|\theta^{k}\|_{L_{t}^{\infty}H_{\vec x}^{N+2}}\|\theta^{k+1}\|_{L_{t}^{\infty}H_{\vec x}^{N+2}}t
\\\leq& C\Big\|D_{x}\Big(\frac{4(\theta^{k})^{3}}{3(1+4(\theta^{k})^{3})}\vec u\Big)\Big\|_{L^{\infty}_{t}L^{\infty}_{\vec x}}^{2}\|\theta^{k}\|_{L_{t}^{\infty}H_{\vec x}^{N+2}}^{2}t+\|\theta^{k+1}\|_{L_{t}^{\infty}H_{\vec x}^{N+2}}^{2}t
\\\leq& C\Big(1+\frac{1}{M_{k}^{l}}\Big)(1+r_{k}^{l})t+\|\theta^{k+1}\|_{L_{t}^{\infty}H_{\vec x}^{N+2}}^{2}t
\end{split}
\end{equation}
and
\begin{equation}\label{D611}
\begin{split}
D_{612}=&\int_{0}^{t}\int_{\mathbb{T}^{3}}\Big(\frac{4(\theta^{k})^{3}}{3(1+4(\theta^{k})^{3})}\vec u\cdot D_{x}^{\gamma-1}\nabla_{x}\theta^{k}\Big)D_{x}^{\gamma+1}\theta^{k+1}\mathrm{d}\vec{x}\mathrm{d}s
\\\leq& \Big\|\frac{4(\theta^{k})^{3}}{3(1+4(\theta^{k})^{3})}\vec u\Big\|_{L^{\infty}_{t}L^{\infty}_{\vec x}}\|\theta^{k}\|_{L_{t}^{\infty}H_{\vec x}^{N+2}}\sqrt{t}\|D_{x}\theta^{k+1}\|_{L_{t}^{2}H_{\vec x}^{N+2}}
\\\leq& C\delta^{-1}\Big\|\frac{4(\theta^{k})^{3}}{3(1+4(\theta^{k})^{3})}\vec u\Big\|_{L^{\infty}_{t}L^{\infty}_{\vec x}}^{2}\|\theta^{k}\|_{L_{t}^{\infty}H_{\vec x}^{N+2}}^{2}t+\delta\|D_{x}\theta^{k+1}\|_{L_{t}^{2}H_{\vec x}^{N+2}}^{2}
\\\leq& C\delta^{-1}\Big(1+\frac{1}{M_{k}^{l}}\Big)(1+r_{k}^{l})t+\delta\|D_{x}\theta^{k+1}\|_{L_{t}^{2}H_{\vec x}^{N+2}}^{2}.
\end{split}
\end{equation}
Finally, 
\begin{equation}\label{D62}
\begin{split}
D_{62}=&\int_{0}^{t}\int_{\mathbb{T}^{3}}
\Big(\Big[D_{x}^{\gamma}, \frac{4(\theta^{k})^{3}}{3(1+4(\theta^{k})^{3})}\vec u\cdot\Big] \nabla_{x}\theta^{k}\Big)D_{x}^{\gamma}\theta^{k+1}\mathrm{d}\vec{x}\mathrm{d}s
\\\leq& \Big\|\frac{4(\theta^{k})^{3}}{3(1+4(\theta^{k})^{3})}\vec u\Big\|_{L_{t}^{\infty}H_{\vec x}^{N+2}}\|\theta^{k}\|_{L_{t}^{\infty}H_{\vec x}^{N+2}}\|\theta^{k+1}\|_{L_{t}^{\infty}H_{\vec x}^{N+2}}t
\\\leq& C\Big\|\frac{4(\theta^{k})^{3}}{3(1+4(\theta^{k})^{3})}\vec u\Big\|_{L_{t}^{\infty}H_{\vec x}^{N+2}}^{2}\|\theta^{k}\|_{L_{t}^{\infty}H_{\vec x}^{N+2}}^{2}t+\|\theta^{k+1}\|_{L_{t}^{\infty}H_{\vec x}^{N+2}}^{2}t
\\\leq& C\Big(1+\frac{1}{M_{k}^{l}}\Big)(1+r_{k}^{l})t+\|\theta^{k+1}\|_{L_{t}^{\infty}H_{\vec x}^{N+2}}^{2}t.
\end{split}
\end{equation}
Then, we can derive the following estimate
\begin{equation}\label{D6-jia}
\begin{split}
D_{6}\leq C(1+\delta^{-1})\Big(1+\frac{1}{M_{k}^{l}}\Big)(1+r_{k}^{l})t+C\|\theta^{k+1}\|_{L_{t}^{\infty}H_{\vec x}^{N+2}}^{2}t+\delta\|D_{x}\theta^{k+1}\|_{L_{t}^{2}H_{\vec x}^{N+2}}^{2}.
\end{split}
\end{equation}

Collecting \eqref{Eq 1}, \eqref{D1}, \eqref{D2}, \eqref{D3-jia}, \eqref{D4-jia}, \eqref{D5-jia} and \eqref{D6-jia} and summing up $\gamma$ from 0 to 6, we have
\begin{equation}
\begin{aligned}
&\|\theta^{k+1}(t)\|_{H_{\vec x}^{N+2}}^{2}+\Big(\frac{M_{k}^{3}}{C(r_{k}+r_{k}^{3}+1)}-4\delta\Big)\int_{0}^{t}\int_{\mathbb{T}^{3}}
\|D_{x}\theta^{k+1}\|_{H_{\vec x}^{N+2}}^{2}
\\\leq& C(1+r_{k}^{2l+2})\Big(1+\frac{1}{M_{k}^{2l}}\Big)(1+\delta^{-1})\|\theta^{k+1}(s)\|_{L_{t}^{2}H_{\vec x}^{N+2}}^{2}t+C(1+r_{k}^{2l+2})\Big(1+\frac{1}{M_{k}^{2l}}\Big)(1+\delta^{-1})t\\&+\|\theta_{0}^{0}\|_{H_{\vec x}^{N+2}}^{2}.
\end{aligned}
\end{equation}
Choosing $\delta=\frac{M_{k}^{3}}{8C(r_{k}+r_{k}^{3}+1)}$, then we obtain $\delta^{-1}\leq C\Big(1+\frac{1}{M_{k}^{3}}\Big)(1+r_{k}^{3})$. 

So,
\begin{equation}\label{theta-k+1}
\begin{aligned}
&\|\theta^{k+1}(t)\|_{H_{\vec x}^{13}}^{2}+\frac{M_{k}^{3}}{2C(r_{k}+r_{k}^{3}+1)}\int_{0}^{t}\int_{\mathbb{T}^{3}}
\|D_{x}\theta^{k+1}\|_{H_{\vec x}^{N+2}}^{2}
\\\leq& C(1+r_{k}^{2l+2})\Big(1+\frac{1}{M_{k}^{2l}}\Big)\|\theta^{k+1}(s)\|_{L_{t}^{\infty}H_{\vec x}^{N+2}}^{2}t+C(1+r_{k}^{2l+5})\Big(1+\frac{1}{M_{k}^{2l+3}}\Big)t\\&+\|\theta_{0}^{0}\|_{H_{\vec x}^{N+2}}^{2}.
\end{aligned}
\end{equation}
Choosing small $T=T(k)$ such that $0<T\leq CT\leq C(1+r_{k}^{2l+2})\Big(1+\frac{1}{M_{k}^{2l}}\Big)T\leq\frac{1}{2}$ and taking the sup over $t\in[0, T]$, then, we have
\begin{equation}\label{theta-k+1-jia}
\begin{aligned}
\|\theta^{k+1}\|_{L_{T}^{\infty}H_{\vec x}^{N+2}}^{2}
\leq C\Big(1+\frac{1}{M_{k}^{2l+3}}\Big)(1+r_{k}^{2l+5})T+2\|\theta_{0}^{0}\|_{H_{\vec x}^{N+2}}^{2}.
\end{aligned}
\end{equation}
Taking $t=T$ in \eqref{theta-k+1} and then substituting \eqref{theta-k+1-jia} into \eqref{theta-k+1}, we can obtain that
\begin{equation}\label{theta-k+1-jiajia}
\begin{aligned}
&\frac{M_{k}^{3}}{2C(r_{k}+r_{k}^{3}+1)}\int_{0}^{T}\int_{\mathbb{T}^{3}}
\|D_{x}\theta^{k+1}\|_{H_{\vec x}^{N+2}}^{2}
\leq C(1+r_{k}^{5l})\Big(1+\frac{1}{M_{k}^{5l}}\Big)T+C\|\theta_{0}^{0}\|_{H_{\vec x}^{N+2}}^{2}.
\end{aligned}
\end{equation}
Furthermore, we have
\begin{equation}\label{theta-k+1-jiajia}
\begin{aligned}
&\int_{0}^{T}\int_{\mathbb{T}^{3}}
\|D_{x}\theta^{k+1}\|_{H_{\vec x}^{N+2}}^{2}
\leq C(\|\theta_{0}^{0}\|_{H_{\vec x}^{N+2}})\Big(1+\frac{1}{M_{k}^{6l}}\Big)(1+r_{k}^{6l}).
\end{aligned}
\end{equation}

Integrating \eqref{research equations-4} from 0 to $t$, we have
\begin{equation}
\begin{aligned}
\theta^{k+1}(t)\geq& \theta^{0}(\vec x)-C(1+r_{k}^{l})\Big(1+\frac{1}{M_{k}^{l}}\Big)(1+r_{k+1})t
\\\geq& b-C(\|\theta_{0}^{0}\|_{H_{\vec x}^{N+2}}^{2})(1+r_{k}^{3l+5})\Big(1+\frac{1}{M_{k}^{3l+3}}\Big)T.
\end{aligned}
\end{equation}

For $k=0$, we have
\begin{equation}
\|\theta^{1}\|_{L_{T}^{\infty}H_{\vec x}^{N+2}}\leq C\Big(1+\frac{1}{b^{2l+3}}\Big)(1+\|\theta_{0}^{0}\|_{H_{\vec x}^{N+2}}^{2l+5})+2\|\theta_{0}^{0}\|_{H_{\vec x}^{N+2}}^{2}:=(r^{0})^{2}
\end{equation}
and 
\begin{equation}\nonumber
\inf_{t, \vec x}\theta^{1}(t, \vec x)\geq b-C(\|\theta_{0}^{0}\|_{H_{\vec x}^{N+2}}^{2})(1+\|\theta_{0}^{0}\|_{H_{\vec x}^{N+2}}^{3l+5})\Big(1+\frac{1}{b^{3l+3}}\Big)T\geq b/2
\end{equation}
for sufficiently small $T_{1}$ and $0<T\leq T_{1}\leq 1$.

We assume $\|\theta^{k}\|_{L_{T}^{\infty}H_{\vec x}^{N+2}}\leq r^{0}$ and $M_{k}\geq \frac{b}{2}$, then, for sufficiently small $T_{2}$ and $0<T\leq T_{2}\leq 1$, we have
\begin{equation}\label{theta-k+1-jiajia}
\begin{aligned}
\|\theta^{k+1}\|_{L_{T}^{\infty}H_{\vec x}^{N+2}}^{2}
\leq C\sup_{\begin{matrix}
  0\leq r_{k}\leq r^{0}, \\
  M_{k}\geq b/2
\end{matrix}}\Big(1+\frac{1}{M_{k}^{2l+3}}\Big)(1+r_{k}^{2l+5})T+2\|\theta_{0}^{0}\|_{H_{\vec x}^{N+2}}^{2}\leq r^{0}.
\end{aligned}
\end{equation}
and 
\begin{equation}
\begin{aligned}
\theta^{k+1}(t)\geq& b-C(\|\theta_{0}^{0}\|_{H_{\vec x}^{N+2}}^{2})\sup_{\begin{matrix}
  0\leq r_{k}\leq r^{0}, \\
  M_{k}\geq b/2
\end{matrix}}(1+r_{k}^{3l+5})\Big(1+\frac{1}{M_{k}^{3l+3}}\Big)T\geq \frac{b}{2}.
\end{aligned}
\end{equation}
Furthermore, we have 
\begin{equation}\label{theta-k+1-jiajiajia}
\begin{aligned}
&\int_{0}^{T}\int_{\mathbb{T}^{3}}
\|D_{x}\theta^{k+1}\|_{H_{\vec x}^{N+2}}^{2}
\leq C(\|\theta_{0}^{0}\|_{H_{\vec x}^{N+2}})\Big(1+\frac{1}{M_{k}^{6l}}\Big)(1+r_{k}^{6l})\leq C(b, r^{0}).
\end{aligned}
\end{equation}

Therefore, we have proved the uniform energy estimate and lower bound of $\theta^{k}$ for any $k\geq 0$. 

Define $P_{k+1}=\theta^{k+1}-\theta^{k}$. Then we get the system about $(P_{k+1}, Q_{k+1})$ as follows
\begin{equation}\label{equations about Pk and Qk}\left\{
\begin{split}
&\partial_{t}P_{k+1}-\frac{4(\theta^{k})^{3}}{(1+4(\theta^{k})^{3})}\Delta_{x}P_{k+1}=J(\theta^{k})-J( \theta^{k-1})+A
,\ \ \mathrm{in} \ \ (0, T]\times\mathbb{T}^{3},\\
&P_{k+1}(0, \vec{x})=0\ \ \mathrm{in} \ \ \mathbb{T}^{3},\\
\end{split}\right.
\end{equation}
where
\begin{equation}\nonumber
\begin{aligned}
J(\theta^{k})=&-\frac{1}{(1
+4(\theta^{k})^{3})}\vec{u}\cdot\nabla_{x}\theta^{k}-\frac{\theta^{k}}{(1+4(\theta^{k})^{3})}\mathrm{div}_{x}\vec{u}
\\&+\frac{4(\theta^{k})^{2}}{3(1+4(\theta^{k})^{3})}|\nabla_{x}\theta^{k}|^{2}+\frac{4(\theta^{k})^{3}}{3(1+4(\theta^{k})^{3})}\nabla_{x}\theta^{k}\cdot\vec u
\end{aligned}
\end{equation}
and
\begin{equation}\nonumber
\begin{aligned}
A=\Big(\frac{4(\theta^{k})^{3}}{(1+4(\theta^{k})^{3})}-\frac{4(\theta^{k-1})^{3}}{(1+4(\theta^{k-1})^{3})}
\Big)\Delta_{x}\theta^{k}.
\end{aligned}
\end{equation}

For $t\in [0,T]$, multiplying $2P_{k+1}$ on both sides of $(\ref{equations about Pk and Qk})_{2}$, we have
\begin{equation}\label{Eq 2}
\begin{split}
&\int_{\mathbb{T}^{3}}|P_{k+1}|^{2}(t)\mathrm{d}\vec{x}
+\int_{0}^{t}\int_{\mathbb{T}^{3}}\frac{8(\theta^{k})^{3}}{(1+4(\theta^{k})^{3})}|\nabla_{x}P_{k+1}|^{2}\mathrm{d}\vec{x}\mathrm{d}s
\\=&-\int_{0}^{t}\int_{\mathbb{T}^{3}}\nabla_{x}\Big(\frac{8(\theta^{k})^{3}}{(1+4(\theta^{k})^{3})}\Big)\nabla_{x}P_{k+1}P_{k+1}
\mathrm{d}\vec{x}\mathrm{d}s
+2\int_{0}^{t}\int_{\mathbb{T}^{3}}(J(\theta^{k})-J( \theta^{k-1}))P_{k+1}\mathrm{d}\vec{x}\mathrm{d}s\\&+2\int_{0}^{t}\int_{\mathbb{T}^{3}}AP_{k+1}\mathrm{d}\vec{x}\mathrm{d}s:=G_{1}+G_{2}+G_{3},
\end{split}
\end{equation}
where 
\begin{equation}\label{G1}
\begin{split}
G_{1}=&-\int_{0}^{t}\int_{\mathbb{T}^{3}}\nabla_{x}\Big(\frac{8(\theta^{k})^{3}}{(1+4(\theta^{k})^{3})}\Big)\nabla_{x}P_{k+1}P_{k+1}
\mathrm{d}\vec{x}\mathrm{d}s\\=&-\frac{1}{2}\int_{0}^{t}\int_{\mathbb{T}^{3}}\Delta_{x}\Big(\frac{4(\theta^{k})^{3}}{(1+4(\theta^{k})^{3})}\Big)
P_{k+1}^{2}\mathrm{d}\vec{x}\mathrm{d}s\\\leq& C(b, r^{0})\|P_{k+1}\|^{2}_{L^{\infty}_{T}L^{2}_{\vec x}}T,
\end{split}
\end{equation}

\begin{equation}\label{G2}
\begin{split}
G_{2}=&2\int_{0}^{t}\int_{\mathbb{T}^{3}}J(\theta^{k})-J( \theta^{k-1}))P_{k+1}\mathrm{d}\vec{x}\mathrm{d}s\\\leq&C(b, r^{0})\int_{0}^{t}\int_{\mathbb{T}^{3}}|P_{k}|
|P_{k+1}|\mathrm{d}\vec{x}\mathrm{d}s\\\leq& C(b, r^{0})(\|P_{k}\|^{2}_{L^{\infty}_{T}L^{2}_{\vec x}}+\|P_{k+1}\|^{2}_{L^{\infty}_{T}L^{2}_{\vec x}})T
\end{split}
\end{equation}
and
\begin{equation}\label{G3}
\begin{split}
G_{3}=&2\int_{0}^{t}\int_{\mathbb{T}^{3}}AP_{k+1}\mathrm{d}\vec{x}\mathrm{d}s\\\leq&C(b, r^{0})\int_{0}^{t}\int_{\mathbb{T}^{3}}|P_{k}|
|P_{k+1}|\mathrm{d}\vec{x}\mathrm{d}s\\\leq& C(b, r^{0})(\|P_{k}\|^{2}_{L^{\infty}_{T}L^{2}_{\vec x}}+\|P_{k+1}\|^{2}_{L^{\infty}_{T}L^{2}_{\vec x}})T.
\end{split}
\end{equation}
So,
\begin{equation}\label{G3-jia}
\begin{split}
&\int_{\mathbb{T}^{3}}|P_{k+1}|^{2}(t)\mathrm{d}\vec{x}\leq C(b, r^{0})(\|P_{k}\|^{2}_{L^{\infty}_{T}L^{2}_{\vec x}}+\|P_{k+1}\|^{2}_{L^{\infty}_{T}L^{2}_{\vec x}})T.
\end{split}
\end{equation}
Taking the sup over $t\in(0, T)$ on both sides of \eqref{G3-jia}, we have
\begin{equation}\label{G3-jia-jia}
\begin{split}
&\|P_{k+1}\|^{2}_{L^{\infty}_{T}L^{2}_{\vec x}}\leq C(b, r^{0})(\|P_{k}\|^{2}_{L^{\infty}_{T}L^{2}_{\vec x}}+\|P_{k+1}\|^{2}_{L^{\infty}_{T}L^{2}_{\vec x}})T.
\end{split}
\end{equation}

Now, we choose small $T=T_{3}$ such that $\sqrt{C(b, r^{0})T}<1$. Then we obtain that
\begin{equation}\nonumber
\begin{split}
&\|P_{k+1}\|_{L^{\infty}_{T}L^{2}_{\vec x}}
\leq C_{1}(b, r^{0})\sqrt{T}\|P_{k}\|_{L^{\infty}_{T}L^{2}_{\vec x}}.
\end{split}
\end{equation}
Furthermore, we choose small $T=T_{4}$ such that $d=C_{1}(b, r^{0})\sqrt{T}<1$.
So, it holds that
\begin{equation}\nonumber
\begin{split}
&\sum_{k=0}^{\infty}\|P_{k}\|_{L^{\infty}_{T}L^{2}_{\vec x}}
<\infty.
\end{split}
\end{equation}
Therefore, we have $\theta^{k}$ converges strongly to $\theta_{0}$ in $L^{\infty}(0, T; L^{2}(\mathbb{T}^{3}))$. So, we can choose
\begin{equation}\nonumber
T=\min\{1, T_{1}, T_{2}, T_{3}, T_{4}\}
\end{equation}
which is independent of $k$ to satisfy all the above uniform estimates and contraction principle. Therefore, we can infer that $\theta_{0}$ is the unique solution of \eqref{equations about f0 epsilon and theta} where $\theta_{0}\in G$. We can also verify that $T^{k+1}\geq T>0$ for $k\geq 0$. Since $f_{0}=\theta_{0}^{4}$, the existence, uniqueness of and uniform estimates of $f_{0}$ can be directly verified.

\noindent\textbf{Step 2.} Construction of first-order terms.

The first-order initial layer $(f_{I, 1}, \theta_{I, 1})$ is defined as
\begin{equation}\label{equations about f I1,theta I1}\left\{
\begin{split}
&\vec w\cdot\nabla_{x}f_{I, 0}+\frac{\partial f_{I, 1}}{\partial\tau}+f_{I, 1}=B_{1}(\theta_{0}^{0}+\theta_{I, 0}; \theta_{1}^{0}+\theta_{I, 1})-B_{1}(\theta_{0}^{0}; \theta_{1}^{0}),\\
&\frac{\partial \theta_{I, 1}}{\partial\tau}=(\overline{f_{I, 1}}-B_{1}(\theta_{0}^{0}+\theta_{I, 0}; \theta_{1}^{0}+\theta_{I, 1})+B_{1}(\theta_{0}^{0}; \theta_{1}^{0})),\\
&f_{I, 1}(0, \vec{x}, \vec w)=-f_{1}^{0}(\vec x, \vec w),\\
&\theta_{I, 1}(0, \vec{x})=-\theta_{1}^{0}(\vec x),\\
&\lim_{\tau\rightarrow\infty}f_{I, 1}(\tau, \vec{x}, \vec{w})=0,\ \ \lim_{\tau\rightarrow\infty}\theta_{I, 1}(\tau, \vec{x})=0,
\end{split}\right.
\end{equation}
where we have used the notation $\theta_{0}^{0}=\theta_{0}|_{t=0}$ and $\theta_{1}^{0}=\theta_{1}|_{t=0}$.

\begin{thm}\label{3.3}
The problem \eqref{equations about f I1,theta I1} has a unique solution $(f_{I, 1}, \theta_{I, 1})\in (C^{1}([0, \infty); L^{\infty}_{\vec w}H_{\vec x}^{N+1})\cap L^{1}([0, \infty); L^{\infty}_{\vec w}H_{\vec x}^{N+1})\times (C^{1}([0, \infty); H_{\vec x}^{N+1})\cap L^{1}([0, \infty); H_{\vec x}^{N+1})$. Furthermore, we have
\begin{equation}\label{ineq-theta00-jia01}
\|\theta_{1}^{0}(\vec x)\|_{H_{\vec x}^{N+1}}\leq C,
\end{equation}
\begin{equation}\label{ineq-theta00-jia11}
\|e^{\sigma\tau}f_{I, 1}\|_{L^{\infty}([0, \infty); L^{\infty}_{\vec w}H_{\vec x}^{N+1})}+\|e^{\sigma\tau}\theta_{I, 1}\|_{L^{\infty}([0, \infty); H_{\vec x}^{N+1})}\leq C
\end{equation}
and
\begin{equation}\label{ineq-theta00-jia21}
\|e^{\sigma\tau}f_{I, 1}\|_{L^{1}([0, \infty); L^{\infty}_{\vec w}H_{\vec x}^{N+1})}+\|e^{\sigma\tau}\theta_{I, 1}\|_{L^{1}([0, \infty); H_{\vec x}^{N+1})}\leq C,
\end{equation}
where $\sigma>0$, suitably small.
\end{thm}
\noindent\textbf{Proof.} From the equations $\eqref{equations about f I1,theta I1}_{1}$ and $\eqref{equations about f I1,theta I1}_{2}$, we can derive 
\begin{equation}
\partial_{\tau}(\overline {f_{I, 1}}+\theta_{I, 1})=-\langle\vec w\cdot\nabla_{x}f_{I, 0}\rangle.
\end{equation}
Then, 
\begin{equation}
\overline {f_{I, 1}}(\tau, \vec x)+\theta_{I, 1}(\tau, \vec x)=-\overline{f_{1}^{0}}(\vec x)-\theta_{1}^{0}(\vec x)-\int_{0}^{\tau}\langle\vec w\cdot\nabla_{x}f_{I, 0}\rangle\mathrm{d}s
\end{equation}
By the fact that $f_{I, 1}, \theta_{I, 1}\rightarrow 0$ as $\tau\rightarrow\infty$, we can obtain
\begin{equation}\label{daoshuweiling-1}
\overline {f_{I, 1}}+\theta_{I, 1}\equiv 0,
\end{equation}
which further implies
\begin{equation}
\overline{f_{1}^{0}}(\vec x)+\theta_{1}^{0}(\vec x)=-\int_{0}^{\infty}\langle\vec w\cdot\nabla_{x}f_{I, 0}\rangle\mathrm{d}s:=l_{1}(\vec x).
\end{equation}
According to \eqref{fk} and the fact that $\theta_{1}(0, \vec x)=4(\theta^{0}_{0})^{3}\theta_{1}^{0}(\vec x)$, we can derive the following equality
\begin{equation}
4(\theta_{0}^{0})^{3}\theta_{1}^{0}+\theta_{1}^{0}=l_{1}(\vec x).
\end{equation}
Then, we can directly get 
\begin{equation}
\theta_{1}^{0}(\vec x)=\frac{1}{4(\theta_{0}^{0})^{3}+1}l_{1}(\vec x)\in H_{\vec x}^{N+1}.
\end{equation}
Simultaneously, we can get $f_{1}^{0}\in L^{\infty}_{\vec w}H_{\vec x}^{N+1}$ by \eqref{fk}.

We can write $\eqref{equations about f I1,theta I1}_{2}$ in the following form
\begin{equation}
\frac{\partial \theta_{I, 1}}{\partial\tau}+(1+4(\theta_{0}^{0}+\theta_{I, 0})^{3})\theta_{I, 1}=-4(\theta_{0}^{0}+\theta_{I, 0})^{3}\theta_{1}^{0}+4(\theta_{0}^{0})^{3}\theta_{1}^{0}.
\end{equation}
Then, we have
\begin{equation}
\theta_{I, 1}=e^{-(1+4(\theta_{0}^{0}+\theta_{I, 0})^{3})\tau}(-\theta_{1}^{0})+\int_{0}^{\tau}e^{-(1+4(\theta_{0}^{0}+\theta_{I, 0})^{3})(\tau-s)}(-4(\theta_{0}^{0}+\theta_{I, 0})^{3}\theta_{1}^{0}+4(\theta_{0}^{0})^{3}\theta_{1}^{0})\mathrm{d}s,
\end{equation}
which implies the global existence of $\theta_{I, 1}$. 

According to $\eqref{equations about f I0,theta I0}_{1}$ and $\eqref{equations about f I0,theta I0}_{3}$, we can write $f_{I, 1}$ in the following form
\begin{equation}
f_{I, 1}=-e^{-\tau}f_{1}^{0}+\int_{0}^{\tau}e^{s-\tau}(B_{1}(\theta_{0}^{0}+\theta_{I, 0}; \theta_{1}^{0}+\theta_{I, 1})-B_{1}(\theta_{0}^{0}; \theta_{1}^{0})-\vec w\cdot\nabla_{x}f_{I, 0})\mathrm{d}s,
\end{equation}
which further implies the existence solution $f_{I, 1}\in C^{1}([0, \infty); L^{\infty}_{\vec w}H_{\vec x}^{N+1}))$. \hfill$\Box$

Using the similar method to Theorem \ref{3.1}, we can get the estimates \eqref{ineq-theta00-jia01}, \eqref{ineq-theta00-jia11} and \eqref{ineq-theta00-jia21}.

We write the equations about $\theta_{1}$ and $f_{1}$ as follows
\begin{equation}\label{equations about theta 1 overline f epsilon 1-jia}\left\{
\begin{split}
&f_{1}(t, \vec{x}, \vec w)=B_{1}(\theta_{0}; \theta_{1})-\vec w\cdot\nabla_{x}B(\theta_{0})\\
&\partial_{t}(\theta_{1}+4\theta_{0}^{3}\theta_{1})+\mathrm{div}(\vec u\theta_{1})-\Delta_{x}\Big(\theta_{1}+\frac{4}{3}\theta_{0}^{3}\theta_{1}\Big)=0,\ \ \mathrm{in}\ \ (0, T]\times\mathbb{T}^{3},\\
&\theta_{1}(0, \vec{x})=\theta_{1}^{0}(\vec x),\;\;f_{1}(0, \vec{x}, \vec w)=B_{1}(\theta_{0}^{0}; \theta_{1}^{0})-\vec w\cdot\nabla_{x}B(\theta_{0}^{0})~~\mathrm{in}  \ \ \mathbb{T}^{3}.
\end{split}\right.
\end{equation}
\begin{thm}\label{3.4}
The problem \eqref{equations about theta 1 overline f epsilon 1-jia} has a unique solution on $[0, T]$, and $(\theta_{1}, f_{1})\in C^{0}([0, T]; H_{\vec x}^{4N+1})\cap L^{2}(0, T; H_{\vec x}^{N+2})\times C^{0}([0, T]; L_{\vec w}^{\infty}H_{\vec x}^{N+1})\cap L^{2}(0, T; L_{\vec w}^{\infty}H_{\vec x}^{N+2})$. Furthermore, we have
\begin{equation}\label{ineq-thetak0}
\|\theta_{1}^{0}(\vec x)\|_{H_{\vec x}^{N+1}}+\|f_{1}^{0}(\vec x, \vec w)\|_{L^{\infty}_{\vec w}H_{\vec x}^{N+1}}\leq C,
\end{equation}
\begin{equation}
\|\theta_{1}\|_{L_{T}^{\infty}H_{\vec x}^{N+1}}+\|\theta_{1}\|_{L_{T}^{2}H_{\vec x}^{N+2}}\leq C
\end{equation}
and
\begin{equation}
\|f_{1}\|_{L_{T}^{\infty}L_{\vec w}^{\infty}H_{\vec x}^{N+1}}+\|f_{1}\|_{L_{T}^{2}L_{\vec w}^{\infty}H_{\vec x}^{N+2}}\leq C.
\end{equation}
\end{thm}
Note that $B_{1}(\theta_{0}; \theta_{1})=4\theta_{0}^{3}\theta_{1}$. We can write system \eqref{equations about theta 1 overline f epsilon 1-jia} as follows
\begin{equation}\label{equations about theta 1 overline f epsilon 1 jia1-jia}\left\{
\begin{split}
&f_{1}(t, x, w)=B_{1}(\theta_{0}; \theta_{1})-\vec w\cdot\nabla_{x}B(\theta_{0}),\\
&\frac{\partial \theta_{1}}{\partial t}-\frac{4\theta_{0}^{3}}{3(1+4\theta_{0}^{3})}\Delta_{x}\theta_{1}=F^{\star}\ \ \mathrm{in}\ \ (0, T]\times\mathbb{T}^{3},\\
&\theta_{1}(0, \vec x)=\theta_{1}^{0}(\vec x),~~\mathrm{in}  \ \ \mathbb{T}^{3},\\
\end{split}\right.
\end{equation}
where 
\begin{equation}
\begin{aligned}
F^{\star}=-\frac{\partial_{t}(1+4\theta_{0}^{3})}{1+4\theta_{0}^{3}}\theta_{1}-\frac{1}{1+4\theta_{0}^{3}}(\mathrm{div}\vec u)\theta_{1}-\frac{1}{1+4\theta_{0}^{3}}\vec u\cdot\nabla_{x}\theta_{1}+\frac{\nabla_{x}(1+\frac{4}{3}\theta_{0}^{3})}{1+4\theta_{0}^{3}}\cdot\nabla_{x}\theta_{1}
-\frac{1+\frac{4}{3}\theta_{0}^{3}}{1+4\theta_{0}^{3}}\theta_{1}.
\end{aligned}
\end{equation}

The local existence of smooth solutions can be established by an iterative method as in Theorem \ref{3.2}. To close Theorem \ref{3.4}, one only needs to establish the a priori estimates.

Applying the operator $D_{x}^{\gamma}$ to the equation $(\ref{equations about theta 1 overline f epsilon 1 jia1-jia})_{2}$, multiplying the resulting equation with $2(D_{x}^{\gamma}\theta_{1})$ in $L^{2}((0, t)\times\mathbb{T}^{3})$, we have
\begin{equation}\label{Eq 1}
\begin{split}
&\int_{\mathbb{T}^{3}}|D_{x}^{\gamma}\theta_{1}|^{2}(t)\mathrm{d}\vec{x}-\int_{\mathbb{T}^{3}}|D_{x}^{\gamma}\theta^{0}_{1}|^{2}(t)\mathrm{d}\vec{x}
\\=&-\int_{0}^{t}\int_{\mathbb{T}^{3}}\frac{3+4\theta_{0}^{3}}{3(1+4\theta_{0}^{3})}|D_{x}^{\gamma}\nabla_{x}\theta_{1}|^{2}\mathrm{d}\vec{x}\mathrm{d}s
+\int_{0}^{t}\int_{\mathbb{T}^{3}}\Big(\Big[D_{x}^{\gamma}, \frac{3+4\theta_{0}^{3}}{3(1+4\theta_{0}^{3})}\Big]\Delta_{x}\theta_{1}\Big)D_{x}^{\gamma}\theta_{1}\mathrm{d}\vec{x}\mathrm{d}s
\\&+\int_{0}^{t}\int_{\mathbb{T}^{3}}D_{x}^{\gamma}F^{\star}D_{x}^{\gamma}\theta_{1}\mathrm{d}\vec{x}\mathrm{d}s
,
\end{split}
\end{equation}
where $0\leq\gamma\leq N+1$ and $t\in [0, \tilde T]$, $\tilde T$ is the local existence time of \eqref{equations about theta 1 overline f epsilon 1 jia1-jia}.

For some $c>0$, we have
\begin{equation}
\begin{aligned}
-\int_{0}^{t}\int_{\mathbb{T}^{3}}\frac{3+4\theta_{0}^{3}}{3(1+4\theta_{0}^{3})}|D_{x}^{\gamma}\nabla_{x}\theta_{1}|^{2}
\mathrm{d}\vec{x}\mathrm{d}s\leq & -\frac{1+\frac{1}{2}a^{3}}{C(1+(r^{0})^{3})}\int_{\mathbb{T}^{3}}|D_{x}^{\gamma+1}\theta_{1}|^{2}
\mathrm{d}\vec{x}\mathrm{d}s\\:=&-c\int_{\mathbb{T}^{3}}|D_{x}^{\gamma+1}\theta_{1}|^{2}
\mathrm{d}\vec{x}\mathrm{d}s.
\end{aligned}
\end{equation}
The remaining two terms can be estimated as follows
\begin{equation}
\begin{aligned}
&\int_{0}^{t}\int_{\mathbb{T}^{3}}\Big(\Big[D_{x}^{\gamma}, \frac{3+4\theta_{0}^{3}}{3(1+4\theta_{0}^{3})}\Big]\Delta_{x}\theta_{1}\Big)D_{x}^{\gamma}\theta_{1}
\mathrm{d}\vec{x}\mathrm{d}s
\\\leq&\delta\int_{0}^{t}\int_{\mathbb{T}^{3}}|D_{x}^{\gamma+1}\theta_{1}|^{2}\mathrm{d}\vec{x}\mathrm{d}s+C(a, r^{0}, \delta)\int_{0}^{t}\int_{\mathbb{T}^{3}}|D_{x}^{\gamma}\theta_{1}|^{2}\mathrm{d}\vec{x}\mathrm{d}s
\end{aligned}
\end{equation}
and
\begin{equation}
\begin{aligned}
&\int_{0}^{t}\int_{\mathbb{T}^{3}}D_{x}^{\gamma}F^{\star}D_{x}^{\gamma}\theta_{1}\mathrm{d}\vec{x}\mathrm{d}s
\\\leq&\delta\int_{0}^{t}\int_{\mathbb{T}^{3}}|D_{x}^{\gamma+1}\theta_{1}|^{2}\mathrm{d}\vec{x}\mathrm{d}s+C(a, r^{0}, \delta)\int_{0}^{t}\int_{\mathbb{T}^{3}}|D_{x}^{\gamma}\theta_{1}|^{2}\mathrm{d}\vec{x}\mathrm{d}s.
\end{aligned}
\end{equation}
Choosing $\delta=\frac{c}{4}$, we have
\begin{equation}\label{Eq 1}
\begin{split}
&\int_{\mathbb{T}^{3}}|D_{x}^{\gamma}\theta_{1}|^{2}(t)\mathrm{d}\vec{x}
+\frac{c}{2}\int_{0}^{t}\int_{\mathbb{T}^{3}}|D_{x}^{\gamma+1}\theta_{1}|^{2}\mathrm{d}\vec{x}\mathrm{d}s
\leq\int_{\mathbb{T}^{3}}|D_{x}^{\gamma}\theta^{0}_{1}|^{2}\mathrm{d}\vec{x}
+C(a, r^{0})\int_{0}^{t}\int_{\mathbb{T}^{3}}|D_{x}^{\gamma}\theta_{1}|^{2}\mathrm{d}\vec{x}\mathrm{d}s.
\end{split}
\end{equation}
By the Gronwall's inequality, we have the following uniform estimate
\begin{equation}
\|\theta_{1}\|_{L_{t}^{\infty}H_{\vec x}^{N+1}}+\|\theta_{1}\|_{L_{t}^{2}H_{\vec x}^{N+2}}\leq C(a, r^{0},t),
\end{equation}
where $C(a, r^{0},t)\leq C$, for $0\leq t\leq T\leq 1$. 

So, the local existence of \eqref{equations about theta 1 overline f epsilon 1 jia1-jia} can be extended to the interval $[0, T]$.

Then, 
\begin{equation}
\|\theta_{1}\|_{L_{T}^{\infty}H_{\vec x}^{N+1}}+\|\theta_{1}\|_{L_{T}^{2}H_{\vec x}^{N+2}}\leq C.
\end{equation}

Since $f_{1}(t, \vec{x}, \vec w)=B_{1}(\theta_{0}; \theta_{1})-\vec w\cdot\nabla_{x}B(\theta_{0})$, we also have
\begin{equation}
\|f_{1}\|_{L_{T}^{\infty}L_{\vec w}^{\infty}H_{\vec x}^{N+1}}+\|f_{1}\|_{L_{T}^{2}L_{\vec w}^{\infty}H_{\vec x}^{N+2}}\leq C.
\end{equation}

\noindent\textbf{Step 3.} Construction of $mth$ order terms.

We define the equation about $f_{I, m}, \theta_{I, m}$, $2\leq m\leq N$, as follows
\begin{equation}\label{equation about inik}\left\{
\begin{split}
&\partial_{\tau}f_{I, m}+\vec{w}\cdot\nabla_{x}f_{I, m-1}+f_{I, m}+f_{I, m-2}-\overline{f_{I, m-2}}-(B_{k}(A_{0}+\theta_{I, 0};A_{m}+\theta_{I, m})-B_{k}(A_{0};A_{m}))=0,\\
&\partial_{\tau}\theta_{I, m}+\mathrm{div}(\vec v_{m}f_{I, m-2})-\overline{f_{I, m}}-\Delta_{x}f_{I, m-2}+(B_{k}(A_{0}+\theta_{I, 0};A_{m}+\theta_{I, m})-B_{k}(A_{0};A_{m}))=0,\\
&(f_{I, m}(0, \vec x), \theta_{I, m}(0, \vec x))=(-f_{m}^{0}(\vec x, \vec w), -\theta_{m}^{0}(\vec x)),\\
&\lim_{\tau\rightarrow\infty}(f_{I, m}(\tau, \vec x), \theta_{I, m}(\tau, \vec x))=(0, 0).\\
\end{split}\right.
\end{equation}

\begin{thm}\label{3.5}
The problem \eqref{equation about inik} has a unique solution $(f_{I, m}, \theta_{I, m})\in (C^{1}([0, \infty); L^{\infty}_{\vec w}H_{\vec x}^{N+2-m})\cap L^{1}([0, \infty); L^{\infty}_{\vec w}H_{\vec x}^{N+2-m}))\times (C^{1}([0, \infty); H_{\vec x}^{N+2-m})\cap L^{1}([0, \infty); H_{\vec x}^{N+2-m}))$. Furthermore, we have
\begin{equation}\label{ineq-thetak0}
\|\theta_{m}^{0}(\vec x)\|_{H_{\vec x}^{N+2-m}}+\|f_{m}^{0}(\vec x, \vec w)\|_{L^{\infty}_{\vec w}H_{\vec x}^{N+2-m}}\leq C,
\end{equation}
\begin{equation}\label{ineq-thetak1}
\|e^{\sigma\tau}f_{I, m}\|_{L^{\infty}([0, \infty); L^{\infty}_{\vec w}H_{\vec x}^{N+2-m})}+\|e^{\sigma\tau}\theta_{I, m}\|_{L^{\infty}([0, \infty); H_{\vec x}^{N+2-m})}\leq C
\end{equation}
and
\begin{equation}\label{ineq-thetak2}
\|e^{\sigma\tau}f_{I, m}\|_{L^{1}([0, \infty); L^{\infty}_{\vec w}H_{\vec x}^{N+2-m})}+\|e^{\sigma\tau}\theta_{I, m}\|_{L^{1}([0, \infty); H_{\vec x}^{N+2-m})}\leq C,
\end{equation}
where $\sigma>0$, suitably small.
\end{thm}
The proof process is similar to Theorem \ref{3.1}, so, we omit the details here.

We write the equations about $\theta_{m}$ and $f_{m}$ as follows
\begin{equation}\label{equations about theta 1 overline f epsilon 1}\left\{
\begin{split}
&f_{m}=B_{m}(\theta_{0};\theta_{m})+\overline{f_{m-2}}-f_{m-2}-\partial_{t}f_{m-2}-\vec{w}\cdot\nabla_{x}f_{m-1},\\
&\partial_{t}(\theta_{m}+4\theta_{0}^{3}\theta_{m})+\mathrm{div}(\vec u\theta_{m})-\Delta_{x}\Big(\theta_{m}+\frac{4}{3}\theta_{0}^{3}\theta_{m}\Big)=\tilde{F}_{m-1},\ \ \mathrm{in}\ \ (0, T]\times\mathbb{T}^{3},\\
&\tilde{F}_{m-1}=-\partial_{t}\langle\overline{f_{m-2}}-f_{m-2}+\vec{w}\cdot\nabla_{x}f_{m-1}-\partial_{t}f_{m-2}\rangle
+\langle\vec{w}\cdot\nabla_{x}f_{m-1}+(\vec{w}\cdot\nabla_{x})^{2}(\overline{f_{m-2}}-f_{m-2})
\\&-(\vec{w}\cdot\nabla_{x})^{3}f_{m-1}
-(\vec{w}\cdot\nabla_{x})^{2}\partial_{t}f_{m-2}+\vec{w}\cdot\nabla_{x}\partial_{t}f_{m-1}\rangle\\&-\frac{4}{3}\Delta_{x}\Big(\sum_{\begin{matrix}
  i+j+l+k=m, \\
  i,j,l,k\geq 1
\end{matrix}}\theta_{i}\theta_{j}\theta_{l}\theta_{k}\Big),\ \ \mathrm{in}\ \ (0, T]\times\mathbb{T}^{3},\\
&\theta_{m}(0, \vec{x})=\theta_{m}^{0}(\vec x),\;\;f_{m}(0, \vec{x}, \vec w)=f_{m}^{0}(\vec x, \vec w)~~\mathrm{in}  \ \ \mathbb{T}^{3}\times\mathbb{S}^{2}.
\end{split}\right.
\end{equation}
\begin{thm}\label{3.6}
The problem \eqref{equations about theta 1 overline f epsilon 1} has a unique solution on $[0, T]$, and $(\theta_{m}, \overline{f_{m}})\in C^{0}([0, T]; H_{\vec x}^{N+2-m})\cap L^{2}(0, T; H_{\vec x}^{N+3-m})\times C^{0}([0, T]; L_{\vec w}^{\infty}H_{\vec x}^{N+2-m})\cap L^{2}(0, T; L_{\vec w}^{\infty}H_{\vec x}^{N+3-m})$. Furthermore, we have
\begin{equation}\label{ineq-thetam1}
\|\theta_{m}\|_{L_{T}^{\infty}H_{\vec x}^{N+2-m}}+\|\theta_{m}\|_{L^{2}_{T}H_{\vec x}^{N+3-m}}\leq C
\end{equation}
and
\begin{equation}\label{ineq-thetam2}
\|f_{m}\|_{L_{T}^{\infty}L_{\vec w}^{\infty}H_{\vec x}^{N+2-m}}+\|f_{m}\|_{L^{2}_{T}L_{\vec w}^{\infty}H_{\vec x}^{N+3-m}}\leq C.
\end{equation}
\end{thm}
The proof process is similar to Theorem \ref{3.2}, so, we omit the details here.

\section{Diffusive Limit}\label{si}

In this section, we prove Theorem \ref{result} by estimating the difference between the solution $(f^{\epsilon}, \theta^{\epsilon})$ to system \eqref{research equations} and the constructed approximate solution $(f^{N}, \theta^{N})$ where 
\begin{equation}\label{expansion of f epsilon}
f^{N}=\sum_{k=0}^{N}\epsilon^{k}f_{k}+\sum_{k=0}^{N}\epsilon^{k}f_{I, k}
\end{equation} 
and
\begin{equation}\label{expansion of theta epsilon-jia}
\theta^{N}=\sum_{k=0}^{N}\epsilon^{k}\theta_{k}+\sum_{k=0}^{N}\epsilon^{k}\theta_{I, k}.
\end{equation} 
The remainder can be defined as
\begin{equation}\nonumber
\begin{split}
f_{r}=f^{\epsilon}-f^{N},\ \ \theta_{r}=\theta^{\epsilon}-\theta^{N},
\end{split}
\end{equation}
functions $(f_{r}, \theta_{r})$ then satisfy
\begin{equation}\label{reeq1}
\epsilon^{2}\partial_{t}f_{r}+\epsilon\vec{w}\cdot\nabla_{x}f_{r}+\epsilon^{2}(f_{r}-\overline{f_{r}})+f_{r}
-(\theta^{N}+\theta_{r})^{4}+(\theta^{N})^{4}=-\mathcal{L}_{1}(f^{N}, \theta^{N}),
\end{equation}
\begin{equation}\label{reeq2}
\epsilon^{2}\partial_{t}\theta_{r}+\epsilon^{2}\mathrm{div}{(\vec u\theta_{r})}-\epsilon^{2}\Delta_{x}\theta_{r}
+(\theta^{N}+\theta_{r})^{4}-(\theta^{N})^{4}-\overline{f_{r}}=-\mathcal{L}_{2}(f^{N}, \theta^{N}),
\end{equation}
with initial conditions
\begin{equation}\label{rein}
f_{r}(0, \vec x, \vec w)=0,\ \ \theta_{r}(0, \vec x)=0,\ \ \mathrm{for}\ \ (\vec x, \vec w)\in\mathbb{T}^{3}\times\mathbb{S}^{2}.
\end{equation}

\begin{thm}\label{errorestimates}
Assume $2\leq N\leq 9$. The composite approximate solution $(f^{N}, \theta^{N})$ constructed in subsection \ref{construction}, satisfies \eqref{reeq1}-\eqref{reeq2} with initial conditions \eqref{rein}. Moreover, the error terms $\mathcal{L}_{1}(f^{N}, \theta^{N}), \mathcal{L}_{2}(f^{N}, \theta^{N})$ satisfy
\begin{equation}
\|\mathcal{L}_{1}(f^{N}, \theta^{N})\|_{L_{T}^{\infty}L_{\vec w}^{\infty}H_{\vec x}^{2}}, \|\mathcal{L}_{2}(f^{N}, \theta^{N})\|_{L_{T}^{\infty}H_{\vec x}^{2}}\leq C\epsilon^{N+1},
\end{equation}
where $C>0$ is a positive constant independent of $\epsilon$.
\end{thm}
\noindent\textbf{Proof.} We first consider $\mathcal{L}_{1}(f^{N}, \theta^{N})$. From \eqref{resi1-zaijia}, we have
\begin{equation}\label{resi1-zaijia}
\begin{aligned}
&\mathcal{L}_{1}(f^{N}, \theta^{N})\\=&\epsilon^{N+1}\partial_{t}f_{N-1}+\epsilon^{N+2}\partial_{t}f_{N}+\epsilon^{N+1}\vec{w}\cdot\nabla_{x}f_{N}
-\sum_{k=N+1}^{4N}\epsilon^{k}B_{k}(\theta_{0};\theta_{k})\\&-\sum_{k=0}^{N}{\epsilon}^{k}(B_{k}(\theta_{0}+\theta_{I, 0};\theta_{k}+\theta_{I, k})-B_{k}(\theta_{0};\theta_{k}))
+\sum_{k=0}^{N}{\epsilon}^{k}(B_{k}(A_{0}+\theta_{I, 0};A_{k}+\theta_{I, k})-B_{k}(A_{0};A_{k}))\\&-\sum_{k=N+1}^{4N}\epsilon^{k}B_{k}(\theta_{0}+\theta_{I, 0};\theta_{k}+\theta_{I, k})
+\epsilon^{N+1}(f_{I, N-1}-\overline{f_{I, N-1}})+\epsilon^{N+2}(f_{I, N}-\overline{f_{I, N}})\\&+\epsilon^{N+1}\vec{w}\cdot\nabla_{x}f_{I, N}:=\mathcal{L}_{11}+\mathcal{L}_{12}+\mathcal{L}_{13},
\end{aligned}
\end{equation}
where
\begin{equation}
\begin{aligned}
\mathcal{L}_{11}=\epsilon^{N+1}\partial_{t}f_{N-1}+\epsilon^{N+2}\partial_{t}f_{N}+\epsilon^{N+1}\vec{w}\cdot\nabla_{x}f_{N}
-\sum_{k=N+1}^{4N}\epsilon^{k}B_{k}(\theta_{0};\theta_{k}),
\end{aligned}
\end{equation}
\begin{equation}
\begin{aligned}
\mathcal{L}_{12}=&-\sum_{k=0}^{N}{\epsilon}^{k}(B_{k}(\theta_{0}+\theta_{I, 0};\theta_{k}+\theta_{I, k})-B_{k}(\theta_{0};\theta_{k}))
+\sum_{k=0}^{N}{\epsilon}^{k}(B_{k}(A_{0}+\theta_{I, 0};A_{k}+\theta_{I, k})\\&-B_{k}(A_{0};A_{k}))
\end{aligned}
\end{equation}
and
\begin{equation}
\begin{aligned}
\mathcal{L}_{13}=&-\sum_{k=N+1}^{4N}\epsilon^{k}B_{k}(\theta_{0}+\theta_{I, 0};\theta_{k}+\theta_{I, k})
+\epsilon^{N+1}(f_{I, N-1}-\overline{f_{I, N-1}})+\epsilon^{N+2}(f_{I, N}-\overline{f_{I, N}})\\&+\epsilon^{N+1}\vec{w}\cdot\nabla_{x}f_{I, N}.
\end{aligned}
\end{equation}
Collecting Theorem \ref{3.1}-\ref{3.6}, we have
\begin{equation}\label{eq1}
\|\mathcal{L}_{11}\|_{L_{T}^{\infty}L_{\vec w}^{\infty}H_{\vec x}^{2}}, \|\mathcal{L}_{13}\|_{L_{T}^{\infty}L_{\vec w}^{\infty}H_{\vec x}^{2}}\leq C\epsilon^{N+1}.
\end{equation}
For $\mathcal{L}_{12}$, from the definition of \eqref{defBk}, 
\begin{equation}
\begin{aligned}
\mathcal{L}_{12}=-\{(\sum_{k=0}^{N}\epsilon^{k}(\theta_{k}+\theta_{I, k}))^{4}-(\sum_{k=0}^{N}\epsilon^{k}\theta_{k})^{4}-(\sum_{k=0}^{N}\epsilon^{k}(A_{k}+\theta_{I, k}))^{4}+(\sum_{k=0}^{N}\epsilon^{k}A_{k})^{4}\}.
\end{aligned}
\end{equation}
Using the formula $c_{1}^{4}-c_{2}^{4}=(c_{1}-c_{2})(c_{1}+c_{2})(c_{1}^{2}+c_{2}^{2})$ and for $c_{1}-c_{2}=c_{3}-c_{4}=c_{0}$,
\begin{equation}
\begin{aligned}
(c_{1}^{4}-c_{2}^{4})-(c_{3}^{4}-c_{4}^{4})=c_{0}(c_{1}-c_{2})(2c_{1}^{2}+2c_{2}^{2}
+(c_{3}+c_{4})(c_{1}+c_{2}+c_{3}+c_{4})),
\end{aligned}
\end{equation}
with $c_{1}=\sum_{k=0}^{N}\epsilon^{k}(\theta_{k}+\theta_{I, k})$, $c_{2}=\sum_{k=0}^{N}\epsilon^{k}\theta_{k}$, 
$c_{3}=\sum_{k=0}^{N}\epsilon^{k}(A_{k}+\theta_{I, k})$, $c_{4}=\sum_{k=0}^{N}\epsilon^{k}A_{k}$, we obtain
\begin{equation}\label{RL}
\begin{aligned}
\mathcal{L}_{12}=(\sum_{k=0}^{N}\epsilon^{k}\theta_{I, k})(\sum_{k=0}^{N}\epsilon^{k}(\theta_{k}-A_{k}))(2c_{1}^{2}+2c_{2}^{2}
+(c_{3}+c_{4})(c_{1}+c_{2}+c_{3}+c_{4})).
\end{aligned}
\end{equation}
Due to the exponential decay estimates \eqref{ineq-theta00-jia1}, \eqref{ineq-theta00-jia11} and \eqref{ineq-thetam1}, we have
\begin{equation}
\|\sum_{k=0}^{N}\epsilon^{k}\theta_{I, k}\|_{L_{T}^{\infty}H_{\vec x}^{2}}\leq Ce^{-\frac{\sigma t}{\epsilon^{2}}}.
\end{equation}
Taylor's formula yields
\begin{equation}
\theta_{k}(t, \vec x)=\sum_{l=0}^{N-k}\frac{t^{l}}{l!}\partial_{t}^{l}\theta_{l}(0, \vec x)+\frac{\partial_{t}^{N-k+1}\theta_{k}(t', \vec x)}{(N-k+1)!}t^{N-k+1},
\end{equation}
with $t'\in [0, t]$.
Using the above formula and \eqref{formula1}, we get
\begin{equation}\label{relation}
\begin{aligned}
\sum_{k=0}^{N}\epsilon^{k}(\theta_{k}-A_{k})=&\sum_{k=0}^{N}\epsilon^{k}\Big(\sum_{l=0}^{N-k}\frac{t^{l}}{l!}\partial_{t}^{l}\theta_{l}(0, \vec x)+\frac{\partial_{t}^{N-k+1}\theta_{k}(t', \vec x)}{(N-k+1)!}t^{N-k+1}\Big)\\&-\sum_{k=0}^{N}\epsilon^{k}\sum_{l=0}^{k}\epsilon^{l}\frac{\tau^{l}}{l!}\partial_{t}^{l}\theta_{k-l}(0, \vec x).
\end{aligned}
\end{equation}
Using the formula
\begin{equation}
\begin{aligned}
\sum_{k=0}^{N}\sum_{l=0}^{k}g(l, k)=\sum_{l=0}^{N}\sum_{k=l}^{N}g(l, k)=\sum_{l=0}^{N}\sum_{s=0}^{N-l}g(l, l+s)=\sum_{k=0}^{N}\sum_{l=0}^{N-k}g(k, k+l),
\end{aligned}
\end{equation}
where we have taken substitutions $s=k-l$ and $l\rightarrow k, s\rightarrow l$ in the above second and third equality.

Then, taking $g(k, k+l)=\epsilon^{k}\frac{t^{l}}{l!}\partial_{t}^{l}\theta_{k}(0, \vec x)$, we get
\begin{equation}
g(l, k)=\epsilon^{l}\frac{t^{k-l}}{(k-l)!}\partial_{t}^{k-l}\theta_{l}(0, \vec x)=\epsilon^{2k-l}\frac{\tau^{k-l}}{(k-l)!}\partial_{t}^{k-l}\theta_{l}(0, \vec x)
\end{equation}
and so
\begin{equation}
\begin{aligned}
\sum_{k=0}^{N}\sum_{l=0}^{N-k}\epsilon^{k}\frac{t^{l}}{l!}\partial_{t}^{l}\theta_{k}(0, \vec x)=\sum_{k=0}^{N}\sum_{l=0}^{N-k}\epsilon^{2k-l}\frac{\tau^{k-l}}{(k-l)!}\partial_{t}^{k-l}\theta_{l}(0, \vec x).
\end{aligned}
\end{equation}
By substitutions $k\rightarrow k, k-l\rightarrow l$, we have
\begin{equation}
\begin{aligned}
\sum_{k=0}^{N}\sum_{l=0}^{N-k}\epsilon^{2k-l}\frac{\tau^{k-l}}{(k-l)!}\partial_{t}^{k-l}\theta_{l}(0, \vec x)=\sum_{k=0}^{N}\sum_{l=0}^{N}\epsilon^{k+l}\frac{\tau^{l}}{l!}\partial_{t}^{l}\theta_{k-l}(0, \vec x).
\end{aligned}
\end{equation}
Taking this relation into \eqref{relation} leads to
\begin{equation}\label{relation1}
\begin{aligned}
\sum_{k=0}^{N}\epsilon^{k}(\theta_{k}-A_{k})=&\sum_{k=0}^{N}\epsilon^{k}\frac{\partial_{t}^{N-k+1}\theta_{k}(t', \vec x)}{(N-k+1)!}t^{N-k+1}.
\end{aligned}
\end{equation}
Combining \eqref{relation1} with \eqref{relation}, \eqref{RL} satisfies
\begin{equation}
\begin{aligned}
\|\mathcal{L}_{12}\|_{L_{T}^{\infty}H_{\vec x}^{2}}\leq C\sum_{k=0}^{N}\epsilon^{k}\frac{1}{(N-k+1)!}t^{N-k+1}e^{-\frac{\sigma t}{\epsilon^{2}}},
\end{aligned}
\end{equation}
where we have used the fact that $\|c_{1}, c_{2}, c_{3}, c_{4}\|_{L_{T}^{\infty}H_{\vec x}^{2}}\leq C$ and $\|\theta_{k}\|_{C^{0}([0, T]; H_{\vec x}^{N+2-k})}\leq C$, $0\leq k\leq N$. Note that the function $h(t)=t^{N-k+1}e^{-\frac{\sigma t}{\epsilon^{2}}}$ attain its maximum at $t^{\star}=\frac{(N-k+1)\epsilon^{2}}{\sigma}$ with the maximum value $h(t^{\star})=(N-k+1)^{N-k+1}/\sigma\cdot e^{-(N-k+1)}$. Therefore,
\begin{equation}
\begin{aligned}
\|\mathcal{L}_{12}\|_{L_{T}^{\infty}H_{\vec x}^{2}}\leq C\sum_{k=0}^{N}\epsilon^{2N+2-k}\frac{(N-k+1)^{N-k+1}}{\sigma^{N-k+1}(N-k+1)!}e^{-(N-k+1)}\leq C\epsilon^{N+2}(\gamma_{N}-1),
\end{aligned}
\end{equation}
where $\gamma_{N}:=\sum_{n=0}^{N+1}n^{n}/(\sigma^{n}n!)>1$ is a constant depending on $N$. Therefore
\begin{equation}\label{eq2}
\|\mathcal{L}_{12}\|_{L_{T}^{\infty}H_{\vec x}^{2}}\leq C\epsilon^{N+2}.
\end{equation}
Collecting \eqref{eq1} and \eqref{eq2}, we can derive
\begin{equation}
\|\mathcal{L}_{1}\|_{L_{T}^{\infty}L_{\vec w}^{\infty}H_{\vec x}^{2}}\leq C\epsilon^{N+1}.
\end{equation}
Recalling \eqref{resi2-zaijia},
\begin{equation}\label{resi2-zaijia}
\begin{aligned}
&\mathcal{L}_{2}(f^{N}, \theta^{N})\\=&\mathcal{L}_{2}\Big(\sum_{k=0}^{N}{\epsilon}^{k}f_{k}, \sum_{k=0}^{N}\epsilon^{k}\theta_{k}\Big)+\sum_{k=0}^{N}{\epsilon}^{k}(B_{k}(\theta_{0}+\theta_{I, 0};\theta_{k}+\theta_{I, k})-B_{k}(\theta_{0};\theta_{k}))
\\&-\sum_{k=0}^{N}{\epsilon}^{k}(B_{k}(A_{0}+\theta_{I, 0};A_{k}+\theta_{I, k})-B_{k}(A_{0};A_{k}))+\sum_{k=N+1}^{4N}\epsilon^{k}B_{k}(\theta_{0}+\theta_{I, 0};\theta_{k}+\theta_{I, k})
\\&+\epsilon^{N+1}(\mathrm{div}(\vec u\theta_{I, k-1})-\Delta_{x}\theta_{I, k-1})+\epsilon^{N+2}(\mathrm{div}(\vec u\theta_{I, k})-\Delta_{x}\theta_{I, k})+\sum_{k=0}^{N}\epsilon^{k}\mathrm{div}((\vec u-\vec v_{k})\theta_{I, k-2}).
\end{aligned}
\end{equation}
Similarly, we can derive
\begin{equation}
\|\mathcal{L}_{2}\|_{L_{T}^{\infty}H_{\vec x}^{2}}\leq C\epsilon^{N+1}.
\end{equation}

In order to prove Theorem \ref{result}, we first derive suitable estimates on a linearized system and then use the Banach fixed point theorem to show the existence of the above problem near zero solutions, leading to the convergence of $(f^{\epsilon}, \theta^{\epsilon})$ to $(f^{N}, \theta^{N})$ as $\epsilon\rightarrow 0$.

\subsection{Linearized system}
We first linearize system \eqref{reeq1}-\eqref{reeq2} around zero and consider the following linear system:
\begin{equation}\label{reeq1-jia}
\epsilon^{2}\partial_{t}f_{r}+\epsilon\vec{w}\cdot\nabla_{x}f_{r}+\epsilon^{2}(f_{r}-\overline{f_{r}})+f_{r}
-4(\theta^{N})^{3}\theta_{r}=R_{1}-R,
\end{equation}
\begin{equation}\label{reeq2-jia}
\epsilon^{2}\partial_{t}\theta_{r}+\epsilon^{2}\mathrm{div}{(\vec u\theta_{r})}-\epsilon^{2}\Delta_{x}\theta_{r}
+4(\theta^{N})^{3}\theta_{r}-\overline{f_{r}}=R_{2}+\langle R\rangle,
\end{equation}
where $R_{1}=R_{1}(t, \vec x, \vec w)$, $R=R(t, \vec x, \vec w)$ and $R_{2}=R_{2}(t, \vec x)$ are given functions and the initial conditions are taken to be 
\begin{equation}\label{rein-jia}
f_{r}(0, \vec x, \vec w)=0,\ \ \theta_{r}(0, \vec x)=0,\ \ \mathrm{for}\ \ (\vec x, \vec w)\in\mathbb{T}^{3}\times\mathbb{S}^{2}.
\end{equation}
The existence of system \eqref{reeq1-jia}-\eqref{rein-jia} will be given later.
\begin{lem}
Let $\epsilon>0$ and $(f^{N}, \theta^{N})$ be the composite approximate solution construct in section \ref{san}. Assume $R_{1}, R\in L_{T}^{2}L_{\vec w}^{2}H_{\vec x}^{2}$ and $R_{2}\in L_{T}^{2}H_{\vec x}^{2}$. Then, there exists a unique solution $(f_{r}, \theta_{r})\in L_{T}^{\infty}L_{\vec x}^{\infty}L_{\vec w}^{\infty}\times L_{T}^{\infty}L_{\vec x}^{\infty}$ to system \eqref{reeq1-jia}-\eqref{reeq2-jia}. Moreover, the solution $(f_{r}, \theta_{r})$ satisfies the following estimates
\begin{equation}\label{bd00}
\begin{aligned}
&\|f_{r}\|_{L_{T}^{\infty}H^{2}_{\vec x}L_{\vec w}^{2}}+\|\theta_{r}\|_{L_{T}^{\infty}H_{\vec x}^{2}}+\frac{1}{\epsilon^{2}}\|f_{r}-4(\theta_{0}^{3}+3\theta_{0}^{2}\theta_{1}\epsilon)\theta_{r}\|_{L_{T}^{2}H_{\vec x}^{2}L_{\vec w}^{2}}+\|\theta_{r}\|_{L_{T}^{2}H_{\vec x}^{3}}
\\\leq& \frac{C(T)}{\epsilon^{4}}\Big(\|R\|_{L_{T}^{2}H_{\vec x}^{2}L_{\vec w}^{2}}+\|R_{1}\|_{L_{T}^{2}H_{\vec x}^{2}L_{\vec w}^{2}}+\|R_{2}\|_{L_{T}^{2}H_{\vec x}^{2}}\Big),
\end{aligned}
\end{equation}
where $C$ is a constant depending on $\kappa$ not depending on $\epsilon$. Here $\kappa$ is a suitably small constant.
\end{lem}
\noindent\textbf{Proof.} We first get the $L^{2}$ type estimate. Finally, the $L^{\infty}$ type estimate is shown.

\noindent\textbf{Step 1.} Dividing $(\theta^{N})^{3}$ into three parts. According to \eqref{expansion of theta epsilon}, we have 
\begin{equation}
\begin{aligned}
(\theta^{N})^{3}=&\theta_{0}^{3}+3\theta_{0}^{2}\theta_{1}\epsilon+\theta_{I, 0}^{3}+3\theta_{I, 1}(\theta_{0}+\theta_{I, 0})^{2}\epsilon+3\theta_{I, 0}^{2}(\theta_{1}+\theta_{I, 1})\epsilon+O(\epsilon^{2})
\\:=&L_{0}+L_{1}+L_{2},
\end{aligned}
\end{equation}
where
\begin{equation}
L_{0}=\theta_{0}^{3}+3\theta_{0}^{2}\theta_{1}\epsilon,
\end{equation}
\begin{equation}
L_{1}=\theta_{I, 0}^{3}+3\theta_{I, 1}(\theta_{0}+\theta_{I, 0})^{2}\epsilon+3\theta_{I, 0}^{2}(\theta_{1}+\theta_{I, 1})\epsilon
\end{equation}
and
\begin{equation}
L_{2}=O(\epsilon^{2}).
\end{equation}
So, we can rewrite \eqref{reeq1-jia} and \eqref{reeq2-jia} as follows
\begin{equation}\label{reeq1-jia-jia}
\epsilon^{2}\partial_{t}f_{r}+\epsilon\vec{w}\cdot\nabla_{x}f_{r}+\epsilon^{2}(f_{r}-\overline{f_{r}})+f_{r}
-4L_{0}\theta_{r}-4L_{1}\theta_{r}-4L_{2}\theta_{r}=R_{1}-R,
\end{equation}
\begin{equation}\label{reeq2-jia-jia}
\epsilon^{2}\partial_{t}\theta_{r}+\epsilon^{2}\mathrm{div}{(\vec u\theta_{r})}-\epsilon^{2}\Delta_{x}\theta_{r}
+4L_{0}\theta_{r}-\overline{f_{r}}+4L_{1}\theta_{r}+4L_{2}\theta_{r}=R_{2}+\langle R\rangle.
\end{equation}

\noindent\textbf{Step 2.} The well-posedness of the linearized system \eqref{reeq1-jia}, \eqref{reeq2-jia} and \eqref{rein-jia}.

For a given $\theta_{r}'\in L^{2}((0, T)\times\Omega)$, we consider the following problem
\begin{equation}\label{research equation-2-1-1}\left\{
\begin{split}
&\epsilon^{2}\partial_{t}f_{r}+\epsilon\vec{w}\cdot\nabla_{x}f_{r}+\epsilon^{2}(f_{r}-\overline{f_{r}})+f_{r}
-4(\theta^{N})^{3}\theta_{r}'=R_{1}-R,\ \ \mathrm{in}\ \ (0, T]\times\Omega,\\
&f_{r}(0, \vec{x}, \vec w)=0\ \ \mathrm{in}\ \ \Omega.\\
\end{split}\right.
\end{equation}
According to Lemma \ref{lianxujieguo}, we obtain that there exists a unique solution $f_{r}\in C^{0}([0, T]; L^{2}(\mathbb{S}^{2}; \\L^{2}(\Omega)))$ of (\ref{research equation-2-1-1}). Define the map
\begin{equation}\nonumber
S:\theta_{r}'\in L^{2}((0, T)\times\Omega)\mapsto f_{r}\in C^{0}([0, T]; L^{2}(\mathbb{S}^{2}; L^{2}(\Omega))).
\end{equation}
Similarly, we consider the parabolic equation
\begin{equation}\label{research equations-2-1-1}\left\{
\begin{split}
&\epsilon^{2}\partial_{t}\theta_{r}^{\tilde\sigma}+\epsilon^{2}\mathrm{div}{(\vec u\theta_{r}^{\tilde\sigma})}-\epsilon^{2}\Delta_{x}\theta_{r}^{\tilde\sigma}
+4(\theta^{N})^{3}\theta_{r}^{\tilde\sigma}=\tilde\sigma(S(\overline {\theta_{r}}')+R_{2}+\langle R\rangle),\\
&\theta_{r}^{\tilde\sigma}(0, \vec{x})=0,\ \ \mathrm{in}  \ \ \Omega.\\
\end{split}\right.
\end{equation}
For sufficiently small $\epsilon$, we can claim that $\theta^{N}>0$. According to the linear parabolic theory, we can obtain that there exists a unique solution $\theta_{r}^{\tilde\sigma}\in L^{2}(0, T; H_{\vec x}^{2})\cap L^{\infty}(0, T; H_{\vec x}^{1})$ for any given $\theta_{r}'\in L^{2}((0, T)\times\Omega), \tilde\sigma\in[0, 1].$ Note that $\partial_{t}\theta_{r}^{\tilde\sigma}\in L^{2}((0, T)\times\Omega)$. According to the Aubin-Lions-Simon lemma, we can define a compact map
\begin{equation}\nonumber
\begin{aligned}
&\Gamma_{\tilde\sigma}:\theta_{r}'\in L^{2}((0, T)\times\Omega)\mapsto&\theta_{r}^{\tilde\sigma}\in \{\theta|\theta\in L^{2}(0, T; H_{\vec x}^{2}),\partial_{t}\theta\in L^{2}((0, T)\times\Omega)\}.
\end{aligned}
\end{equation}

Now, we assume that $\overline {f_{1}^{\tilde\sigma}}$ and $\theta_{1}=S(\overline {f^{\tilde\sigma}_{1}})$ satisfy
\begin{equation}\label{research equations-2-1-2}\left\{
\begin{split}
&\epsilon^{2}\partial_{t}\theta_{r}^{\sigma}+\epsilon^{2}\mathrm{div}{(\vec u\theta_{r}^{\tilde\sigma})}-\epsilon^{2}\Delta_{x}\theta_{r}^{\tilde\sigma}
+4(\theta^{N})^{3}\theta_{r}^{\tilde\sigma}=\tilde\sigma(\overline {f_{r}}+R_{2}+\langle R\rangle),\\
&\epsilon^{2}\partial_{t}f_{r}+\epsilon\vec{w}\cdot\nabla_{x}f_{r}+\epsilon^{2}(f_{r}-\overline{f_{r}})+f_{r}
-4(\theta^{N})^{3}\theta_{r}^{\tilde\sigma}=R_{1}-R,\ \ \mathrm{in}\ \ (0, T]\times\Omega,\\
&\theta_{r}^{\tilde\sigma}(0, \vec{x})=0,\mathrm{in}\ \ \Omega,\\
&f_{r}(0, \vec{x}, \vec w)=0,\ \ \mathrm{in}\ \ \Omega.
\end{split}\right.
\end{equation}
Multiplying $\theta_{r}^{\tilde\sigma}$ and $f_{r}$ on both sides of equation $(\ref{research equations-2-1-2})_{1}$ and $(\ref{research equations-2-1-2})_{2}$ respectively, for small $\delta$, we derive
\begin{equation}\nonumber
\begin{split}
&\epsilon^{2}\int_{\Omega}|\theta_{r}^{\tilde\sigma}|^{2}\mathrm{d}\vec{x}(t)+\epsilon^{2}\int_{\mathbb{S}^{2}}\int_{\Omega}|f_{r}|^{2}\mathrm{d}\vec{x}(t)
+\epsilon^{2}\int_{0}^{t}\int_{\Omega}|\nabla_{x}\theta_{r}^{\tilde\sigma}|^{2}\mathrm{d}\vec{x}\mathrm{d}\vec{w}\mathrm{d}s
+\int_{0}^{t}\int_{\Omega}4(\theta^{N})^{3}|\theta_{r}^{\tilde\sigma}|^{2}\mathrm{d}\vec{x}\mathrm{d}s
\\&+\epsilon^{2}\int_{0}^{t}\int_{\Omega}\int_{\mathbb{S}^{2}}|f_{r}-\overline{f_{r}}|^{2}\mathrm{d}\vec{w}\mathrm{d}\vec{x}\mathrm{d}s
\end{split}
\end{equation}
\begin{equation}\nonumber
\begin{split}
\\\leq& C(\tilde\sigma+1)\Big(\int_{0}^{t}\int_{\mathbb{S}^{2}}\int_{\Omega}|f_{r}|^{2}\mathrm{d}\vec{x}\mathrm{d}\vec{w}\mathrm{d}s
+\int_{0}^{t}\int_{\Omega}|\theta_{r}^{\tilde\sigma}|^{2}\mathrm{d}\vec{x}\mathrm{d}s
+\delta\epsilon^{2}\int_{0}^{t}\int_{\Omega}|\nabla_{x}\theta_{r}^{\tilde\sigma}|^{2}\mathrm{d}\vec{w}\mathrm{d}\vec{x}\mathrm{d}s \\&+\int_{0}^{t}\int_{\Omega}\int_{\mathbb{S}^{2}}|R|^{2}\mathrm{d}\vec{w}\mathrm{d}\vec{x}\mathrm{d}s
+\int_{0}^{t}\int_{\Omega}\int_{\mathbb{S}^{2}}|R_{1}|^{2}\mathrm{d}\vec{w}\mathrm{d}\vec{x}\mathrm{d}s
+\int_{0}^{t}\int_{\Omega}|R_{2}|^{2}\mathrm{d}\vec{x}\mathrm{d}s\Big).
\end{split}
\end{equation}
Since $s\in (0, 1)$, we have
\begin{equation}\nonumber
\begin{aligned}
&\epsilon^{2}\int_{\Omega}|\theta_{r}^{\tilde\sigma}|^{2}\mathrm{d}\vec{x}(t)+\epsilon^{2}\int_{\mathbb{S}^{2}}\int_{\Omega}|f_{r}|^{2}\mathrm{d}\vec{x}(t)
+\epsilon^{2}\int_{0}^{t}\int_{\Omega}|\nabla_{x}\theta_{r}^{\tilde\sigma}|^{2}\mathrm{d}\vec{x}\mathrm{d}\vec{w}\mathrm{d}s
+\int_{0}^{t}\int_{\Omega}4(\theta^{N})^{3}|\theta_{r}^{\tilde\sigma}|^{2}\mathrm{d}\vec{x}\mathrm{d}s
\\&+\epsilon^{2}\int_{0}^{t}\int_{\Omega}\int_{\mathbb{S}^{2}}|f_{r}-\overline{f_{r}}|^{2}\mathrm{d}\vec{w}\mathrm{d}\vec{x}\mathrm{d}s
\\\leq& C\Big(\int_{0}^{t}\int_{\mathbb{S}^{2}}\int_{\Omega}|f_{r}|^{2}\mathrm{d}\vec{x}\mathrm{d}\vec{w}\mathrm{d}s
+\int_{0}^{t}\int_{\Omega}|\theta_{r}^{\tilde\sigma}|^{2}\mathrm{d}\vec{x}\mathrm{d}s+\int_{0}^{t}\int_{\Omega}\int_{\mathbb{S}^{2}}|R|^{2}\mathrm{d}\vec{w}\mathrm{d}\vec{x}\mathrm{d}s
\\&+\int_{0}^{t}\int_{\Omega}\int_{\mathbb{S}^{2}}|R_{1}|^{2}\mathrm{d}\vec{w}\mathrm{d}\vec{x}\mathrm{d}s
+\int_{0}^{t}\int_{\Omega}|R_{2}|^{2}\mathrm{d}\vec{x}\mathrm{d}s\Big).
\end{aligned}
\end{equation}
According to the Gronwall's inequality, we obtain $\theta_{r}^{\tilde\sigma}\in L^{\infty}_{T}L^{2}_{\vec x}\cap L^{2}_{T}H^{1}_{\vec x}$ and $f_{r}\in L^{\infty}_{T}L^{2}_{\vec w}L^{2}_{\vec x}\cap L^{2}_{T}L^{2}_{\vec w}L^{2}_{\vec x}$. Note that $R_{1}, R\in L_{T}^{2}L_{\vec w}^{2}H_{\vec x}^{2}$ and $R_{2}\in L_{T}^{2}H_{\vec x}^{2}$. After a bootstrap, we have $\theta_{r}^{\tilde\sigma}\in L^{\infty}_{T}H^{2}_{\vec x}\cap L^{2}_{T}H^{3}_{\vec x}$, $f_{r}\in L^{\infty}_{T}L^{2}_{\vec w}H^{2}_{\vec x}\cap L^{2}_{T}L^{2}_{\vec w}H^{2}_{\vec x}$ and $\partial_{t}\theta_{r}^{\tilde\sigma}\in L^{\infty}_{T}L^{2}_{\vec x}\cap L^{2}_{T}H^{1}_{\vec x}$, and these estimates are independent of $s$. Therefore, we can obtain that there exists a solution $(\theta_{r}^{\tilde\sigma}, f_{r})$ of problem (\ref{research equations-2-1-2}) such that $\theta_{r}^{\tilde\sigma}\in L^{\infty}_{T}H^{2}_{\vec x}\cap L^{2}_{T}H^{3}_{\vec x}$ and $\theta_{r}^{\tilde\sigma}=S(f_{r})\in L^{\infty}_{T}L^{2}_{\vec w}H^{2}_{\vec x}\cap L^{2}_{T}L^{2}_{\vec w}H^{2}_{\vec x}$ according to the Schaefer fixed point theorem.

\noindent\textbf{Step 3.} The energy estimate. We apply the operator $D_{x}^{\alpha}$, $0\leq\alpha\leq 2$ on both sides of \eqref{reeq1-jia-jia} and \eqref{reeq2-jia-jia} and multiply \eqref{reeq1-jia-jia} by $D_{x}^{\alpha}f_{r}$ and \eqref{reeq2-jia-jia} by $D_{x}^{\alpha}(4L_{0}\theta_{r})$, and integrate over $s\in(0, t)$, $\vec x\in\mathbb{T}^{3}$ and $\vec w\in\mathbb{S}^{2}$ to get
\begin{equation}
\begin{aligned}
&\epsilon^{2}\int_{\mathbb{T}^{3}}\int_{\mathbb{S}^{2}}|D_{x}^{\alpha}f_{r}|^{2}(t)\mathrm{d}\vec w\mathrm{d}\vec x
+\epsilon^{2}\int_{\mathbb{T}^{3}}\int_{\mathbb{S}^{2}}4L_{0}|D_{x}^{\alpha}\theta_{r}|^{2}(t)
\mathrm{d}\vec w\mathrm{d}\vec x
\\&+\epsilon^{2}\int_{0}^{t}\int_{\mathbb{T}^{3}}\int_{\mathbb{S}^{2}}|D_{x}^{\alpha}f_{r}-D_{x}^{\alpha}\overline{f_{r}}|^{2}
\mathrm{d}\vec w\mathrm{d}\vec x\mathrm{d}s
+\int_{0}^{t}\int_{\mathbb{T}^{3}}\int_{\mathbb{S}^{2}}|D_{x}^{\alpha}(4L_{0}\theta_{r})
-D_{x}^{\alpha}f_{r}|^{2}\mathrm{d}\vec w\mathrm{d}\vec x\mathrm{d}s\\&+\epsilon^{2}\int_{0}^{t}\int_{\mathbb{T}^{3}}\int_{\mathbb{S}^{2}}4L_{0}|\nabla_{x}D_{x}^{\alpha}\theta_{r}
|^{2}\mathrm{d}\vec w\mathrm{d}\vec x\mathrm{d}s
\\=&-\epsilon^{2}\int_{0}^{t}\int_{\mathbb{T}^{3}}\int_{\mathbb{S}^{2}}D_{x}^{\alpha}\mathrm{div}{(\vec u\theta_{r})}D_{x}^{\alpha}(4L_{0}\theta_{r})\mathrm{d}\vec w\mathrm{d}\vec x\mathrm{d}s
+\epsilon^{2}\int_{0}^{t}\int_{\mathbb{T}^{3}}\int_{\mathbb{S}^{2}}\partial_{t}(4L_{0})|D_{x}^{\alpha}\theta_{r}|^{2}(t)
\mathrm{d}\vec w\mathrm{d}\vec x\mathrm{d}s\\&-\epsilon^{2}\int_{0}^{t}\int_{\mathbb{T}^{3}}\int_{\mathbb{S}^{2}}[D_{x}^{\alpha}, 4L_{0}]\theta_{r}\partial_{t}D_{x}^{\alpha}\theta_{r}\mathrm{d}\vec w\mathrm{d}\vec x\mathrm{d}t-\epsilon^{2}\int_{0}^{t}\int_{\mathbb{T}^{3}}\int_{\mathbb{S}^{2}}
\nabla_{x}(4L_{0})D_{x}^{\alpha}\theta_{r}\nabla_{x}D_{x}^{\alpha}\theta_{r}
\mathrm{d}\vec w\mathrm{d}\vec x\mathrm{d}t\\&-\epsilon^{2}\int_{0}^{t}\int_{\mathbb{T}^{3}}\int_{\mathbb{S}^{2}}[D_{x}^{\alpha}, 4L_{0}]\theta_{r}\Delta_{x}D_{x}^{\alpha}\theta_{r}\mathrm{d}\vec w\mathrm{d}\vec x\mathrm{d}s+\int_{0}^{t}\int_{\mathbb{T}^{3}}\int_{\mathbb{S}^{2}}D_{x}^{\alpha}(4L_{1}\theta_{r})D_{x}^{\alpha}f_{r}\mathrm{d}\vec w\mathrm{d}\vec x\mathrm{d}s\\&+\int_{0}^{t}\int_{\mathbb{T}^{3}}\int_{\mathbb{S}^{2}}D_{x}^{\alpha}(4L_{2}\theta_{r})D_{x}^{\alpha}f_{r}\mathrm{d}\vec w\mathrm{d}\vec x\mathrm{d}s-\int_{0}^{t}\int_{\mathbb{T}^{3}}\int_{\mathbb{S}^{2}}D_{x}^{\alpha}(4L_{1}\theta_{r})D_{x}^{\alpha}(4L_{0}\theta_{r})\mathrm{d}\vec w\mathrm{d}\vec x\mathrm{d}s
\end{aligned}
\end{equation}
\begin{equation}\label{bd0}
\begin{aligned}
\\&-\int_{0}^{t}\int_{\mathbb{T}^{3}}\int_{\mathbb{S}^{2}}D_{x}^{\alpha}(4L_{2}\theta_{r})D_{x}^{\alpha}(4L_{0}\theta_{r})\mathrm{d}\vec w\mathrm{d}\vec x\mathrm{d}s
+\int_{0}^{t}\int_{\mathbb{T}^{3}}\int_{\mathbb{S}^{2}}D_{x}^{\alpha}(R_{1}-R)D_{x}^{\alpha}f_{r}\mathrm{d}\vec w\mathrm{d}\vec x\mathrm{d}t\\&+\int_{0}^{t}\int_{\mathbb{T}^{3}}\int_{\mathbb{S}^{2}}D_{x}^{\alpha}(R_{2}+\langle R\rangle)D_{x}^{\alpha}(4L_{0}\theta_{r})\mathrm{d}\vec w\mathrm{d}\vec x\mathrm{d}s
=:\sum_{i=1}^{11}H_{i},
\end{aligned}
\end{equation}
where
\begin{equation}\label{bd1}
\begin{aligned}
H_{1}=&-\epsilon^{2}\int_{0}^{t}\int_{\mathbb{T}^{3}}\int_{\mathbb{S}^{2}}D_{x}^{\alpha}\mathrm{div}{(\vec u\theta_{r})}D_{x}^{\alpha}(4(\theta^{N})^{3}\theta_{r})\mathrm{d}\vec w\mathrm{d}\vec x\mathrm{d}s
\\\leq& C(\delta)\epsilon^{2}\|\vec u\|_{C_{t, \vec x}^{3}}^{2}\|\theta^{N}\|_{L_{t}^{\infty}H_{\vec x}^{2}}^{6}\int_{0}^{t}\int_{\mathbb{T}^{3}}\int_{\mathbb{S}^{2}}|D_{x}^{\alpha}\theta_{r}|^{2}\mathrm{d}\vec w\mathrm{d}\vec x\mathrm{d}s+\delta\epsilon^{2}\int_{0}^{t}\int_{\mathbb{T}^{3}}\int_{\mathbb{S}^{2}}|D_{x}^{\alpha+1}\theta_{r}|^{2}\mathrm{d}\vec w\mathrm{d}\vec x\mathrm{d}s,
\end{aligned}
\end{equation}

\begin{equation}\label{bd2}
\begin{aligned}
H_{2}=&\epsilon^{2}\int_{0}^{t}\int_{\mathbb{T}^{3}}\int_{\mathbb{S}^{2}}\partial_{t}(4L_{0})|D_{x}^{\alpha}\theta_{r}|^{2}(t)
\mathrm{d}\vec w\mathrm{d}\vec x\mathrm{d}s
\\\leq& C\epsilon^{2}\int_{0}^{t}\int_{\mathbb{T}^{3}}\int_{\mathbb{S}^{2}}|D_{x}^{\alpha}\theta_{r}|^{2}\mathrm{d}\vec w\mathrm{d}\vec x\mathrm{d}s,
\end{aligned}
\end{equation}

\begin{equation}\label{bd3}
\begin{aligned}
H_{4}=&-\epsilon^{2}\int_{0}^{t}\int_{\mathbb{T}^{3}}\int_{\mathbb{S}^{2}}\nabla_{x}(4(\theta^{N})^{3})D_{x}^{\alpha}\theta_{r}\nabla_{x}D_{x}^{\alpha}\theta_{r}
\mathrm{d}\vec w\mathrm{d}\vec x\mathrm{d}s
\\\leq& C(\delta)\epsilon^{2}\|\theta^{N}\|_{L_{t}^{\infty}H_{\vec x}^{2}}^{6}\int_{0}^{t}\int_{\mathbb{T}^{3}}\int_{\mathbb{S}^{2}}|D_{x}^{\alpha}\theta_{r}|^{2}\mathrm{d}\vec w\mathrm{d}\vec x\mathrm{d}t+\delta\epsilon^{2}\int_{0}^{t}\int_{\mathbb{T}^{3}}\int_{\mathbb{S}^{2}}|D_{x}^{\alpha+1}\theta_{r}|^{2}\mathrm{d}\vec w\mathrm{d}\vec x\mathrm{d}s,
\end{aligned}
\end{equation}

\begin{equation}\label{bd4}
\begin{aligned}
H_{5}=&-\epsilon^{2}\int_{0}^{t}\int_{\mathbb{T}^{3}}\int_{\mathbb{S}^{2}}[D_{x}^{\alpha}, 4(\theta^{N})^{3}]\theta_{r}\Delta_{x}D_{x}^{\alpha}\theta_{r}\mathrm{d}\vec w\mathrm{d}\vec x\mathrm{d}s
\\=&\epsilon^{2}\int_{0}^{t}\int_{\mathbb{T}^{3}}\int_{\mathbb{S}^{2}}\nabla_{x}([D_{x}^{\alpha}, 4(\theta^{N})^{3}]\theta_{r})\cdot\nabla_{x}D_{x}^{\alpha}\theta_{r}\mathrm{d}\vec w\mathrm{d}\vec x\mathrm{d}s
\\\leq& C(\delta)\epsilon^{2}\|\theta^{N}\|_{L_{t}^{\infty}H_{\vec x}^{2}}^{6}\int_{0}^{t}\int_{\mathbb{T}^{3}}\int_{\mathbb{S}^{2}}|D_{x}^{\alpha}\theta_{r}|^{2}\mathrm{d}\vec w\mathrm{d}\vec x\mathrm{d}s+\delta\epsilon^{2}\int_{0}^{t}\int_{\mathbb{T}^{3}}\int_{\mathbb{S}^{2}}|D_{x}^{\alpha+1}\theta_{r}|^{2}\mathrm{d}\vec w\mathrm{d}\vec x\mathrm{d}s,
\end{aligned}
\end{equation}

\begin{equation}\label{bd5}
\begin{aligned}
H_{10}=&\int_{0}^{t}\int_{\mathbb{T}^{3}}\int_{\mathbb{S}^{2}}D_{x}^{\alpha}(R_{1}-R)D_{x}^{\alpha}f_{r}\mathrm{d}\vec w\mathrm{d}\vec x\mathrm{d}s
\\\leq& \frac{1}{2\epsilon^{2}}\int_{0}^{t}\int_{\mathbb{T}^{3}}\int_{\mathbb{S}^{2}}|D_{x}^{\alpha}(R_{1}-R)|^{2}\mathrm{d}\vec w\mathrm{d}\vec x\mathrm{d}s+\frac{\epsilon^{2}}{2}\int_{0}^{t}\int_{\mathbb{T}^{3}}\int_{\mathbb{S}^{2}}|D_{x}^{\alpha}f_{r}|^{2}\mathrm{d}\vec w\mathrm{d}\vec x\mathrm{d}s
\end{aligned}
\end{equation}
and
\begin{equation}\label{bd6}
\begin{aligned}
H_{11}=&\int_{0}^{t}\int_{\mathbb{T}^{3}}\int_{\mathbb{S}^{2}}D_{x}^{\alpha}(R_{2}+\langle R\rangle)D_{x}^{\alpha}(4(\theta^{N})^{3}\theta_{r})\mathrm{d}\vec w\mathrm{d}\vec x\mathrm{d}s
\\\leq& C\|\theta^{N}\|_{L_{t}^{\infty}H_{\vec x}^{2}}^{6}\Big(\frac{1}{2\epsilon^{2}}\Big)\int_{0}^{t}\int_{\mathbb{T}^{3}}\int_{\mathbb{S}^{2}}|D_{x}^{\alpha}(R_{2}+\langle R\rangle)|^{2}\mathrm{d}\vec w\mathrm{d}\vec x\mathrm{d}s+\frac{\epsilon^{2}}{2}\int_{0}^{t}(\|\theta_{r}\|_{L_{\vec x}^{2}}^{2}+\|D_{x}^{\alpha}\theta_{r}\|_{L_{\vec x}^{2}}^{2})\mathrm{d}s.
\end{aligned}
\end{equation}
Since \eqref{reeq2-jia-jia}, we have
\begin{equation}\label{bd7}
\begin{aligned}
H_{3}=&-\epsilon^{2}\int_{0}^{t}\int_{\mathbb{T}^{3}}\int_{\mathbb{S}^{2}}[D_{x}^{\alpha}, 4L_{0}]\theta_{r}\partial_{t}D_{x}^{\alpha}\theta_{r}\mathrm{d}\vec w\mathrm{d}\vec x\mathrm{d}s
\\=&\int_{0}^{t}\int_{\mathbb{T}^{3}}\int_{\mathbb{S}^{2}}[D_{x}^{\alpha}, 4L_{0}]\theta_{r}\{\epsilon^{2}D_{x}^{\alpha}\mathrm{div}(\vec u\theta_{r})-\epsilon^{2}\Delta_{x}D_{x}^{\alpha}\theta_{r}
+(D_{x}^{\alpha}(4L_{0}\theta_{r})-D_{x}^{\alpha}\overline{f_{r}})
\\&+D_{x}^{\alpha}
(4L_{1}\theta_{r}+4L_{2}\theta_{r})-D_{x}^{\alpha}R_{2}-D_{x}^{\alpha}\langle R\rangle\}\mathrm{d}\vec w\mathrm{d}\vec x\mathrm{d}s
\\:=&H_{3}^{1}+H_{3}^{2}+H_{3}^{3}+H_{3}^{4}+H_{3}^{5}+H_{3}^{6},
\end{aligned}
\end{equation}
where
\begin{equation}\label{bd8}
\begin{aligned}
H_{3}^{1}=&\epsilon^{2}\int_{0}^{t}\int_{\mathbb{T}^{3}}\int_{\mathbb{S}^{2}}[D_{x}^{\alpha}, 4L_{0}]\theta_{r}D_{x}^{\alpha}\mathrm{div}(\vec u\theta_{r})\mathrm{d}\vec w\mathrm{d}\vec x\mathrm{d}s
\\=&-\epsilon^{2}\int_{0}^{t}\int_{\mathbb{T}^{3}}\int_{\mathbb{S}^{2}}D_{x}[D_{x}^{\alpha}, 4L_{0}]\theta_{r}D_{x}^{\alpha-1}\mathrm{div}(\vec u\theta_{r})\mathrm{d}\vec w\mathrm{d}\vec x\mathrm{d}s
\\\leq& C\epsilon^{2}\int_{0}^{t}\int_{\mathbb{T}^{3}}\int_{\mathbb{S}^{2}}(|D_{x}^{\alpha}\theta_{r}|^{2}
+|\theta_{r}|^{2})\mathrm{d}\vec w\mathrm{d}\vec x\mathrm{d}s,
\end{aligned}
\end{equation}

\begin{equation}\label{bd9}
\begin{aligned}
H_{3}^{2}=&-\epsilon^{2}\int_{0}^{t}\int_{\mathbb{T}^{3}}\int_{\mathbb{S}^{2}}[D_{x}^{\alpha}, 4L_{0}]\theta_{r}D_{x}^{\alpha}\Delta_{x}\theta_{r}\mathrm{d}\vec w\mathrm{d}\vec x\mathrm{d}s
\\=&-\epsilon^{2}\int_{0}^{t}\int_{\mathbb{T}^{3}}\int_{\mathbb{S}^{2}}\nabla_{x}[D_{x}^{\alpha}, 4L_{0}]\theta_{r}\nabla_{x}D_{x}^{\alpha}\theta_{r}\mathrm{d}\vec w\mathrm{d}\vec x\mathrm{d}s
\\\leq& C(\delta)\epsilon^{2}\int_{0}^{t}\int_{\mathbb{T}^{3}}\int_{\mathbb{S}^{2}}|D_{x}^{\alpha}\theta_{r}|^{2}
\mathrm{d}\vec w\mathrm{d}\vec x\mathrm{d}s+\delta\epsilon^{2}\int_{0}^{t}\int_{\mathbb{T}^{3}}\int_{\mathbb{S}^{2}}|D_{x}^{\alpha+1}\theta_{r}|^{2}
\mathrm{d}\vec w\mathrm{d}\vec x\mathrm{d}s,
\end{aligned}
\end{equation}

\begin{equation}\label{bd10}
\begin{aligned}
H_{3}^{3}=&\int_{0}^{t}\int_{\mathbb{T}^{3}}\int_{\mathbb{S}^{2}}[D_{x}^{\alpha}, 4L_{0}]\theta_{r}(D_{x}^{\alpha}(4L_{0}\theta_{r})-D_{x}^{\alpha}\overline{f_{r}})\mathrm{d}\vec w\mathrm{d}\vec x\mathrm{d}s
\\\leq& C(\delta)\int_{0}^{t}\int_{\mathbb{T}^{3}}\int_{\mathbb{S}^{2}}(|D_{x}^{\alpha-1}\theta_{r}|^{2}
+|D_{x}^{\alpha-2}\theta_{r}|^{2})
\mathrm{d}\vec w\mathrm{d}\vec x\mathrm{d}s
\\&+\delta\int_{0}^{t}\int_{\mathbb{T}^{3}}\int_{\mathbb{S}^{2}}
|D_{x}^{\alpha}(4L_{0}\theta_{r})-D_{x}^{\alpha}f_{r}|^{2}
\mathrm{d}\vec w\mathrm{d}\vec x\mathrm{d}s,
\end{aligned}
\end{equation}

\begin{equation}\label{bd11}
\begin{aligned}
H_{3}^{4}=&\int_{0}^{t}\int_{\mathbb{T}^{3}}\int_{\mathbb{S}^{2}}[D_{x}^{\alpha}, 4L_{0}]\theta_{r}D_{x}^{\alpha}
(4L_{1}\theta_{r}+4L_{2}\theta_{r})\mathrm{d}\vec w\mathrm{d}\vec x\mathrm{d}s
\\\leq& C(\|\theta_{r}\|_{L_{t}^{\infty}L_{\vec x}^{2}}^{2}+\|D_{x}^{\alpha}\theta_{r}\|_{L_{t}^{\infty}L_{\vec x}^{2}}^{2})\int_{0}^{\frac{t}{\epsilon^{2}}}\|L_{1}\|_{H_{\vec x}^{2}}\epsilon^{2}\mathrm{d}\tau+C\epsilon^{2}\int_{0}^{t}(\|\theta_{r}\|_{L_{\vec x}^{2}}^{2}+\|D_{x}^{\alpha}\theta_{r}\|_{L_{\vec x}^{2}}^{2})\mathrm{d}s
\\\leq& C\eta\epsilon^{2}(\|\theta_{r}\|_{L_{t}^{\infty}L_{\vec x}^{2}}^{2}+\|D_{x}^{\alpha}\theta_{r}\|_{L_{t}^{\infty}L_{\vec x}^{2}}^{2})+C\epsilon^{2}\int_{0}^{t}(\|\theta_{r}\|_{L_{\vec x}^{2}}^{2}+\|D_{x}^{\alpha}\theta_{r}\|_{L_{\vec x}^{2}}^{2})\mathrm{d}s,
\end{aligned}
\end{equation}
and
\begin{equation}\label{bd12}
\begin{aligned}
H_{3}^{5}+H_{3}^{6}=&-\int_{0}^{t}\int_{\mathbb{T}^{3}}\int_{\mathbb{S}^{2}}[D_{x}^{\alpha}, 4L_{0}]\theta_{r}(D_{x}^{\alpha}R_{2}+D_{x}^{\alpha}\langle R\rangle)\mathrm{d}\vec w\mathrm{d}\vec x\mathrm{d}s
\\\leq& \Big(\frac{1}{2\epsilon^{2}}\Big)\int_{0}^{t}\int_{\mathbb{T}^{3}}\int_{\mathbb{S}^{2}}|D_{x}^{\alpha}(R_{2}+\langle R\rangle)|^{2}\mathrm{d}\vec w\mathrm{d}\vec x\mathrm{d}s+C\epsilon^{2}\int_{0}^{t}(\|\theta_{r}\|_{L_{\vec x}^{2}}^{2}+\|D_{x}^{\alpha-1}\theta_{r}\|_{L_{\vec x}^{2}}^{2})\mathrm{d}s,
\end{aligned}
\end{equation}
which further implies
\begin{equation}\label{bd13}
\begin{aligned}
H_{3}\leq& C\epsilon^{2}\int_{0}^{t}\int_{\mathbb{T}^{3}}\int_{\mathbb{S}^{2}}(|D_{x}^{\alpha}\theta_{r}|^{2}
+|\theta_{r}|^{2})\mathrm{d}\vec w\mathrm{d}\vec x\mathrm{d}s+\delta\epsilon^{2}\int_{0}^{t}\int_{\mathbb{T}^{3}}\int_{\mathbb{S}^{2}}|D_{x}^{\alpha+1}\theta_{r}|^{2}
\mathrm{d}\vec w\mathrm{d}\vec x\mathrm{d}s\\&+C\int_{0}^{t}\int_{\mathbb{T}^{3}}\int_{\mathbb{S}^{2}}(|D_{x}^{\alpha-1}\theta_{r}|^{2}
+|D_{x}^{\alpha-2}\theta_{r}|^{2})
\mathrm{d}\vec w\mathrm{d}\vec x\mathrm{d}s
+\delta\int_{0}^{t}\int_{\mathbb{T}^{3}}\int_{\mathbb{S}^{2}}
|D_{x}^{\alpha}(4L_{0}\theta_{r})-D_{x}^{\alpha}f_{r}|^{2}
\mathrm{d}\vec w\mathrm{d}\vec x\mathrm{d}s
\\&+C\eta\epsilon^{2}(\|\theta_{r}\|_{L_{t}^{\infty}L_{\vec x}^{2}}^{2}+\|D_{x}^{\alpha}\theta_{r}\|_{L_{t}^{\infty}L_{\vec x}^{2}}^{2}).
\end{aligned}
\end{equation}
For $H_{6}-H_{9}$, we have
\begin{equation}\label{bd14}
\begin{aligned}
H_{6}=&\int_{0}^{t}\int_{\mathbb{T}^{3}}\int_{\mathbb{S}^{2}}D_{x}^{\alpha}(4L_{1}\theta_{r})D_{x}^{\alpha}f_{r}\mathrm{d}\vec w\mathrm{d}\vec x\mathrm{d}s
\\\leq& C(\|\theta_{r}\|_{L_{t}^{\infty}L_{\vec x}^{2}}+\|D_{x}^{\alpha}\theta_{r}\|_{L_{t}^{\infty}L_{\vec x}^{2}})\|D_{x}^{\alpha}f_{r}\|_{L_{t}^{\infty}L_{\vec w}^{2}L_{\vec x}^{2}}\int_{0}^{t}\|L_{1}\|_{H_{\vec x}^{2}}\mathrm{d}s
\\\leq& C(\|\theta_{r}\|_{L_{t}^{\infty}L_{\vec x}^{2}}+\|D_{x}^{\alpha}\theta_{r}\|_{L_{t}^{\infty}L_{\vec x}^{2}})\|D_{x}^{\alpha}f_{r}\|_{L_{t}^{\infty}L_{\vec w}^{2}L_{\vec x}^{2}}\int_{0}^{\frac{t}{\epsilon^{2}}}\|L_{1}\|_{H_{\vec x}^{2}}\epsilon^{2}\mathrm{d}\tau
\\\leq& C\eta\epsilon^{2}(\|\theta_{r}\|_{L_{t}^{\infty}L_{\vec x}^{2}}^{2}+\|D_{x}^{\alpha}\theta_{r}\|_{L_{t}^{\infty}L_{\vec x}^{2}}^{2}+\|D_{x}^{\alpha}f_{r}\|_{L_{t}^{\infty}L_{\vec w}^{2}L_{\vec x}^{2}}^{2})
,
\end{aligned}
\end{equation}

\begin{equation}
\begin{aligned}
H_{7}=&\int_{0}^{t}\int_{\mathbb{T}^{3}}\int_{\mathbb{S}^{2}}D_{x}^{\alpha}(4L_{2}\theta_{r})D_{x}^{\alpha}f_{r}\mathrm{d}\vec w\mathrm{d}\vec x\mathrm{d}s
\\\leq& C\epsilon^{2}\int_{0}^{t}(\|\theta_{r}\|_{L_{\vec x}^{2}}^{2}+\|D_{x}^{\alpha}\theta_{r}\|_{L_{\vec x}^{2}}^{2})\mathrm{d}s+C\epsilon^{2}\int_{0}^{t}\|D_{x}^{\alpha}f_{r}\|_{L_{\vec w}^{2}L_{\vec x}^{2}}^{2}\mathrm{d}s,
\end{aligned}
\end{equation}

\begin{equation}\label{bd15}
\begin{aligned}
H_{8}=&-\int_{0}^{t}\int_{\mathbb{T}^{3}}\int_{\mathbb{S}^{2}}D_{x}^{\alpha}(4L_{1}\theta_{r})D_{x}^{\alpha}(4L_{0}\theta_{r})\mathrm{d}\vec w\mathrm{d}\vec x\mathrm{d}s
\\\leq& C(\|\theta_{r}\|_{L_{t}^{\infty}L_{\vec x}^{2}}+\|D_{x}^{\alpha}\theta_{r}\|_{L_{t}^{\infty}L_{\vec x}^{2}})^{2}\int_{0}^{t}\|L_{1}\|_{H_{\vec x}^{2}}\mathrm{d}s
\\\leq& C\epsilon^{2}(\|\theta_{r}\|_{L_{t}^{\infty}L_{\vec x}^{2}}+\|D_{x}^{\alpha}\theta_{r}\|_{L_{t}^{\infty}L_{\vec x}^{2}})^{2}\int_{0}^{\frac{t}{\epsilon^{2}}}\|L_{1}\|_{H_{\vec x}^{2}}\mathrm{d}\tau
\\\leq& C\eta\epsilon^{2}(\|\theta_{r}\|_{L_{t}^{\infty}L_{\vec x}^{2}}^{2}+\|D_{x}^{\alpha}\theta_{r}\|_{L_{t}^{\infty}L_{\vec x}^{2}}^{2}),
\end{aligned}
\end{equation}
and
\begin{equation}
\begin{aligned}
H_{9}=&-\int_{0}^{t}\int_{\mathbb{T}^{3}}\int_{\mathbb{S}^{2}}D_{x}^{\alpha}(4L_{2}\theta_{r})D_{x}^{\alpha}(4L_{0}\theta_{r})\mathrm{d}\vec w\mathrm{d}\vec x\mathrm{d}s
\\\leq& C\epsilon^{2}\int_{0}^{t}(\|\theta_{r}\|_{L_{\vec x}^{2}}^{2}+\|D_{x}^{\alpha}\theta_{r}\|_{L_{\vec x}^{2}}^{2})\mathrm{d}s.
\end{aligned}
\end{equation}

For sufficiently small $\epsilon$, we have $4L_{0}\geq\kappa>0$, where $\kappa$ is a given positive constant. Then, collecting \eqref{bd0} and the estimates of $H_{1}-H_{11}$, we have
\begin{equation}\label{bd0-jia1}
\begin{aligned}
&\epsilon^{2}\int_{\mathbb{T}^{3}}\int_{\mathbb{S}^{2}}|D_{x}^{\alpha}f_{r}|^{2}(t)\mathrm{d}\vec w\mathrm{d}\vec x
+\kappa\epsilon^{2}\int_{\mathbb{T}^{3}}\int_{\mathbb{S}^{2}}|D_{x}^{\alpha}\theta_{r}|^{2}(t)
\mathrm{d}\vec w\mathrm{d}\vec x
\\&+\epsilon^{2}\int_{0}^{t}\int_{\mathbb{T}^{3}}\int_{\mathbb{S}^{2}}|D_{x}^{\alpha}f_{r}-D_{x}^{\alpha}\overline{f_{r}}|^{2}
\mathrm{d}\vec w\mathrm{d}\vec x\mathrm{d}s
+\int_{0}^{t}\int_{\mathbb{T}^{3}}\int_{\mathbb{S}^{2}}|D_{x}^{\alpha}(4L_{0}\theta_{r})
-D_{x}^{\alpha}f_{r}|^{2}\mathrm{d}\vec w\mathrm{d}\vec x\mathrm{d}s\\&+\kappa\epsilon^{2}\int_{0}^{t}\int_{\mathbb{T}^{3}}\int_{\mathbb{S}^{2}}|\nabla_{x}D_{x}^{\alpha}\theta_{r}
|^{2}\mathrm{d}\vec w\mathrm{d}\vec x\mathrm{d}s
\\\leq& C\epsilon^{2}\int_{0}^{t}\int_{\mathbb{T}^{3}}\int_{\mathbb{S}^{2}}(|D_{x}^{\alpha}\theta_{r}|^{2}
+|\theta_{r}|^{2})\mathrm{d}\vec w\mathrm{d}\vec x\mathrm{d}s+C\delta\epsilon^{2}\int_{0}^{t}\int_{\mathbb{T}^{3}}\int_{\mathbb{S}^{2}}|D_{x}^{\alpha+1}\theta_{r}|^{2}
\mathrm{d}\vec w\mathrm{d}\vec x\mathrm{d}s\\&+C\int_{0}^{t}\int_{\mathbb{T}^{3}}\int_{\mathbb{S}^{2}}(|D_{x}^{\alpha-1}\theta_{r}|^{2}
+|D_{x}^{\alpha-2}\theta_{r}|^{2})
\mathrm{d}\vec w\mathrm{d}\vec x\mathrm{d}s
+\delta\int_{0}^{t}\int_{\mathbb{T}^{3}}\int_{\mathbb{S}^{2}}
|D_{x}^{\alpha}(4L_{0}\theta_{r})-D_{x}^{\alpha}f_{r}|^{2}
\mathrm{d}\vec w\mathrm{d}\vec x\mathrm{d}s
\\&+C\eta\epsilon^{2}(\|\theta_{r}\|_{L_{t}^{\infty}L_{\vec x}^{2}}^{2}+\|D_{x}^{\alpha}\theta_{r}\|_{L_{t}^{\infty}L_{\vec x}^{2}}^{2}+\|D_{x}^{\alpha}f_{r}\|_{L_{t}^{\infty}L_{\vec w}^{2}L_{\vec x}^{2}}^{2})+C\epsilon^{2}\int_{0}^{t}\|D_{x}^{\alpha}f_{r}\|_{L_{\vec w}^{2}L_{\vec x}^{2}}^{2}\mathrm{d}s
\\&+\frac{C}{\epsilon^{2}}\int_{0}^{t}\int_{\mathbb{T}^{3}}\int_{\mathbb{S}^{2}}(|D_{x}^{\alpha}R_{1}|^{2}
+|D_{x}^{\alpha}R_{2}|^{2}
+|D_{x}^{\alpha}R|^{2})\mathrm{d}\vec w\mathrm{d}\vec x\mathrm{d}s.
\end{aligned}
\end{equation}
Taking the sup from 0 to $t$ on the left, choosing $\delta=\min\{\frac{\kappa}{2C}, \frac{1}{2}\}$ and $C\eta=\min\{\frac{\kappa}{4}, \frac{1}{2}\}$, we have
\begin{equation}\nonumber
\begin{aligned}
&\frac{1}{2}\epsilon^{2}\|D_{x}^{\alpha}f_{r}\|_{L_{t}^{\infty}L_{\vec w}^{2}L_{\vec x}^{2}}^{2}
+\frac{3\kappa}{4}\epsilon^{2}\|D_{x}^{\alpha}\theta_{r}\|_{L_{t}^{\infty}L_{\vec x}^{2}}^{2}
+\epsilon^{2}\int_{0}^{t}\int_{\mathbb{T}^{3}}\int_{\mathbb{S}^{2}}|D_{x}^{\alpha}f_{r}-D_{x}^{\alpha}\overline{f_{r}}|^{2}
\mathrm{d}\vec w\mathrm{d}\vec x\mathrm{d}s
\\&+\frac{1}{2}\int_{0}^{t}\int_{\mathbb{T}^{3}}\int_{\mathbb{S}^{2}}|D_{x}^{\alpha}(4L_{0}\theta_{r})
-D_{x}^{\alpha}f_{r}|^{2}\mathrm{d}\vec w\mathrm{d}\vec x\mathrm{d}s+\frac{\kappa}{2}\epsilon^{2}\int_{0}^{t}\int_{\mathbb{T}^{3}}\int_{\mathbb{S}^{2}}|\nabla_{x}D_{x}^{\alpha}\theta_{r}
|^{2}\mathrm{d}\vec w\mathrm{d}\vec x\mathrm{d}s
\end{aligned}
\end{equation}
\begin{equation}\label{bd0-jia2}
\begin{aligned}
\\\leq& C\epsilon^{2}\int_{0}^{t}\int_{\mathbb{T}^{3}}\int_{\mathbb{S}^{2}}(|D_{x}^{\alpha}\theta_{r}|^{2}
+|\theta_{r}|^{2})\mathrm{d}\vec w\mathrm{d}\vec x\mathrm{d}s+C\int_{0}^{t}\int_{\mathbb{T}^{3}}\int_{\mathbb{S}^{2}}(|D_{x}^{\alpha-1}\theta_{r}|^{2}
+|D_{x}^{\alpha-2}\theta_{r}|^{2})
\mathrm{d}\vec w\mathrm{d}\vec x\mathrm{d}s
\\&+\frac{\kappa}{4}\epsilon^{2}\|\theta_{r}\|_{L_{t}^{\infty}L_{\vec x}^{2}}^{2}+C\epsilon^{2}\int_{0}^{t}\|D_{x}^{\alpha}f_{r}\|_{L_{\vec w}^{2}L_{\vec x}^{2}}^{2}\mathrm{d}t+\frac{C}{\epsilon^{2}}\int_{0}^{t}\int_{\mathbb{T}^{3}}\int_{\mathbb{S}^{2}}(|D_{x}^{\alpha}R_{1}|^{2}
+|D_{x}^{\alpha}R_{2}|^{2}
\\&+|D_{x}^{\alpha}R|^{2})\mathrm{d}\vec w\mathrm{d}\vec x\mathrm{d}s.
\end{aligned}
\end{equation}

In the following, we will get the estimates for $\alpha=0, 1, 2$, iteratively.

When $\alpha=0$, we can write \eqref{bd0-jia2} in the following form:
\begin{equation}\label{bd0-jia3}
\begin{aligned}
&\frac{1}{2}\epsilon^{2}\|f_{r}\|_{L_{t}^{\infty}L_{\vec w}^{2}L_{\vec x}^{2}}^{2}
+\frac{\kappa}{2}\epsilon^{2}\|\theta_{r}\|_{L_{t}^{\infty}L_{\vec x}^{2}}^{2}
+\epsilon^{2}\int_{0}^{t}\int_{\mathbb{T}^{3}}\int_{\mathbb{S}^{2}}|f_{r}-\overline{f_{r}}|^{2}
\mathrm{d}\vec w\mathrm{d}\vec x\mathrm{d}s
\\&+\frac{1}{2}\int_{0}^{t}\int_{\mathbb{T}^{3}}\int_{\mathbb{S}^{2}}|4L_{0}\theta_{r}
-f_{r}|^{2}\mathrm{d}\vec w\mathrm{d}\vec x\mathrm{d}s+\frac{\kappa}{2}\epsilon^{2}\int_{0}^{t}\int_{\mathbb{T}^{3}}\int_{\mathbb{S}^{2}}|\nabla_{x}\theta_{r}
|^{2}\mathrm{d}\vec w\mathrm{d}\vec x\mathrm{d}s
\\\leq& C\epsilon^{2}\int_{0}^{t}\int_{\mathbb{T}^{3}}\int_{\mathbb{S}^{2}}
|\theta_{r}|^{2}\mathrm{d}\vec w\mathrm{d}\vec x\mathrm{d}s
+C\epsilon^{2}\int_{0}^{t}\|f_{r}\|_{L_{\vec w}^{2}L_{\vec x}^{2}}^{2}\mathrm{d}s\\&+\frac{C}{\epsilon^{2}}\int_{0}^{t}\int_{\mathbb{T}^{3}}\int_{\mathbb{S}^{2}}(|R_{1}|^{2}
+|R_{2}|^{2}
+|R|^{2})\mathrm{d}\vec w\mathrm{d}\vec x\mathrm{d}s,
\end{aligned}
\end{equation}
which further implies
\begin{equation}\label{bd0-jia4}
\begin{aligned}
&\epsilon^{2}\|f_{r}(t)\|_{L_{\vec w}^{2}L_{\vec x}^{2}}^{2}
+\epsilon^{2}\|\theta_{r}(t)\|_{L_{\vec x}^{2}}^{2}
+\epsilon^{2}\int_{0}^{t}\int_{\mathbb{T}^{3}}\int_{\mathbb{S}^{2}}|f_{r}-\overline{f_{r}}|^{2}
\mathrm{d}\vec w\mathrm{d}\vec x\mathrm{d}s
\\&+\int_{0}^{t}\int_{\mathbb{T}^{3}}\int_{\mathbb{S}^{2}}|4L_{0}\theta_{r}
-f_{r}|^{2}\mathrm{d}\vec w\mathrm{d}\vec x\mathrm{d}t+\epsilon^{2}\int_{0}^{t}\int_{\mathbb{T}^{3}}\int_{\mathbb{S}^{2}}|\nabla_{x}\theta_{r}
|^{2}\mathrm{d}\vec w\mathrm{d}\vec x\mathrm{d}s
\\\leq& C\epsilon^{2}\int_{0}^{t}\int_{\mathbb{T}^{3}}\int_{\mathbb{S}^{2}}
|\theta_{r}|^{2}\mathrm{d}\vec w\mathrm{d}\vec x\mathrm{d}s
+C\epsilon^{2}\int_{0}^{t}\|f_{r}\|_{L_{\vec w}^{2}L_{\vec x}^{2}}^{2}\mathrm{d}s\\&+\frac{C}{\epsilon^{2}}\int_{0}^{t}\int_{\mathbb{T}^{3}}\int_{\mathbb{S}^{2}}(|R_{1}|^{2}
+|R_{2}|^{2}
+|R|^{2})\mathrm{d}\vec w\mathrm{d}\vec x\mathrm{d}s.
\end{aligned}
\end{equation}
By the Gronwall's inequality, we have
\begin{equation}\label{bd0-jia5}
\begin{aligned}
&\epsilon^{2}\|f_{r}\|_{L_{t}^{\infty}L_{\vec w}^{2}L_{\vec x}^{2}}^{2}
+\epsilon^{2}\|\theta_{r}\|_{L_{t}^{\infty}L_{\vec x}^{2}}^{2}
+\epsilon^{2}\int_{0}^{t}\int_{\mathbb{T}^{3}}\int_{\mathbb{S}^{2}}|f_{r}-\overline{f_{r}}|^{2}
\mathrm{d}\vec w\mathrm{d}\vec x\mathrm{d}s
\\&+\int_{0}^{t}\int_{\mathbb{T}^{3}}\int_{\mathbb{S}^{2}}|4L_{0}\theta_{r}
-f_{r}|^{2}\mathrm{d}\vec w\mathrm{d}\vec x\mathrm{d}s+\epsilon^{2}\int_{0}^{t}\int_{\mathbb{T}^{3}}\int_{\mathbb{S}^{2}}|\nabla_{x}\theta_{r}
|^{2}\mathrm{d}\vec w\mathrm{d}\vec x\mathrm{d}s\\&+C\epsilon^{2}\int_{0}^{t}\int_{\mathbb{T}^{3}}\int_{\mathbb{S}^{2}}
|\theta_{r}|^{2}\mathrm{d}\vec w\mathrm{d}\vec x\mathrm{d}s
+C\epsilon^{2}\int_{0}^{t}\|f_{r}\|_{L_{\vec w}^{2}L_{\vec x}^{2}}^{2}\mathrm{d}s
\\\leq& \frac{C(t)}{\epsilon^{2}}\int_{0}^{t}\int_{\mathbb{T}^{3}}\int_{\mathbb{S}^{2}}(|R_{1}|^{2}
+|R_{2}|^{2}
+|R|^{2})\mathrm{d}\vec w\mathrm{d}\vec x\mathrm{d}s.
\end{aligned}
\end{equation}

When $\alpha=1$, we can write \eqref{bd0-jia2} in the following form:
\begin{equation}\label{bd0-jia6}
\begin{aligned}
&\frac{1}{2}\epsilon^{2}\|D_{x}f_{r}\|_{L_{t}^{\infty}L_{\vec w}^{2}L_{\vec x}^{2}}^{2}
+\frac{3\kappa}{4}\epsilon^{2}\|D_{x}\theta_{r}\|_{L_{t}^{\infty}L_{\vec x}^{2}}^{2}
+\epsilon^{2}\int_{0}^{t}\int_{\mathbb{T}^{3}}\int_{\mathbb{S}^{2}}|D_{x}f_{r}-D_{x}\overline{f_{r}}|^{2}
\mathrm{d}\vec w\mathrm{d}\vec x\mathrm{d}s
\\&+\frac{1}{2}\int_{0}^{t}\int_{\mathbb{T}^{3}}\int_{\mathbb{S}^{2}}|D_{x}(4L_{0}\theta_{r})
-D_{x}f_{r}|^{2}\mathrm{d}\vec w\mathrm{d}\vec x\mathrm{d}s+\frac{\kappa}{2}\epsilon^{2}\int_{0}^{t}\int_{\mathbb{T}^{3}}\int_{\mathbb{S}^{2}}|\nabla_{x}D_{x}\theta_{r}
|^{2}\mathrm{d}\vec w\mathrm{d}\vec x\mathrm{d}s
\\\leq& C\epsilon^{2}\int_{0}^{t}\int_{\mathbb{T}^{3}}\int_{\mathbb{S}^{2}}(|D_{x}\theta_{r}|^{2}
+|\theta_{r}|^{2})\mathrm{d}\vec w\mathrm{d}\vec x\mathrm{d}s+C\int_{0}^{t}\int_{\mathbb{T}^{3}}\int_{\mathbb{S}^{2}}|\theta_{r}|^{2}
\mathrm{d}\vec w\mathrm{d}\vec x\mathrm{d}s
\\&+\frac{\kappa}{4}\epsilon^{2}\|\theta_{r}\|_{L_{t}^{\infty}L_{\vec x}^{2}}^{2}+C\epsilon^{2}\int_{0}^{t}\|D_{x}f_{r}\|_{L_{\vec w}^{2}L_{\vec x}^{2}}^{2}\mathrm{d}s+\frac{C}{\epsilon^{2}}\int_{0}^{t}\int_{\mathbb{T}^{3}}\int_{\mathbb{S}^{2}}(|D_{x}R_{1}|^{2}
+|D_{x}R_{2}|^{2}
\\&+|D_{x}R|^{2})\mathrm{d}\vec w\mathrm{d}\vec x\mathrm{d}s
\\\leq& C\epsilon^{2}\int_{0}^{t}\int_{\mathbb{T}^{3}}\int_{\mathbb{S}^{2}}|D_{x}\theta_{r}|^{2}
\mathrm{d}\vec w\mathrm{d}\vec x\mathrm{d}s+C\epsilon^{2}\int_{0}^{t}\|D_{x}f_{r}\|_{L_{\vec w}^{2}L_{\vec x}^{2}}^{2}\mathrm{d}s+\frac{C}{\epsilon^{2}}\int_{0}^{t}\int_{\mathbb{T}^{3}}\int_{\mathbb{S}^{2}}(|D_{x}R_{1}|^{2}
+|D_{x}R_{2}|^{2}
\\&+|D_{x}R|^{2})\mathrm{d}\vec w\mathrm{d}\vec x\mathrm{d}s+\frac{C(t)}{\epsilon^{4}}\int_{0}^{t}\int_{\mathbb{T}^{3}}\int_{\mathbb{S}^{2}}(|R_{1}|^{2}
+|R_{2}|^{2}
+|R|^{2})\mathrm{d}\vec w\mathrm{d}\vec x\mathrm{d}s.
\end{aligned}
\end{equation}
By the Gronwall's inequality, we have
\begin{equation}\label{bd0-jia7}
\begin{aligned}
&\epsilon^{2}\|D_{x}f_{r}\|_{L_{t}^{\infty}L_{\vec w}^{2}L_{\vec x}^{2}}^{2}
+\epsilon^{2}\|D_{x}\theta_{r}\|_{L_{t}^{\infty}L_{\vec x}^{2}}^{2}
+\epsilon^{2}\int_{0}^{t}\int_{\mathbb{T}^{3}}\int_{\mathbb{S}^{2}}|D_{x}f_{r}-D_{x}\overline{f_{r}}|^{2}
\mathrm{d}\vec w\mathrm{d}\vec x\mathrm{d}s
\\&+\int_{0}^{t}\int_{\mathbb{T}^{3}}\int_{\mathbb{S}^{2}}|D_{x}(4L_{0}\theta_{r})
-D_{x}f_{r}|^{2}\mathrm{d}\vec w\mathrm{d}\vec x\mathrm{d}s+\epsilon^{2}\int_{0}^{t}\int_{\mathbb{T}^{3}}\int_{\mathbb{S}^{2}}|\nabla_{x}D_{x}\theta_{r}
|^{2}\mathrm{d}\vec w\mathrm{d}\vec x\mathrm{d}s\\&+C\epsilon^{2}\int_{0}^{t}\int_{\mathbb{T}^{3}}\int_{\mathbb{S}^{2}}
|D_{x}\theta_{r}|^{2}\mathrm{d}\vec w\mathrm{d}\vec x\mathrm{d}s
+C\epsilon^{2}\int_{0}^{t}\|D_{x}f_{r}\|_{L_{\vec w}^{2}L_{\vec x}^{2}}^{2}\mathrm{d}s
\\\leq& \frac{C(t)}{\epsilon^{2}}\int_{0}^{t}\int_{\mathbb{T}^{3}}\int_{\mathbb{S}^{2}}(|D_{x}R_{1}|^{2}
+|D_{x}R_{2}|^{2}
+|D_{x}R|^{2})\mathrm{d}\vec w\mathrm{d}\vec x\mathrm{d}s\\&+\frac{C(t)}{\epsilon^{4}}\int_{0}^{t}\int_{\mathbb{T}^{3}}\int_{\mathbb{S}^{2}}(|R_{1}|^{2}
+|R_{2}|^{2}
+|R|^{2})\mathrm{d}\vec w\mathrm{d}\vec x\mathrm{d}s.
\end{aligned}
\end{equation}
For $\alpha=2$, we can finally get
\begin{equation}\label{bd0-jia8}
\begin{aligned}
&\epsilon^{2}\|D_{x}^{2}f_{r}\|_{L_{t}^{\infty}L_{\vec w}^{2}L_{\vec x}^{2}}^{2}
+\epsilon^{2}\|D_{x}^{2}\theta_{r}\|_{L_{t}^{\infty}L_{\vec x}^{2}}^{2}
+\epsilon^{2}\int_{0}^{t}\int_{\mathbb{T}^{3}}\int_{\mathbb{S}^{2}}|D_{x}^{2}f_{r}-D_{x}^{2}\overline{f_{r}}|^{2}
\mathrm{d}\vec w\mathrm{d}\vec x\mathrm{d}s
\\&+\int_{0}^{t}\int_{\mathbb{T}^{3}}\int_{\mathbb{S}^{2}}|D_{x}^{2}(4L_{0}\theta_{r})
-D_{x}^{2}f_{r}|^{2}\mathrm{d}\vec w\mathrm{d}\vec x\mathrm{d}s+\epsilon^{2}\int_{0}^{t}\int_{\mathbb{T}^{3}}\int_{\mathbb{S}^{2}}|\nabla_{x}D_{x}^{2}\theta_{r}
|^{2}\mathrm{d}\vec w\mathrm{d}\vec x\mathrm{d}s\\&+C\epsilon^{2}\int_{0}^{t}\int_{\mathbb{T}^{3}}\int_{\mathbb{S}^{2}}
|D_{x}^{2}\theta_{r}|^{2}\mathrm{d}\vec w\mathrm{d}\vec x\mathrm{d}s
+C\epsilon^{2}\int_{0}^{t}\|D_{x}^{2}f_{r}\|_{L_{\vec w}^{2}L_{\vec x}^{2}}^{2}\mathrm{d}s
\\\leq& \frac{C(t)}{\epsilon^{2}}\int_{0}^{t}\int_{\mathbb{T}^{3}}\int_{\mathbb{S}^{2}}(|D_{x}^{2}R_{1}|^{2}
+|D_{x}^{2}R_{2}|^{2}
+|D_{x}^{2}R|^{2})\mathrm{d}\vec w\mathrm{d}\vec x\mathrm{d}s\\&+\frac{C(t)}{\epsilon^{4}}\int_{0}^{t}\int_{\mathbb{T}^{3}}\int_{\mathbb{S}^{2}}(|D_{x}R_{1}|^{2}
+|D_{x}R_{2}|^{2}
+|D_{x}R|^{2})\mathrm{d}\vec w\mathrm{d}\vec x\mathrm{d}s\\&+\frac{C(t)}{\epsilon^{6}}\int_{0}^{t}\int_{\mathbb{T}^{3}}\int_{\mathbb{S}^{2}}(|R_{1}|^{2}
+|R_{2}|^{2}
+|R|^{2})\mathrm{d}\vec w\mathrm{d}\vec x\mathrm{d}s,
\end{aligned}
\end{equation}
as the procedure of $\alpha=0$ and $\alpha=1$.

Collecting \eqref{bd0-jia3}, \eqref{bd0-jia7} and \eqref{bd0-jia8}, we can derive \eqref{bd00}.

\subsection{Nonlinear system}

We now show the existence and uniqueness of solutions to system \eqref{research equations} around the constructed composite approximate solution $(f^{N}, \theta^{N})$ and finish the proof of Theorem \ref{result}. Due to the equivalence between system \eqref{reeq1}-\eqref{rein} and \eqref{research equations}, we only need to show the existence and uniqueness of solutions for \eqref{reeq1}-\eqref{rein} in the neighborhood of zero.

\noindent\textbf{Proof of Theorem \ref{result}.} The proof of existence and uniqueness is obtained using the Banach fixed point theorem. We first construct a sequence of functions and then show the sequence is a contraction sequence. Finally we show the convergence of $(f_{r}, \theta_{r})$ to zero as $\epsilon\rightarrow 0$.

\noindent\textbf{Step 1.}Construction of sequence of functions. Let $\{f_{r}^{0}, \theta_{r}^{0}\}$ be zero functions
\begin{equation}
f_{r}^{0}=0, \ \ \theta_{r}^{0}=0,
\end{equation}
and for $k\geq 1,$ $\{f_{r}^{k+1}, \theta_{r}^{k+1}\}$ are defined recursively by
\begin{equation}\label{reeq1-jia}
\begin{aligned}
&\epsilon^{2}\partial_{t}f_{r}^{k+1}+\epsilon\vec{w}\cdot\nabla_{x}f_{r}^{k+1}+\epsilon^{2}(f_{r}^{k+1}
-\overline{f_{r}^{k+1}})+f_{r}^{k+1}
-4(\theta^{N})^{3}\theta_{r}^{k+1}\\=&-\mathcal{L}_{1}(f^{N}, \theta^{N})+6(\theta^{N})^{2}(\theta_{r}^{k})^{2}+4\theta^{N}(\theta_{r}^{k})^{3}+(\theta_{r}^{k})^{4},
\end{aligned}
\end{equation}
\begin{equation}\label{reeq2-jia}
\begin{aligned}
&\epsilon^{2}\partial_{t}\theta_{r}^{k+1}+\epsilon^{2}\mathrm{div}{(\vec u\theta_{r}^{k+1})}-\epsilon^{2}\Delta_{x}\theta_{r}^{k+1}
+4(\theta^{N})^{3}\theta_{r}^{k+1}-\overline{f_{r}^{k+1}}\\=&-\mathcal{L}_{2}(f^{N}, \theta^{N})-6(\theta^{N})^{2}(\theta_{r}^{k})^{2}-4\theta^{N}(\theta_{r}^{k})^{3}-(\theta_{r}^{k})^{4},
\end{aligned}
\end{equation}
with initial conditions
\begin{equation}\label{rein-jia}
f_{r}^{k+1}(0, \vec x, \vec w)=0,\ \ \theta_{r}^{k+1}(0, \vec x)=0,\ \ \mathrm{for}\ \ (\vec x, \vec w)\in\mathbb{T}^{3}\times\mathbb{S}^{2}.
\end{equation}
The above system defines a mapping $\mathcal{T}$ with $(f_{r}^{k+1}, \theta_{r}^{k+1})=\mathcal{T}((f_{r}^{k}, \theta_{r}^{k}))$.

\noindent\textbf{Step 2.} The contraction mapping. We consider the solution in the function space
\begin{equation}
O_{q}:=\{(f_{r}, \theta_{r})\in L^{\infty}_{T}L^{2}_{\vec w}H^{2}_{\vec x}\times L^{\infty}_{T}H^{2}_{\vec x}: \|f_{r}\|_{L^{\infty}_{T}L^{2}_{\vec w}H^{2}_{\vec x}}+\|\theta_{r}\|_{L^{\infty}_{T}H^{2}_{\vec x}}\leq \epsilon^{q}\},
\end{equation}
where $q>0$ is a constant to be chosen later.

First, we show $\mathcal{T}$ maps the space $O_{q}$ into itself. Assume the residuals satisfy 
\begin{equation}
\|\mathcal{L}_{1}(f^{N}, \theta^{N})\|_{L^{\infty}_{T}L^{2}_{\vec w}H^{2}_{\vec x}}, \|\mathcal{L}_{2}(f^{N}, \theta^{N})\|_{L^{\infty}_{T}L^{2}_{\vec w}H^{2}_{\vec x}}=O(\epsilon^{N+1}). 
\end{equation}
By \eqref{bd00} with $R_{1}=-\mathcal{L}_{1}(f^{N}, \theta^{N})$, $R_{2}=-\mathcal{L}_{2}(f^{N}, \theta^{N})$ and
$R=-6(\theta^{N})^{2}(\theta_{r}^{k})^{2}-4\theta^{N}(\theta_{r}^{k})^{3}-(\theta_{r}^{k})^{4}$, the following estimate holds (assuming $q\geq 1$):
\begin{equation}\label{bd00-jia1}
\begin{aligned}
&\|f_{r}^{k+1}\|_{L_{T}^{\infty}H^{2}_{\vec x}L_{\vec w}^{2}}+\|\theta_{r}^{k+1}\|_{L_{T}^{\infty}H_{\vec x}^{2}}+\frac{1}{\epsilon^{2}}\|f_{r}^{k+1}-4(\theta^{N})^{3}\theta_{r}^{k+1}\|_{L_{T}^{2}H_{\vec x}^{2}L_{\vec w}^{2}}+\|\theta_{r}^{k+1}\|_{L_{T}^{2}H_{\vec x}^{3}}
\\\leq& \frac{C(T)}{\epsilon^{4}}\Big(\|6(\theta^{N})^{2}(\theta_{r}^{k})^{2}+4\theta^{N}(\theta_{r}^{k})^{3}
+(\theta_{r}^{k})^{4}\|_{L_{T}^{2}H_{\vec x}^{2}L_{\vec w}^{2}}+\|\mathcal{L}_{1}(f^{N}, \theta^{N})\|_{L_{T}^{2}H_{\vec x}^{2}L_{\vec w}^{2}}\\&+\|\mathcal{L}_{2}(f^{N}, \theta^{N})\|_{L_{T}^{2}H_{\vec x}^{2}}\Big)
\\\leq& \frac{C(T)}{\epsilon^{4}}(\|\theta_{r}^{k}\|_{L_{T}^{\infty}H^{2}_{\vec x}}^{2}+\|\theta_{r}^{k}\|_{L_{T}^{\infty}H^{2}_{\vec x}}^{4}+\epsilon^{N+1})
\\\leq& \frac{C(T)}{\epsilon^{4}}(\epsilon^{2q}+\epsilon^{N+1})
\\\leq& C(T)(\epsilon^{2q-4}+\epsilon^{N-3}).
\end{aligned}
\end{equation}
Assuming $2q-4>q$ and $N-3>q$, i.e. $q>4$ and $N>3+q$, the above inequality implies
\begin{equation}
\|f_{r}^{k+1}\|_{L_{T}^{\infty}L_{\vec w}^{2}H^{2}_{\vec x}}+\|\theta_{r}^{k+1}\|_{L_{T}^{\infty}H_{\vec x}^{2}}\leq C(T)(\epsilon^{2m-4}+\epsilon^{N-3})\leq \epsilon^{q},
\end{equation}
for sufficiently small $\epsilon$. Thus we obtain that $(f_{r}^{k+1}, \theta_{r}^{k+1})\in O_{q}$ and therefore $\mathcal{T}$ maps $O_{q}$ into itself.

Next, we show the map $\mathcal{T}$ is a contraction mapping. Let $\varphi^{k+1}=f_{r}^{k+1}-f_{r}^{k}$, $h^{k+1}=\theta_{r}^{k+1}-\theta_{r}^{k}$, then they satisfy
\begin{equation}\label{reeq1-jia-new}
\begin{aligned}
&\epsilon^{2}\partial_{t}\varphi^{k+1}+\epsilon\vec{w}\cdot\nabla_{x}\varphi^{k+1}+\epsilon^{2}(\varphi^{k+1}
-\overline{\varphi^{k+1}})+\varphi^{k+1}
-4(\theta^{N})^{3}h^{k+1}\\=&6(\theta^{N})^{2}(\theta_{r}^{k}+\theta_{r}^{k-1})h^{k}
+4\theta^{N}((\theta_{r}^{k})^{2}+\theta_{r}^{k}\theta_{r}^{k-1}+(\theta_{r}^{k-1})^{2})h^{k}
+((\theta_{r}^{k})^{2}+(\theta_{r}^{k-1})^{2})(\theta_{r}^{k}+\theta_{r}^{k-1})h^{k},
\end{aligned}
\end{equation}
\begin{equation}\label{reeq2-jia-new}
\begin{aligned}
&\epsilon^{2}\partial_{t}h^{k+1}+\epsilon^{2}\mathrm{div}{(\vec uh^{k+1})}-\epsilon^{2}\Delta_{x}h^{k+1}
+4(\theta^{N})^{3}h^{k+1}-\overline{f_{r}^{k+1}}\\=&-6(\theta^{N})^{2}(\theta_{r}^{k}+\theta_{r}^{k-1})h^{k}
-4\theta^{N}((\theta_{r}^{k})^{2}+\theta_{r}^{k}\theta_{r}^{k-1}+(\theta_{r}^{k-1})^{2})h^{k}
-((\theta_{r}^{k})^{2}+(\theta_{r}^{k-1})^{2})(\theta_{r}^{k}+\theta_{r}^{k-1})h^{k},
\end{aligned}
\end{equation}
with initial conditions
\begin{equation}\label{rein-jia-new}
\varphi^{k+1}(0, \vec x, \vec w)=0,\ \ h^{k+1}(0, \vec x)=0,\ \ \mathrm{for}\ \ (\vec x, \vec w)\in\mathbb{T}^{3}\times\mathbb{S}^{2}.
\end{equation}
Using \eqref{bd00} with $R_{1}=R_{2}=0$, $R=-6(\theta^{N})^{2}(\theta_{r}^{k}+\theta_{r}^{k-1})h^{k}
-4\theta^{N}((\theta_{r}^{k})^{2}+\theta_{r}^{k}\theta_{r}^{k-1}+(\theta_{r}^{k-1})^{2})h^{k}
-((\theta_{r}^{k})^{2}+(\theta_{r}^{k-1})^{2})(\theta_{r}^{k}+\theta_{r}^{k-1})h^{k}$, we obtain
\begin{equation}
\begin{aligned}
\|\varphi^{k+1}\|_{L_{T}^{\infty}L_{\vec w}^{2}H^{2}_{\vec x}}+\|h^{k+1}\|_{L_{T}^{\infty}H^{2}_{\vec x}}\leq C(T)\epsilon^{q-4}\|h^{k}\|_{L_{T}^{\infty}H^{2}_{\vec x}}\leq C(T)\epsilon^{q-4}(\|\varphi^{k}\|_{L_{T}^{\infty}L_{\vec w}^{2}H^{2}_{\vec x}}+\|h^{k}\|_{L_{T}^{\infty}H^{2}_{\vec x}}).
\end{aligned}
\end{equation}
Assume $q>4$, then $C(T)\epsilon^{q-4}<1$ for $\epsilon$ sufficiently small, the above inequality implies
\begin{equation}
\begin{aligned}
\|\varphi^{k+1}\|_{L_{T}^{\infty}L_{\vec w}^{2}H^{2}_{\vec x}}+\|h^{k+1}\|_{L_{T}^{\infty}H^{2}_{\vec x}}\leq C_{1}(\|\varphi^{k}\|_{L_{T}^{\infty}L_{\vec w}^{2}H^{2}_{\vec x}}+\|h^{k}\|_{L_{T}^{\infty}H^{2}_{\vec x}}),
\end{aligned}
\end{equation}
for some constant $0<C_{1}<1$. Therefore, for $q>4$ and $N>3+q$, $\mathcal{T}$ is a contraction mapping. By the Banach fixed point theorem, there exist a unique fixed point $(f_{r}, \theta_{r})$ such that $(f_{r}, \theta_{r})=\mathcal{T}((f_{r}, \theta_{r}))$. Therefore, there exists a unique solution to \eqref{reeq1}-\eqref{rein} in $O_{q}$.

Taking $q=5$ and $N=9$, we can conclude that
\begin{equation}
\|f_{r}\|_{L_{T}^{\infty}L_{\vec w}^{2}H^{2}_{\vec x}}+\|\theta_{r}\|_{L_{T}^{\infty}H^{2}_{\vec x}}\leq C\epsilon^{5}.
\end{equation}
In the following, we consider 
\begin{equation}\label{remainder equation-2}\left\{
\begin{split}
&\epsilon^{2}\partial_{t}f_{r}+\epsilon\vec{w}\cdot\nabla_{x}f_{r}+\epsilon^{2}(f_{r}-\overline{f_{r}})+f_{r}
=(\theta^{N}+\theta_{r})^{4}-(\theta^{N})^{4}-\mathcal{L}_{1}(f^{N}, \theta^{N})\ \ \mathrm{in} \ \ (0, T]\times \mathbb{T}^{3}\times \mathbb{S}^{2},\\
&f_{r}(0, \vec{x}, \vec{w})=0\ \ \mathrm{in}\ \ \mathbb{T}^{3}\times \mathbb{S}^{2}.\\
\end{split}\right.
\end{equation} 
According to Lemma \ref{main result}, we have
\begin{equation}
\begin{aligned}
\|f_{r}\|_{L_{T}^{\infty}L_{\vec w}^{\infty}L^{\infty}_{\vec x}}\leq& C\|(\theta^{N}+\theta_{r})^{4}-(\theta^{N})^{4}-\mathcal{L}_{1}(f^{N}, \theta^{N})\|_{L_{T}^{\infty}L_{\vec w}^{\infty}L^{\infty}_{\vec x}}
\\\leq& C(\|\theta_{r}\|_{L_{T}^{\infty}L^{\infty}_{\vec x}}+\epsilon^{N+1})
\\\leq& C(\|\theta_{r}\|_{L_{T}^{\infty}H^{2}_{\vec x}}+\epsilon^{N+1})
\\\leq& C\epsilon^{5}.
\end{aligned}
\end{equation}
So, 
\begin{equation}
\|f^{\epsilon}-\sum_{k=0}^{9}(f_{k}+f_{I, k})\|_{L_{T}^{\infty}L_{\vec w}^{\infty}L^{\infty}_{\vec x}}
\\\leq C\epsilon^{5}
\end{equation}
and
\begin{equation}
\|\theta^{\epsilon}-\sum_{k=0}^{9}(\theta_{k}+\theta_{I, k})\|_{L_{T}^{\infty}H^{2}_{\vec x}}
\\\leq C\epsilon^{5}.
\end{equation}
Therefore, we derive
\begin{equation}
\|f^{\epsilon}-\theta_{0}^{4}-f_{I, 0}\|_{L_{T}^{\infty}L_{\vec w}^{\infty}L^{\infty}_{\vec x}}
\\\leq C\epsilon
\end{equation}
and
\begin{equation}
\|\theta^{\epsilon}-\theta_{0}-\theta_{I, 0}\|_{L_{T}^{\infty}H^{2}_{\vec x}}
\\\leq C\epsilon.
\end{equation}

\vskip 5mm

\noindent{\bf Acknowledgments}
\noindent{This research was partially supported by the Scientic Research Project of Fuyang Normal University(No. 2024KYQD0106).}

\vskip 2mm
\noindent{\bf Conflict of Interest} The authors declare no conflict of interest.

\end{CJK}
\end{document}